\documentclass[11pt]{article}
\usepackage[letterpaper, portrait, margin=1in]{geometry}

\usepackage{graphicx}
\usepackage{xcolor}
\usepackage{float}

\usepackage[round]{natbib}

\usepackage{amsmath}
\usepackage{amssymb}
\usepackage{amsfonts}
\usepackage{amsthm}

\usepackage{csquotes}

\allowdisplaybreaks[1]

\newtheorem{theorem}{Theorem}[section]

\newtheorem{lemma}[theorem]{Lemma}
\newtheorem{remark}{Remark}[section]
\newtheorem{example}{Example}[section]

\newtheorem{manualtheoreminner}{Theorem}
\newenvironment{manualtheorem}[1]{%
  \renewcommand\themanualtheoreminner{#1}%
  \manualtheoreminner
}{\endmanualtheoreminner}

\newtheorem{manuallemmainner}{Lemma}
\newenvironment{manuallemma}[1]{%
  \renewcommand\themanuallemmainner{#1}%
  \manuallemmainner
}{\endmanuallemmainner}

\newcommand{\ignore}[1]{}

\def\defi{:=}
\def\diff{\mathop{}\!\mathrm{d}}

\title{{\bf Tractable Profit Maximization over Multiple Attributes under Discrete Choice Models}}

\author{
{\bf Anton J. Kleywegt},  {\bf Hongzhang Shao} \\
School of Industrial and Systems Engineering \\
Georgia Institute of Technology \\
Atlanta, Georgia 30332-0205 \\
}

\date{
First Version: July 2020\\
Current Version: December 2021
}

\begin{document}

\maketitle

\begin{abstract}
A fundamental problem in revenue management is to optimally choose the attributes of products, such that the total profit or revenue or market share is maximized. 
Usually, these attributes can affect both a product's market share (probability to be chosen) and its profit margin. 
For example, if a smart phone has a better battery, then it is more costly to be produced, but is more likely to be purchased by a customer. 
The decision maker then needs to choose an optimal vector of attributes for each product that balances this trade-off. 
In spite of the importance of such problems, there is not yet a method to solve it efficiently in general. 
Past literature in revenue management and discrete choice models focus on pricing problems, where price is the only attribute to be chosen for each product. 
Existing approaches to solve pricing problems tractably cannot be generalized to the optimization problem with multiple product attributes as decision variables. 
On the other hand, papers studying product line design with multiple attributes all result in intractable optimization problems. 
Then we found a way to reformulate the static multi-attribute optimization problem, as well as the multi-stage fluid optimization problem with both resource constraints and upper and lower bounds of attributes, as a \emph{tractable} convex conic optimization problem. 
Our result applies to optimization problems under the multinomial logit (MNL) model, the Markov chain (MC) choice model, and with certain conditions, the nested logit (NL) model. 
\end{abstract}

\newpage

%------------------------------------------------------------------------
%	Introduction
%------------------------------------------------------------------------

\section{Introduction}
\label{sec:introduction}

Consider an assortment with multiple products. 
A classical optimization problem in revenue management is to find the optimal prices of these products to be provided to customers. 
In general, when the price of a product increases, its market share shrinks, while its profit margin becomes larger. 
The goal of price optimization is to find a balance between market shares and profit margins, such that the overall profit or revenue is maximized.

In many applications, however, price is not the only variable that affects both market share and profit margin. 
Here we list some common scenarios, under which there are multiple attributes need to be determined, in addition to price:
\begin{itemize}
  \item A ride-sharing service operator wants customers to walk to a pick-up location and to wait for some time before being picked up. With more waiting time and walking distance, the operator can make better matchings, and reduce the operational cost. However, customers are less likely to accept rides if waiting times and walking distances are long.
  \item A health insurance company needs to design both the price and the coverage of a plan. A insurance plan that covers more health care providers is more attractive to customers, but is also more costly to the company.
  \item An airline company sets (1) price, (2) cancelation fee, (3) baggage allowance, and (4) frequent flyer credits for each type of airline ticket (fare class). Lower fees and better services can attract more customers, but reduce the profit margin.
  \item A hotel decides (1) price, (2) cancellation policy, (2) earliest check-in time, and (4) latest check-out time for its rooms. Stricter cancelation policy and check-in / check-out time hurt customer satisfaction (and thus total sales), but also reduce the operational cost.
  \item An e-commerce platform chooses a combination of (1) price, (2) delivery time, and (3) return policy for each product. By offering better delivery time and return policy, the platform is able to charge a higher price without losing sales.
\end{itemize}

Literature that characterize products as a bundle of multiple attributes usually focus on the product line design problem (e.g. \cite{bertsimas2017robust} and \cite{akccakucs2021exact}).
 The objective of a product line design problem is usually to maximize the market share, but papers like \cite{akccakucs2021exact} extend the discussion to profit maximization as well. 
 Under the profit maximization scenario, the attributes also play a role in the profit margin. 
 Papers long this line assumed product attributes to be discrete (usually binary). 
 The solution method presented in these works are usually mixed-integer approaches, which is intractable.

In this paper, we present a \emph{tractable} approach to solve the profit (or revenue) maximization problem, where for each product, there are multiple attributes that need to be determined simultaneously. 
We formulate both the static optimization problem \ref{eqn:static general} and the multi-stage fluid problem \ref{eqn:fluid general}, and consider both resource constraints and upper and lower bounds of attributes. 
More concretely, we discuss problems \ref{eqn:static general} and \ref{eqn:fluid general} under the multinomial logit (MNL) model, the Markov chain (MC) choice model, and the nested logit (NL) model. 
We show that under these choice models, \ref{eqn:static general} and \ref{eqn:fluid general} can be solved by solving a convex conic program (for the NL model, additional conditions are needed to establish the results). 
We also give the necessary and sufficient conditions that an optimal solution exists for each of the problems. 
The formulations, results and discussions are provided in Section \ref{sec:formulation}.

Our work is built upon existing literature on pricing problems under discrete choice models. 
Papers on pricing problems consider price as the only attribute for each product. 
Thus, our formulation of \ref{eqn:static general} and \ref{eqn:fluid general} is a generalization to the pricing problems discussed in these papers. 
Existing methods to solve the pricing problem cannot be applied to \ref{eqn:static general} and \ref{eqn:fluid general}: 
First, the solution approaches in \cite{song2007demand}, \cite{dong2009dynamic}, \cite{li2011pricing}, \cite{keller2014efficient} and \cite{dong2019pricing} relies on the ability of writing prices as inverse functions of market shares (the expected overall profit is concave in market shares). 
However, when each product has multiple attributes instead of price alone, it is no longer possible to replace the attributes by close-form functions of the market shares. 
Second, the optimality condition approach in \cite{wang2012joint}, \cite{gallego2014multiproduct} and \cite{zhang2018multiproduct} only applies to unconstrained optimization problems. 
However, having constrains in \ref{eqn:static general} and \ref{eqn:fluid general} is not only important but also necessary. 
When a product have multiple attributes that can control both the profit margin and the choice probability, it is very likely that we can "scarifies" one attribute for another. 
That is, we may be able to increase one attribute and decrease another attribute at the same time, such that the choice probability of the product remain unchanged, while the profit margin increases. 
If there is no constraint in attributes, the problem is certainly unbounded. 
Third, although discretizing attributes (like in \cite{davis2013assortment} and \cite{gallego2014constrained}) may be a theoretically valid approach, the number of candidate attribute vector grows exponentially with the number of attributes. 
As we embed more attributes into the problem, such a method would easily be impractical. A review of these papers are provided in Section \ref{sec:literature}.

Most results in the papers above can be reproduced with our method. Thus, besides solving the generalized optimization problem \ref{eqn:static general} and \ref{eqn:fluid general}, our method also provides a more unified approach to solve the price (single attribute) optimization problems. We discuss these results as well as some potential extension directions in Section \ref{sec:conclusion}.

%------------------------------------------------------------------------
%	Introduction
%------------------------------------------------------------------------

\section{Related Literature}
\label{sec:literature}

Multi-product price optimization problems under discrete choice models have been the subject of active research in the last two decades. 
Comparing to the past literature, These research use choice models to capture the price-driven substitution behavior from customers. 
Pricing problems under discrete choice models are usually not easy to solve in its original form. 
\cite{hanson1996optimizing} study this problem of maximizing the expected revenue under the multinomial logit (MNL) model, and observe that the expected profit is not concave in price. 
(They then proposed a path-following approach to perturb the objective function into a concave one.) 
Since then, many efficient solution methods have been introduced. We give a review to some of the most important works here:

\emph{\underline{Market Share Reformulation}:} An important approach to solve the pricing problems efficiently is to use market shares as decision variables. 
\cite{song2007demand} and \cite{dong2009dynamic} considered the same problem as in \cite{hanson1996optimizing}, and showed that the expected profit is concave in the market shares. 
In order words, the pricing problem under the MNL model can be solved efficiently by using the market shares of the products as decision variables. 
(Both works assumed that the price sensitivity parameters are constant for all products.) 
\cite{li2011pricing} extended the concavity results to the nested logit (NL) model (which subsumes the MNL model), and showed that the expected profit is concave in the market shares if (1) the price-sensitivity parameters are identical for all the products within a nest, and (2) the nest coefficients are restricted to be in the unit interval. 
As a special case, the expected profit under a general MNL model (with asymmetric price-sensitivity parameters) is concave in the market shares. 
\cite{keller2014efficient} discussed similar concavity results under the general attraction demand model (which also subsumes the MNL model), and showed that constraints such as price bounds and joint price constraints can be added to the pricing problem as linear constraints in the market shares. 
\cite{zhang2018multiproduct} discussed the multi-product pricing problem under the generalized extreme value models (which subsumes the NL model). 
They showed that with homogeneous price sensitivity parameters, the problem can formulation as a convex program. 
Recently, \cite{dong2019pricing} formulate the pricing problem under the Markov chain (MC) choice model. 
They showed that the pricing problem can be solved as a dynamic program. They also presented a market-share based approach to solve the dynamic pricing problem with a single resource.

\emph{\underline{First-order Conditions}:} Another approach to solve the problems is to study the first-order condition. 
(Results that are built upon first-order conditions usually only apply to unconstrained optimization problems.) 
\cite{anderson1992multiproduct}, \cite{besanko1998logit} and \cite{aydin2000product} observex that under a MNL model with homogeneous price sensitivity parameters, the profit margin (price minus cost) is constant across all the products at the optimality of the expected profit. 
\cite{aydin2008joint} and \cite{akcay2010joint} then pointed out that under such scenarios, the profit function is uni-modal with respect to the markup, and its unique optimal solution can be found by solving the first-order conditions. 
\cite{gallego2014multiproduct} extended this result and showed that under the general NL model, the adjusted markup (price minus cost minus the reciprocal of price sensitivity) is constant for all the products within a nest at optimality. 
They then defined an adjusted nest-level markup and showed it to be invariant for all nests, which reduces the pricing problem to a single dimensional optimization problem. 
This one-dimensional problem has a single local maxima when either (1) the difference between price sensitivity parameters are bounded by a certain value, or (2) the dissimilarity parameter is greater than 1. 
These conditions are more general than \cite{li2011pricing}. 
(\cite{wang2012joint} also showed the constant adjusted markup property for pricing problem under the MNL model.) 
Huh and Li (2015) further extended this result to the pricing problem under the multistage nested attraction model. 
They showed that the problem can be reduced to single dimensional, and the optimal solution is unique under generalized conditions. 
\cite{zhang2018multiproduct} discussed the multi-product pricing problem under the generalized extreme value models with homogeneous price sensitivity parameters, and provided similar results. 

\emph{\underline{Discretization of Price}:} A third approach is to limit price candidates to a small, finite set. 
Such methods are usually very similar to approaches used for solving assortment planning. 
For instance, \cite{davis2013assortment} studied assortment planning problems under the MNL model with totally uni-modular constraints. 
They showed that such problems can be solved efficiently as linear programs. 
As a consequence, pricing problems under the MNL model with finite possible prices can also be formulated as a linear program. 
\cite{gallego2014constrained} followed a similar approach and extend the results to pricing problems under the nested logit model. 
Researchers also used this approach to model constraints in the pricing problem in a tractable way. 
As an example, \cite{davis2017pricing} show that under the NL model, quality consistency constraints (pairwise inequalities in prices) can be added to the pricing problem, and the problem can be solved as a linear program.

Research on discrete choice models dates back to \cite{mcfadden1973conditional}. The work integrates the random utility framework from \cite{thurstone1927law} and the choice axioms introduced by \cite{luce1959individual}, and provides a general econometric procedure to the analysis of individual choice behavior. We point the readers to \cite{ben1985discrete} and \cite{train2009discrete} for more in-depth discussions on the (parameterized) random utility choice models, such as the MNL model and the NL model mentioned above. Our paper also covers the Markov chain (MC) choice model. The MC model was recently analyzed in \cite{blanchet2016markov}. A discussion on different choice-modeling approaches and their mathematical equivalence can be found in \cite{feng2017relation}

%------------------------------------------------------------------------
%	Section
%------------------------------------------------------------------------

\section{Optimization of Product Attributes}
\label{sec:formulation}

A seller offers an assortment of $J$ products, indexed by $j \in \mathcal{J}$.
Each product has a number of features, such as price and service characteristics, the values of which can be chosen by the seller.
Different products may have different numbers of these features.
In this paper, \textit{attribute} refers to a feature for which the seller has to choose values for one or more products, and that has unique parameter values, as specified below.
The value of a feature for a product affects the seller's profit in two ways:
(1)~The feature's value affects the demand for the product (and the demand for other products).
(2)~The feature's value affects the seller's revenue and/or cost, and thus profit margin, per unit product sold.
Both of these effects are specified by a model with parameters that may depend on the product and the feature.
For example, in a discrete choice model, the systematic utility associated with a product~$j$ and feature~$\ell$ may be given by $\beta_{j\ell} y_{j\ell}$, where $\beta_{j\ell}$ is a parameter and $y_{j\ell}$ is the value of feature~$\ell$ for product~$j$.
Similarly, the profit margin associated with a product~$j$ and feature~$\ell$ may be given by $\check{\phi}_{j\ell} y_{j\ell}$, where $\check{\phi}_{j\ell}$ is a parameter.
We will simplify the notation by defining an attribute to be a feature with unique parameter values.
Thus, if two products~$j$ and~$j'$ have different parameter values associated with a feature~$\ell$, for example $\beta_{j\ell} \neq \beta_{j'\ell}$ or $\check{\phi}_{j\ell} \neq \check{\phi}_{j'\ell}$, then there are separate attributes~$k$ and~$k'$ in the model for feature~$\ell$ when used in combination with product~$j$, and for feature~$\ell$ when used in combination with product~$j'$.
Also, if two products have the same parameter values associated with a feature, then the two products have the same attribute, but they may have different values for the attribute.
That way, parameter values are uniquely specified by the attribute.
Let $\mathcal{K}$ denote the resulting index set of attributes, and for each $j \in \mathcal{J}$, let $\mathcal{K}_{j} \subset \mathcal{K}$ denote the index set of the attributes for product~$j$ for which the seller has to choose a value.
Thus, $\mathcal{K}_{j} \neq \varnothing$ ($K_{j} \defi |\mathcal{K}_{j}| \geq 1$) for each $j \in \mathcal{J}$, $\mathcal{K} = \cup_{j \in \mathcal{J}} \mathcal{K}_{j}$, and $\mathcal{K}_{j} \cap \mathcal{K}_{j'}$ may be nonempty for two products~$j$ and~$j'$.

For each $j \in \mathcal{J}$ and $k \in \mathcal{K}_{j}$, the value chosen for attribute~$k$ for product~$j$ is denoted with~$y_{jk}$.
Let $y \defi (y_{jk}, \, j \in \mathcal{J}, \, k \in \mathcal{K}_{j})$ denote the list of attribute values, and let $\check{P}_{j}(y)$ denote the probability that a customer chooses product~$j$.
We assume that the profit margin per unit of product~$j$ sold depends on $y_{jk}, k \in \mathcal{K}_{j}$, and that this dependence is specified by an affine function $\sum_{k \in \mathcal{K}_{j}} \check{\phi}_{k} y_{jk} - \check{\psi}_{j}$, where $\check{\phi}_{k}, k \in \mathcal{K}$ and $\check{\psi}_{j}, j \in \mathcal{J}$ are parameters.
For example, if attribute~$k$ represents product price, then $\check{\phi}_{k} = 1$.
In general, $\sum_{k \in \mathcal{K}_{j}} \check{\phi}_{k} y_{jk} - \check{\psi}_{j}$ can be regarded as an affine approximation of the relation between the attribute values of product~$j$ and the profit margin per unit of product~$j$ sold.

%Attribute values may be unbounded below (in which case the constraint $y_{jk} \ge \underline{y}_{jk}$ is dropped) or unbounded above (in which case the constraint $y_{jk} \le \overline{y}_{jk}$ is dropped).
It is important to make provision for bounds on attribute values, for the following reasons:
\begin{enumerate}
\item
Usually discrete choice models and profit margin models are calibrated with a limited range of attribute value data, and it would be unwise to allow selection of attribute values much outside this range.
\item
As explained in Remark~\ref{rem:attribute bounds} below, if multiple attributes are unbounded, then one attribute can be traded off for another attribute to obtain unbounded profit.
\end{enumerate}
Let $\underline{\mathcal{K}}_{j}$ denote the set of attributes $k \in \mathcal{K}_{j}$ such that $y_{jk}$ has a lower bound $\underline{y}_{jk}$, and let $\overline{\mathcal{K}}_{j}$ denote the set of attributes $k \in \mathcal{K}_{j}$ such that $y_{jk}$ has an upper bound $\overline{y}_{jk}$.
It is assumed that $\underline{y}_{jk} \le \overline{y}_{jk}$ for all $j \in \mathcal{J}$ and $k \in \underline{\mathcal{K}}_{j} \cap \overline{\mathcal{K}}_{j}$.

The resulting static attribute optimization problem is
\begin{align}
\max_{y} \quad & \sum_{j \in \mathcal{J}} \left(\sum_{k \in \mathcal{K}_{j}} \check{\phi}_{k} y_{jk} - \check{\psi}_{j}\right) \check{P}_{j}(y)
\tag{\textsf{SP}}
\label{eqn:static general} \\
\text{s.t.} \quad
& y_{jk} \ \ \geq \ \ \underline{y}_{jk}
& \forall \ j \in \mathcal{J}, \ k \in \underline{\mathcal{K}}_{j} \nonumber \\
& y_{jk} \ \ \leq \ \ \overline{y}_{jk}
& \forall \ j \in \mathcal{J}, \ k \in \overline{\mathcal{K}}_{j} \nonumber
\end{align}

We will also consider the following revenue management problem.
Products require the use of resources.
Let $\mathcal{R}$ denote the set of resources, and for each product $j \in \mathcal{J}$ and each resource $r \in \mathcal{R}$, let $a_{rj}$ denote the amount of resource~$r$ needed per unit of product~$j$.
For each resource $r \in \mathcal{R}$, there are $b_{r}$ units of the resource available for use over a selling horizon.
The selling horizon is discretized into time periods indexed by $t = 0,1,\ldots,T$.
In each time period~$t$, the expected number of customer arrivals is denoted by $\lambda_{t}$.
Parameters and decisions may depend on $t$, which is reflected in the notation.
We will consider the following fluid optimization problem:
\begin{align}
\max_{y} \quad & \sum_{t = 0}^{T} \lambda_{t} \sum_{j \in \mathcal{J}} \left(\sum_{k \in \mathcal{K}_{jt}} \check{\phi}_{kt} y_{jkt} - \check{\psi}_{jt}\right) \check{P}_{jt}(y)
\tag{\textsf{FP}}
\label{eqn:fluid general} \\
\text{s.t.} \quad
& \sum_{t = 0}^{T} \lambda_{t} \sum_{j \in \mathcal{J}} a_{rj} \check{P}_{jt}(y) \ \ \le \ \ b_{r}
& \forall \ r \in \mathcal{R} \nonumber \\
& y_{jkt} \ \ \geq \ \ \underline{y}_{jkt}
& \forall \ j \in \mathcal{J}, \ k \in \underline{\mathcal{K}}_{jt}, \ t = 0,1,\ldots,T \nonumber \\
& y_{jkt} \ \ \leq \ \ \overline{y}_{jkt}
& \forall \ j \in \mathcal{J}, \ k \in \overline{\mathcal{K}}_{jt}, \ t = 0,1,\ldots,T \nonumber
\end{align}

In the remainder of the section, we consider problems~\ref{eqn:static general} and~\ref{eqn:fluid general} for various widely used discrete choice models.

%------------------------------------------------------------------------

\subsection{\ref{eqn:static general} under the MNL Model}

First we consider problem~\ref{eqn:static general} in which $\check{P}_{j}(y)$ is given by a multinomial logit model:
\[
\check{P}_{j}(y) \ \ = \ \ \frac{\exp\left(\alpha_{j} - \sum_{k \in \mathcal{K}_{j}} \beta_{k} y_{jk}\right)}{1 + \sum_{j' \in \mathcal{J}} \exp\left(\alpha_{j'} - \sum_{k \in \mathcal{K}_{j'}} \beta_{k} y_{j'k}\right)}
\]
where $\alpha_{j}, j \in \mathcal{J}$ are the product ``baseline attractiveness'' parameters, and $\beta_{k}, k \in \mathcal{K}$ are the attribute ``sensitivity'' parameters.

Without loss of generality, we assume that $\beta_{k} > 0$ (because the sign of each $y_{jk}$ can be chosen such that $\beta_{k} > 0$).
Therefore we also assume that $\check{\phi}_{k} > 0$, that is, there is a trade-off between the effect of attribute~$k$ on the demand for products and the effect of attribute~$k$ on the profit margin of products --- as $y_{jk}$ increases the demand for product~$j$ decreases and the profit margin per unit of product~$j$ increases.
To further simplify notation, consider the scaled attributes $x_{jk} \defi \beta_{k} y_{jk} - \alpha_{j} / K_{j}$, the scaled attribute lower bounds $\underline{x}_{jk} \defi \beta_{k} \underline{y}_{jk} - \alpha_{j} / K_{j}$, the scaled attribute upper bounds $\overline{x}_{jk} \defi \beta_{k} \overline{y}_{jk} - \alpha_{j} / K_{j}$, the scaled profit margins $\phi_{k} \defi \beta_{k} \check{\phi}_{k} > 0$, the scaled item fixed costs $\psi_{j} \defi \check{\psi}_{j} + (\alpha_{j} / K_{j}) \sum_{k \in \mathcal{K}_{j}} \check{\phi}_{k}$, for every $j \in \mathcal{J}$ and $k \in \mathcal{K}_{j}$.
Let $x \defi (x_{jk}, \, j \in \mathcal{J}, \, k \in \mathcal{K}_{j})$.
Then the choice probabilities are given by
\[
P_{j}(x) \ \ = \ \ \frac{\exp\left(- \sum_{k \in \mathcal{K}_{j}} x_{jk}\right)}{1 + \sum_{j' \in \mathcal{J}} \exp\left(- \sum_{k \in \mathcal{K}_{j'}} x_{j'k}\right)}
\]
and the expected profit of product~$j$ becomes $\left(\sum_{k \in \mathcal{K}_{j}} \phi_{k} x_{jk} - \psi_{j}\right) P_{j}(x)$.
Thus, we consider the static attribute optimization problem
\begin{align}
\max_{d, \, x} \quad & \sum_{j \in \mathcal{J}} \left(\sum_{k \in \mathcal{K}_{j}} \phi_{k} x_{jk} - \psi_{j}\right) d_{j}
\tag{$\mathsf{SP^{MNL}_{1}}$}
\label{eqn:static MNL1} \\
\text{s.t.} \quad
& d_{j} \ \ = \ \ \frac{\exp\left(- \sum_{k \in \mathcal{K}_{j}} x_{jk}\right)}{1 + \sum_{j' \in \mathcal{J}} \exp\left(- \sum_{k \in \mathcal{K}_{j'}} x_{j'k}\right)}
& \forall \ j \in \mathcal{J} \label{eqn:spmnl1 dj} \\
& x_{jk} \ \ \geq \ \ \underline{x}_{jk}
& \forall \ j \in \mathcal{J}, \ k \in \underline{\mathcal{K}}_{j} \label{eqn:spmnl1-lower-bound} \\
& x_{jk} \ \ \leq \ \ \overline{x}_{jk}
& \forall \ j \in \mathcal{J}, \ k \in \overline{\mathcal{K}}_{j} \label{eqn:spmnl1-upper-bound}
\end{align}
In \ref{eqn:static MNL1}, $d \defi (d_{j}, \, j \in \mathcal{J})$ denotes the vector of market shares.

\begin{remark}
\label{rem:attribute bounds}
Bounds on attribute values are important when multiple attributes are chosen for a product.
Suppose that $\phi_{k} > \phi_{k'}$ and that $k, k' \in \mathcal{K}_{j}$ for some $j \in \mathcal{J}$.
If attribute~$k$ is not bounded above and attribute~$k'$ is not bounded below, then $x_{jk}$ can be increased and $x_{jk'}$ can be decreased by the same amount, while keeping all other attribute values constant.
Thereby all demands remain constant, while the profit contribution of product~$j$ grows without bound.
In an application, it makes sense that attributes can be traded off to some extent, but arbitrarily large or small attribute values are not sensible.
\end{remark}

Problem~\ref{eqn:static MNL1} is not a convex optimization problem.
One of the reasons is that nonlinear constraint~\eqref{eqn:spmnl1 dj} is an equality constraint.
Also, if some products have multiple attributes ($K_{j} > 1$ for some~$j$), then the approach in \cite{song2007demand}, \cite{dong2009dynamic}, \cite{li2011pricing}, \cite{keller2014efficient}, and \cite{dong2019pricing} of inverting the constraint to write~$x$ as a function of~$d$ does not apply here.
We will relax the equality constraints to inequality constraints, show that an optimal solution of the relaxation provides an optimal solution of~\ref{eqn:static MNL1}, then reformulate the relaxation as a convex optimization problem, and then provide necessary and sufficient conditions for this convex optimization problem to have an optimal solution.

First, note that a decision variable $d_{0}$ and the constraints $d_{0} + \sum_{j \in \mathcal{J}} d_{j} = 1$ and $d_{0} > 0$ can be added to problem~\ref{eqn:static MNL1} without changing the set of feasible $(d, x)$-values or the objective value for any feasible $(d, x)$-value, because the constraints will force
\begin{equation}
\label{eqn:d0}
d_{0} \ \ = \ \ \frac{1}{1 + \sum_{j \in \mathcal{J}} \exp\left(- \sum_{k \in \mathcal{K}_{j}} x_{jk}\right)}
\end{equation}
without affecting the objective value.
Next, we rewrite constraints~\eqref{eqn:spmnl1 dj} to obtain the following static attribute optimization problem
\begin{align}
\max_{d, \, d_{0}, \, x} \quad & \sum_{j \in \mathcal{J}} \left(\sum_{k \in \mathcal{K}_{j}} \phi_{k} x_{jk} - \psi_{j}\right) d_{j}
\tag{$\mathsf{SP^{MNL}_{2}}$}
\label{eqn:static MNL2} \\
\text{s.t.} \quad
& \ln\left(\frac{d_{j}}{d_{0}}\right) \ \ = \ \ - \sum_{k \in \mathcal{K}_{j}} x_{jk}
& \forall \ j \in \mathcal{J} \label{mp2-mnl} \\
& x_{jk} \ \ \geq \ \ \underline{x}_{jk}
& \forall \ j \in \mathcal{J}, \ k \in \underline{\mathcal{K}}_{j} \label{eqn:spmnl2-lower-bound} \\
& x_{jk} \ \ \leq \ \ \overline{x}_{jk}
& \forall \ j \in \mathcal{J}, \ k \in \overline{\mathcal{K}}_{j} \label{eqn:spmnl2-upper-bound} \\
& d_{0} + \sum_{j \in \mathcal{J}} d_{j} \ \ = \ \ 1 \label{mp2-totalmarket} \\
& d \ > \ 0, \quad d_{0} \ > \ 0 \nonumber
\end{align}
Note that for any $(d, x)$ feasible for~\ref{eqn:static MNL1}, it holds that $(d, d_{0}, x)$ with $d_{0}$ given by~\eqref{eqn:d0} is feasible for~\ref{eqn:static MNL2} and has the same objective value.
Conversely, consider any $(d, d_{0}, x)$ feasible for~\ref{eqn:static MNL2}.
It follows from~\eqref{mp2-mnl} that $d_{j} = d_{0} \exp\left(- \sum_{k \in \mathcal{K}_{j}} x_{jk}\right)$, then it follows from~\eqref{mp2-totalmarket} that $d_{0}$ satisfies~\eqref{eqn:d0}, and finally it follows that $d_{j}$ satisfies~\eqref{eqn:spmnl1 dj}.
Thus, for any $(d, d_{0}, x)$ feasible for~\ref{eqn:static MNL2}, it holds that $(d, x)$ is feasible for~\ref{eqn:static MNL1} and has the same objective value.

Note that constraint~\eqref{mp2-mnl} is equivalent to $- \sum_{k \in \mathcal{K}_{j}} x_{jk} \le \ln\left(d_{j} / d_{0}\right) \le - \sum_{k \in \mathcal{K}_{j}} x_{jk}$.
Next, we relax constraint~\eqref{mp2-mnl} as follows.
Let $\overline{\mathcal{J}} \defi \left\{j \in \mathcal{J} \, : \, \overline{\mathcal{K}}_{j} = \mathcal{K}_{j}\right\}$ denote the set of products for which all attributes are upper bounded.
For each $j \in \overline{\mathcal{J}}$, constraint~\eqref{mp2-mnl} is relaxed to $- \sum_{k \in \mathcal{K}_{j}} \overline{x}_{jk} \le \ln\left(d_{j} / d_{0}\right) \le - \sum_{k \in \mathcal{K}_{j}} x_{jk}$, and for each $j \in \mathcal{J} \setminus \overline{\mathcal{J}}$, constraint~\eqref{mp2-mnl} is relaxed to $\ln\left(d_{j} / d_{0}\right) \le - \sum_{k \in \mathcal{K}_{j}} x_{jk}$.
%(If $x_{jk}$ is unbounded above for some $k \in \mathcal{K}_{j}$, then the constraint on the left side is omitted for~$j$.)
Thus we consider the following relaxation of~\ref{eqn:static MNL2}:
\begin{align}
\max_{d, \, d_{0}, \, x} \quad & \sum_{j \in \mathcal{J}} \left(\sum_{k \in \mathcal{K}_{j}} \phi_{k} x_{jk} - \psi_{j}\right) d_{j}
\tag{$\mathsf{SP^{MNL}_{3}}$}
\label{eqn:static MNL3} \\
\text{s.t.} \quad
& \ln\left(\frac{d_{j}}{d_{0}}\right) \ \ \le \ \ - \sum_{k \in \mathcal{K}_{j}} x_{jk}
& \forall \ j \in \mathcal{J} \label{mp3-mnl} \\
& \ln\left(\frac{d_{0}}{d_{j}}\right) \ \ \le \ \ \sum_{k \in \mathcal{K}_{j}} \overline{x}_{jk}
& \forall \ j \in \overline{\mathcal{J}} \label{mp3-total upper bound} \\
& x_{jk} \ \ \geq \ \ \underline{x}_{jk}
& \forall \ j \in \mathcal{J}, \ k \in \underline{\mathcal{K}}_{j} \label{eqn:spmnl3-lower-bound} \\
& x_{jk} \ \ \leq \ \ \overline{x}_{jk}
& \forall \ j \in \mathcal{J}, \ k \in \overline{\mathcal{K}}_{j} \label{eqn:spmnl3-upper-bound} \\
& d_{0} + \sum_{j \in \mathcal{J}} d_{j} \ \ = \ \ 1 \label{mp3-totalmarket} \\
& d \ > \ 0, \quad d_{0} \ > \ 0 \nonumber
\end{align}
Next we show that an optimal solution of the relaxation gives an optimal solution of~\ref{eqn:static MNL2}.
Consider any $(d, d_{0}, x)$ feasible for~\ref{eqn:static MNL3}.
If $\ln\left(d_{j} / d_{0}\right) = - \sum_{k \in \mathcal{K}_{j}} x_{jk}$ for all~$j$, then $(d, d_{0}, x)$ is feasible for~\ref{eqn:static MNL2}, and has the same objective value.
Otherwise, $\ln\left(d_{j} / d_{0}\right) < - \sum_{k \in \mathcal{K}_{j}} x_{jk}$ for some~$j$.
Then consider 2~cases: (1)~$j \in \mathcal{J} \setminus \overline{\mathcal{J}}$ or $x_{jk} < \overline{x}_{jk}$ for some $k \in \overline{\mathcal{K}}_{j}$, or (2)~$j \in \overline{\mathcal{J}}$ and $x_{jk} = \overline{x}_{jk}$ for all $k \in \mathcal{K}_{j}$.
In case~(1), $x_{jk}$ can be increased in~\ref{eqn:static MNL3} while keeping $(d,d_{0})$ unchanged.
Since $\phi_{k} > 0$, the objective value will improve, and thus such $(d, d_{0}, x)$ cannot be optimal for~\ref{eqn:static MNL3}.
In case~(2), it follows that $\ln\left(d_{j} / d_{0}\right) < - \sum_{k \in \mathcal{K}_{j}} \overline{x}_{jk}$, that is, $\ln\left(d_{0} / d_{j}\right) > \sum_{k \in \mathcal{K}_{j}} \overline{x}_{jk}$, which violates constraint~\eqref{mp3-total upper bound}.
Thus, an optimal solution for~\ref{eqn:static MNL3} satisfies $\ln\left(d_{j} / d_{0}\right) = - \sum_{k \in \mathcal{K}_{j}} x_{jk}$ for all~$j$, and hence is feasible and therefore optimal for~\ref{eqn:static MNL2}.

Problem~\ref{eqn:static MNL3} can be reformulated as a convex optimization problem as follows.
Note that the functions $(d_{j},d_{0}) \mapsto d_{j} \ln\left(d_{j} / d_{0}\right)$ and $(d_{j},d_{0}) \mapsto d_{0} \ln\left(d_{0} / d_{j}\right)$ are convex on $(0,\infty)^2$.
Also, introduce a new variable $u_{jk} = d_{j} x_{jk}$, and let $u \defi (u_{jk}, \, j \in \mathcal{J}, \, k \in \mathcal{K}_{j})$.
Thus we consider the following convex reformulation of~\ref{eqn:static MNL3}:
\begin{align}
\max_{d, \, d_{0}, \, u} \quad & \sum_{j \in \mathcal{J}} \left(\sum_{k \in \mathcal{K}_{j}} \phi_{k} u_{jk} - \psi_{j} d_{j}\right)
\tag{$\mathsf{SP^{MNL}_{4}}$}
\label{eqn:static MNL4} \\
\text{s.t.} \quad
& d_{j} \ln\left(\frac{d_{j}}{d_{0}}\right) \ \ \le \ \ - \sum_{k \in \mathcal{K}_{j}} u_{jk}
& \forall \ j \in \mathcal{J} \label{mp4-mnl} \\
& d_{0} \ln\left(\frac{d_{0}}{d_{j}}\right) \ \ \le \ \ \sum_{k \in \mathcal{K}_{j}} \overline{x}_{jk} d_{0}
& \forall \ j \in \overline{\mathcal{J}} \label{mp4-total upper bound} \\
& u_{jk} \ \ \geq \ \ \underline{x}_{jk} d_{j}
& \forall \ j \in \mathcal{J}, \ k \in \underline{\mathcal{K}}_{j} \label{eqn:spmnl4-lower-bound} \\
& u_{jk} \ \ \leq \ \ \overline{x}_{jk} d_{j}
& \forall \ j \in \mathcal{J}, \ k \in \overline{\mathcal{K}}_{j} \label{eqn:spmnl4-upper-bound} \\
& d_{0} + \sum_{j \in \mathcal{J}} d_{j} \ \ = \ \ 1 \label{mp4-totalmarket} \\
& d \ > \ 0, \quad d_{0} \ > \ 0 \label{eqn:mp4-positive d}
\end{align}
Note that for any $(d, d_{0}, x)$ feasible for~\ref{eqn:static MNL3}, $(d, d_{0}, u)$ with $u_{jk} = d_{j} x_{jk}$ for all $j \in \mathcal{J}$, $k \in \mathcal{K}_{j}$ is feasible for~\ref{eqn:static MNL4}, and has the same objective value.
Conversely, for any $(d, d_{0}, u)$ feasible for~\ref{eqn:static MNL4}, $(d, d_{0}, x)$ with $x_{jk} = u_{jk} / d_{j}$ for all $j \in \mathcal{J}$, $k \in \mathcal{K}_{j}$ is feasible for~\ref{eqn:static MNL3}, and has the same objective value.

The feasible set of~\ref{eqn:static MNL4} may not be closed.
For example, suppose that $j \in \mathcal{J} \setminus \overline{\mathcal{J}}$ and that $\underline{x}_{jk} \le 0$ for all $k \in \underline{\mathcal{K}}_{j}$.
Consider a sequence $\{(d^{n}, d^{n}_{0}, u^{n})\}_{n=0}^{\infty}$ of feasible solutions with $d^{0}_{j} \le d^{0}_{0}$, $d^{n}_{j} = d^{0}_{j} / n$, $d^{n}_{0} = d^{0}_{0} + d^{0}_{j} (1 - 1 / n)$, $u_{jk} = 0$ for all $k \in \mathcal{K}_{j}$, and $(d^{n}_{j'}, u^{n}_{j'})$ remain constant for all $j' \neq j$.
Then $(d^{n}, d^{n}_{0}, u^{n}) \to (\bar{d}, \bar{d}_{0}, \bar{u})$ with $\bar{d}_{j} = 0$, and thus the limit $(\bar{d}, \bar{d}_{0}, \bar{u})$ is not in the feasible set of~\ref{eqn:static MNL4}.

Next we relax~\ref{eqn:static MNL4} to make the feasible set closed, and then we show that optimal solutions are not affected by the relaxation.
%The feasible set is made closed by extending the dimension of the domain, and not by taking the closure of the feasible set.
%For each $j \in \mathcal{J} \setminus \overline{\mathcal{J}}$, there is a free variable~$w_{j}$.
%Let $w \defi (w_{j}, \, j \in \mathcal{J} \setminus \overline{\mathcal{J}})$.
Let
\begin{align*}
\mathcal{K}_{\exp} \ \ & \defi \ \ \mbox{closure}\big\{(a_{1}, a_{2}, a_{3}) \; : \; a_{3} \leq a_{2} \ln(a_{1} / a_{2}), \; a_{1} > 0, \; a_{2} > 0\big\} \\
& = \ \ \big\{(a_{1}, a_{2}, a_{3}) \; : \; a_{3} \leq a_{2} \ln(a_{1} / a_{2}), \; a_{1} > 0, \; a_{2} > 0\big\} \cup \big\{(a_{1}, 0, a_{3}) \; : \; a_{1} \geq 0, \; a_{3} \leq 0\big\}
\end{align*}
denote the exponential cone.
Then consider the following convex conic relaxation of~\ref{eqn:static MNL4}:
\begin{align}
\max_{d, \, d_{0}, \, u} \quad & \sum_{j \in \mathcal{J}} \left(\sum_{k \in \mathcal{K}_{j}} \phi_{k} u_{jk} - \psi_{j} d_{j}\right)
\tag{$\mathsf{SP^{MNL}_{5}}$}
\label{eqn:static MNL5} \\
\text{s.t.} \quad
& \left(d_{0}, d_{j}, \sum_{k \in \mathcal{K}_{j}} u_{jk}\right) \ \ \in \ \ \mathcal{K}_{\exp}
& \forall \ j \in \mathcal{J} \label{mp5-mnl} \\
& \left(d_{j}, d_{0}, - \sum_{k \in \mathcal{K}_{j}} \overline{x}_{jk} d_{0}\right) \ \ \in \ \ \mathcal{K}_{\exp}
& \forall \ j \in \overline{\mathcal{J}} \label{mp5-total upper bound} \\
%& \left(d_{j}, d_{0}, w_{j}\right) \ \ \in \ \ \mathcal{K}_{\exp}
%& \forall \ j \in \mathcal{J} \setminus \overline{\mathcal{J}} \label{mp5-free upper bound} \\
& u_{jk} \ \ \geq \ \ \underline{x}_{jk} d_{j}
& \forall \ j \in \mathcal{J}, \ k \in \underline{\mathcal{K}}_{j} \label{eqn:spmnl5-lower-bound} \\
& u_{jk} \ \ \leq \ \ \overline{x}_{jk} d_{j}
& \forall \ j \in \mathcal{J}, \ k \in \overline{\mathcal{K}}_{j} \label{eqn:spmnl5-upper-bound} \\
& d_{0} + \sum_{j \in \mathcal{J}} d_{j} \ \ = \ \ 1 \label{mp5-totalmarket}
\end{align}
%Consider any $(d, d_{0}, u)$ feasible for~\ref{eqn:static MNL4}.
%For each $j \in \mathcal{J} \setminus \overline{\mathcal{J}}$, set $w_{j} = d_{0} \ln(d_{j} / d_{0})$.
%Then $(d, d_{0}, u, w)$ is feasible for~\ref{eqn:static MNL5} and has the same objective value.
%Next, consider any $(d, d_{0}, u, w)$ feasible for~\ref{eqn:static MNL5}.
%We show by contradiction that $d > 0$, $d_{0} > 0$, and thus $(d, d_{0}, u)$ is feasible for~\ref{eqn:static MNL4}.
%Suppose that $d_{0} = 0$.
%Then it follows from constraint~\eqref{mp5-mnl} that $d_{j} = 0$ for all~$j$, which would violate constraint~\eqref{mp5-totalmarket}.
%Similarly, suppose that $d_{j} = 0$ for some~$j$.
%If $j \in \overline{\mathcal{J}}$, then it follows from constraint~\eqref{mp5-total upper bound} that $d_{0} = 0$, which as shown above cannot be feasible for~\ref{eqn:static MNL5}.
%If $j \in \mathcal{J} \setminus \overline{\mathcal{J}}$, then it follows from constraint~\eqref{mp5-free upper bound} that $d_{0} = 0$, which cannot be feasible for~\ref{eqn:static MNL5}.
%Thus $(d, d_{0}, u)$ is feasible for~\ref{eqn:static MNL4} and has the same objective value.

The following results hold for \ref{eqn:static MNL5}:

\begin{lemma}
\label{lem:spmnl5 optimal}
\ref{eqn:static MNL5} has an optimal solution if and only if for each $j \in \mathcal{J}$, and for all $k_{1} \in \mathcal{K}_{j} \setminus \overline{\mathcal{K}}_{j}$, $k_{2} \in \mathcal{K}_{j} \setminus \underline{\mathcal{K}}_{j}$, it holds that $\phi_{k_{1}} \le \phi_{k_{2}}$; otherwise, \ref{eqn:static MNL5} is unbounded.
\end{lemma}

\begin{lemma}
\label{lem:spmnl4=spmnl5}
Every optimal solution for~\ref{eqn:static MNL5} is also optimal for~\ref{eqn:static MNL4}.
\end{lemma}

Next we summarize the results for~\ref{eqn:static general} under the MNL choice model.

\begin{theorem}
\label{lem:spmnl1 optimal}
\ref{eqn:static MNL1} has an optimal solution if and only if for each $j \in \mathcal{J}$, and for all $k_{1} \in \mathcal{K}_{j} \setminus \overline{\mathcal{K}}_{j}$, $k_{2} \in \mathcal{K}_{j} \setminus \underline{\mathcal{K}}_{j}$, it holds that $\phi_{k_{1}} \le \phi_{k_{2}}$, in which case an optimal solution can be found by solving the convex conic program~\ref{eqn:static MNL5}.
Otherwise, \ref{eqn:static MNL1} is unbounded.
\end{theorem}

\subsection{\ref{eqn:fluid general} under the MNL Model}

Next, we consider the fluid revenue management problem~\ref{eqn:fluid general} under the following MNL model:
\[
\check{P}_{jt}(y) \ \ = \ \ \frac{\exp\left(\alpha_{jt} - \sum_{k \in \mathcal{K}_{jt}} \beta_{kt} y_{jkt}\right)}{1 + \sum_{j' \in \mathcal{J}} \exp\left(\alpha_{j't} - \sum_{k \in \mathcal{K}_{j't}} \beta_{kt} y_{j'kt}\right)}
\]
Let $x_{jkt} \defi \beta_{kt} y_{jkt} - \alpha_{jt} / K_{jt}$, $\underline{x}_{jkt} \defi \beta_{kt} \underline{y}_{jkt} - \alpha_{jt} / K_{jt}$, $\overline{x}_{jkt} \defi \beta_{kt} \overline{y}_{jkt} - \alpha_{jt} / K_{jt}$, $\phi_{kt} \defi \beta_{kt} \check{\phi}_{kt} > 0$, $\psi_{jt} \defi \check{\psi}_{jt} + (\alpha_{jt} / K_{jt}) \sum_{k \in \mathcal{K}_{jt}} \check{\phi}_{kt}$, for every $j \in \mathcal{J}$, $t = 0,1,\ldots,T$, and $k \in \mathcal{K}_{jt}$.
Let $x \defi (x_{jkt}, \, j \in \mathcal{J}, \, k \in \mathcal{K}_{jt}, \, t = 0,1,\ldots,T)$, $d \defi (d_{jt}, \, j \in \mathcal{J}, \, t = 0,1,\ldots,T)$, $d_{0} \defi (d_{0t}, \, t = 0,1,\ldots,T)$, and $u \defi (u_{jkt}, \, j \in \mathcal{J}, \, k \in \mathcal{K}_{jt}, \, t = 0,1,\ldots,T)$.
Thus, we consider the optimization problem
\begin{align}
\max_{d, \, x} \quad & \sum_{t = 0}^{T} \lambda_{t} \sum_{j \in \mathcal{J}} \left(\sum_{k \in \mathcal{K}_{jt}} \phi_{kt} x_{jkt} - \psi_{jt}\right) d_{jt}
\tag{\textsf{$\mathsf{FP_1^{MNL}}$}}
\label{eqn:fluid MNL1} \\
\text{s.t.} \quad
& d_{jt} \ \ = \ \ \frac{\exp\left(- \sum_{k \in \mathcal{K}_{jt}} x_{jkt}\right)}{1 + \sum_{j' \in \mathcal{J}} \exp\left(- \sum_{k \in \mathcal{K}_{j't}} x_{j'kt}\right)}
& \forall \ j \in \mathcal{J}, \ t = 0,1,\ldots,T \nonumber \\
& \sum_{t = 0}^{T} \lambda_{t} \sum_{j \in \mathcal{J}} a_{rj} d_{jt} \ \ \le \ \ b_{r}
& \forall \ r \in \mathcal{R} \nonumber \\
& x_{jkt} \ \ \geq \ \ \underline{x}_{jkt}
& \forall \ j \in \mathcal{J}, \ k \in \underline{\mathcal{K}}_{jt}, \ t = 0,1,\ldots,T \nonumber \\
& x_{jkt} \ \ \leq \ \ \overline{x}_{jkt}
& \forall \ j \in \mathcal{J}, \ k \in \underline{\mathcal{K}}_{jt}, \ t = 0,1,\ldots,T \nonumber
\end{align}

Consider the following convex conic relaxation of~\ref{eqn:fluid MNL1}:
\begin{align}
\max_{d, \, d_{0}, \, u} \quad & \sum_{t = 0}^{T} \lambda_{t} \sum_{j \in \mathcal{J}} \left(\sum_{k \in \mathcal{K}_{jt}} \phi_{kt} u_{jkt} - \psi_{jt} d_{jt}\right)
\tag{\textsf{$\mathsf{FP_2^{MNL}}$}}
\label{eqn:fluid MNL2} \\
\text{s.t.} \quad
& \sum_{t = 0}^{T} \lambda_{t} \sum_{j \in \mathcal{J}} a_{rj} d_{jt} \ \ \le \ \ b_{r}
& \forall \ r \in \mathcal{R}
\label{eqn:fluid MNL2 resource} \\
& \left(d_{0t}, d_{jt}, \sum_{k \in \mathcal{K}_{jt}} u_{jkt}\right) \ \ \in \ \ \mathcal{K}_{\exp}
& \forall \ j \in \mathcal{J}, \ t = 0,1,\ldots,T
\label{eqn:fluid MNL2 cone} \\
& \left(d_{jt}, d_{0t}, - \sum_{k \in \mathcal{K}_{jt}} \overline{x}_{jkt} d_{0t}\right) \ \ \in \ \ \mathcal{K}_{\exp}
& \forall \ j \in \overline{\mathcal{J}}_{t}, \ t = 0,1,\ldots,T
\label{eqn:fluid MNL2 cone dummy} \\
& u_{jkt} \ \ \geq \ \ \underline{x}_{jkt} d_{jt}
& \forall \ j \in \mathcal{J}, \ k \in \underline{\mathcal{K}}_{jt}, \ t = 0,1,\ldots,T
\label{eqn:fluid MNL2 lower} \\
& u_{jkt} \ \ \leq \ \ \overline{x}_{jkt} d_{jt}
& \forall \ j \in \mathcal{J}, \ k \in \overline{\mathcal{K}}_{jt}, \ t = 0,1,\ldots,T
\label{eqn:fluid MNL2 upper} \\
& d_{0t} + \sum_{j \in \mathcal{J}} d_{jt} \ \ = \ \ 1
& \forall \ t = 0,1,\ldots,T
\label{eqn:fluid MNL2 unit}
\end{align}

Unlike~\ref{eqn:static MNL5}, problem~\ref{eqn:fluid MNL2} may not be feasible.
For example, suppose that $\mathcal{J} = \{j\}$ and $\overline{\mathcal{K}}_{jt} = \mathcal{K}_{jt}$ for all~$t$.
Let
\[
\underline{d}_{jt} \ \ \defi \ \ \frac{\exp\left(- \sum_{k \in \mathcal{K}_{jt}} \overline{x}_{jkt}\right)}{1 + \exp\left(- \sum_{k \in \mathcal{K}_{jt}} \overline{x}_{jkt}\right)}
\]
denote the minimum feasible demand.
If $\sum_{t = 0}^{T} \lambda_{t} a_{rj} \underline{d}_{jt} > b_{r}$ for some~$r$, then problem~\ref{eqn:fluid MNL2} is not feasible, and hence problem~\ref{eqn:fluid MNL1} is not feasible.

Examples in Section \ref{sec:remarks} show that it is possible for problem~\ref{eqn:fluid MNL1} to be infeasible while problem~\ref{eqn:fluid MNL2} is feasible, either with~\ref{eqn:fluid MNL2} having an optimal solution or being unbounded.
Next we show that the status of problems~\ref{eqn:fluid MNL1} and~\ref{eqn:fluid MNL2} is easy to determine by solving only~\ref{eqn:fluid MNL2}.

\begin{lemma}
\label{lem:fpmnl2 optimal}
If \ref{eqn:fluid MNL2} is feasible, then \ref{eqn:fluid MNL2} has an optimal solution if and only if for each $j \in \mathcal{J}$, $t \in \{0, 1, \ldots, T\}$, $k_{1} \in \mathcal{K}_{jt} \setminus \overline{\mathcal{K}}_{jt}$, and $k_{2} \in \mathcal{K}_{jt} \setminus \underline{\mathcal{K}}_{jt}$, it holds that $\phi_{k_{1},t} \le \phi_{k_{2},t}$.
\end{lemma}

\begin{lemma}
\label{lem:fpmnl1=fpmnl2}
Suppose that \ref{eqn:fluid MNL2} has a feasible solution $(d', d_{0}', u')$ with $d' > 0$.
Then, any optimal solution $(d^*, d_{0}^*, u^*)$ for \ref{eqn:fluid MNL2} satisfies $d^* > 0$, $d_{0}^* > 0$.
\end{lemma}

Next we summarize the results for~\ref{eqn:fluid general} under the MNL choice model.

\begin{theorem}
\label{thm:fpmnl1=fpmnl2}
\ref{eqn:fluid MNL1} can be solved by solving \ref{eqn:fluid MNL2} and taking into account the following possibilities:
\begin{enumerate}
\item
If \ref{eqn:fluid MNL2} is infeasible, then \ref{eqn:fluid MNL1} is infeasible.
\item
If \ref{eqn:fluid MNL2} is feasible, then the following possibilities hold:
\begin{enumerate}
\item
If \ref{eqn:fluid MNL2} is unbounded (which happens if and only if for some $j \in \mathcal{J}$, $t \in \{0, 1, \ldots, T\}$, $k_{1} \in \mathcal{K}_{jt} \setminus \overline{\mathcal{K}}_{jt}$, and $k_{2} \in \mathcal{K}_{jt} \setminus \underline{\mathcal{K}}_{jt}$, it holds that $\phi_{k_{1},t} > \phi_{k_{2},t}$), then the following possibilities hold:
\begin{enumerate}
\item
If \ref{eqn:fluid MNL2} does not have a feasible solution $(d', d_{0}', u')$ with $d' > 0$, then \ref{eqn:fluid MNL1} is infeasible.
\item
If \ref{eqn:fluid MNL2} has a feasible solution $(d', d_{0}', u')$ with $d' > 0$, then \ref{eqn:fluid MNL1} is unbounded.
\end{enumerate}
\item
If \ref{eqn:fluid MNL2} has an optimal solution $(d^*, d_{0}^*, u^*)$ (which happens if and only if for each $j \in \mathcal{J}$, $t \in \{0, 1, \ldots, T\}$, $k_{1} \in \mathcal{K}_{jt} \setminus \overline{\mathcal{K}}_{jt}$, and $k_{2} \in \mathcal{K}_{jt} \setminus \underline{\mathcal{K}}_{jt}$, it holds that $\phi_{k_{1},t} \le \phi_{k_{2},t}$), then the following possibilities hold:
\begin{enumerate}
\item
If $d^* > 0$, then $(d^*, x^*)$ with $x^*$ given by $x_{jkt}^* = u_{jkt}^* / d_{jt}^*$ for every $j \in \mathcal{J}$, $t \in \{0,1,\ldots,T\}$, and $k \in \mathcal{K}_{jt}$, is an optimal solution for \ref{eqn:fluid MNL1}.
\item
If $d_{jt}^* = 0$ for some $j \in \mathcal{J}$ and $t \in \{0,1,\ldots,T\}$, then \ref{eqn:fluid MNL1} is infeasible.
\end{enumerate}
\end{enumerate}
\end{enumerate}
\end{theorem}

%------------------------------------------------------------------------

\subsection{\ref{eqn:static general} and \ref{eqn:fluid general} under the MC Model}

Next we consider problem \ref{eqn:fluid general} (and as a special case, \ref{eqn:static general}) in which the choice probabilities $\check{P}_{jt}(y_{t})$ are given by a Markov chain choice model. 
The Markov chain choice model here is similar to the model considered in \cite{dong2019pricing}, adjusted to make provision for multiple attributes.
The choice probabilities are given by $\check{P}_{jt}(y_{t}) = \check{Q}_{jt}(y_{t}) \check{V}_{jt}(y_{t})$ , where $\check{V}_{t}(y_{t}) = (\check{V}_{jt}(y_{t}) \, : \, j \in \mathcal{J})$ is the \emph{unique} solution of the linear system
\[
\check{V}_{jt}(y_{t}) \ \ = \ \ \theta_{jt} + \sum_{i \in \mathcal{J}_{t}} \left[1 - \check{Q}_{it}(y_{t})\right] \rho_{ijt} \check{V}_{it}(y_{t}) \quad , \quad j \in \mathcal{J}
\tag{$\mathsf{MC}$}
\label{eqn:sys-mc}
\]
for every $t = 0,1,\ldots,T$. 
Here, the model parameters $\theta_{t} \defi (\theta_{jt} \, : \, j \in \mathcal{J})$, $\rho_{t} \defi (\rho_{ijt} \, : \, i,j \in \mathcal{J})$ and functions $(\check{Q}_{jt}(y_{t}) \ : \ j \in \mathcal{J})$ should satisfy the following assumptions:
\begin{enumerate}
\item
$\theta_{t} \geq 0$ and $\sum_{j \in \mathcal{J}} \theta_{jt} > 0$ (otherwise the problem is trivial); 
\item
$\rho_{t} \geq 0$ and $I - \rho_{t}$ is nonsingular ($I$ is an identity matrix with proper dimensionality); 
\item
$\check{Q}_{jt}(y_{t}) \in [0,1]$ for all feasible $y_{t}$.
\end{enumerate}
Note that $(I - \rho_{t})^{-1} = I + \sum_{n=1}^{\infty} \rho_{t}^{n}$. Thus, $I - \rho_{t}$ is non-singular if and only if $\sum_{n=1}^{\infty} \rho_{t}^{n} < \infty$.
Let $\rho_{ijt}'(y_{t}) \defi \left[1 - \check{Q}_{it}(y_{t})\right] \rho_{ijt}$ for every $i,j \in \mathcal{J}$, and let $\rho_{t}'(y_{t}) \defi (\rho_{ijt}'(y_{t}) \, : \, i,j \in \mathcal{J})$. 
Since $\check{Q}_{jt}(y_{t}) \in [0,1]$, we have $\sum_{n=1}^{\infty} (\rho_{t}'(y_{t}))^{n} \le \sum_{n=1}^{\infty} \rho_{t}^{n} < \infty$, which implies that \ref{eqn:sys-mc} has a unique solution $\check{V}_{t}(y_{t}) = \left(I - (\rho_{t}'(y_{t}))^{\top}\right)^{-1} \theta_{t}$ for any feasible $y_{t}$.

\begin{remark}
\cite{dong2019pricing} provided the following interpretation to the MC model:
Each arriving customer first visits product $j$ with probability $\theta_{jt}$.
A customer visiting product $j$ purchases the product with probability $\check{Q}_{jt}(y_{t})$, or transition from the product with probability $1 - \check{Q}_{jt}(y_{t})$.
When a customer choose to transition from product $j$, she either transition to another product $i$ with probability $\rho_{jit}$, or transition to the no purchase option (and leaves the system) with probability $1 - \sum_{i \in \mathcal{J}} \rho_{jit}$.
In this way, each arriving customer transitions between different products until purchasing one of the products, or deciding to leave without making a purchase.
Thus, $\check{V}_{jt}(y)$ gives the expected number of times that a customer visits product $j$.

\noindent However, as \cite{dong2019pricing} stated, we should not view the MC model as a faithful model of the mental thought process of the customers when they purchase a product. The value of the MC model comes from the facts that (1) this model is compatible with the random utility maximization principle (see Theorem 1 in \cite{dong2019pricing}), (2) its parameter values can be set so that the purchase probabilities under it become identical to the generalized attraction model, which subsumes the MNL model (see Lemma 2 in \cite{dong2019pricing}), and (3) it yields tractable optimization problems (see \cite{blanchet2016markov}, \cite{feldman2017revenue} and \cite{dong2019pricing}). 
\end{remark}

\begin{remark}
\cite{dong2019pricing} showed that parameter values of the MC model can be set so that the purchase probabilities under it become identical to the MNL model. 
The result is established assuming that price attribute is lower bounded for all products. 
We can extend this result to the case with multiple attributes, assuming that all attributes are lower bounded. 
See Section \ref{sec:remarks}. 
\end{remark}

In this section, we let
\[
\check{Q}_{jt}(y_{t}) \ \ = \ \ 
		\exp\left( \alpha_{jt} - \sum_{k \in \mathcal{K}_{jt}} \beta_{kt} y_{jkt} \right)
		\quad , \quad j \in \mathcal{J} \ , \ t = 0,1,\ldots,T 
\]
Similar to the previous discussion, we let $x_{jkt} \defi \beta_{kt} y_{jkt} - \alpha_{jt} / K_{jt}$, $\underline{x}_{jkt} \defi \beta_{kt} \underline{y}_{jkt} - \alpha_{jt} / K_{jt}$, $\overline{x}_{jkt} \defi \beta_{kt} \overline{y}_{jkt} - \alpha_{jt} / K_{jt}$, $\phi_{kt} \defi \beta_{kt} \check{\phi}_{kt} > 0$, $\psi_{jt} \defi \check{\psi}_{jt} + (\alpha_{jt} / K_{jt}) \sum_{k \in \mathcal{K}_{jt}} \check{\phi}_{kt}$, for every $j \in \mathcal{J}$, $t = 0,1,\ldots,T$, and $k \in \mathcal{K}_{jt}$. In addition, to make Assumption 3 above hold, we need 
\[
\sum_{k \in \underline{\mathcal{K}}_{j}} \underline{x}_{jkt} \geq 0
		\quad , \quad j \in \mathcal{J} \ , \ t = 0,1,\ldots,T 
\]
and $\mathcal{K}_{j} = \underline{\mathcal{K}}_{j}$ (all attributes are lower bounded). 

Let $x \defi (x_{jkt}, \, j \in \mathcal{J}, \, k \in \mathcal{K}_{jt}, \, t = 0,1,\ldots,T)$, $d \defi (d_{jt}, \, j \in \mathcal{J}, \, t = 0,1,\ldots,T)$, $v \defi (v_{jt}, \, t = 0,1,\ldots,T)$, and $u \defi (u_{jkt}, \, j \in \mathcal{J}, \, k \in \mathcal{K}_{jt}, \, t = 0,1,\ldots,T)$. 
We consider the problem
\begin{align}
\max_{v, \, d, \, x} \quad & \sum_{t = 0}^{T} \lambda_{t} \sum_{j \in \mathcal{J}} \left(\sum_{k \in \mathcal{K}_{jt}} \phi_{kt} x_{jkt} - \psi_{jt}\right) d_{jt}
\tag{\textsf{$\mathsf{FP_1^{MC}}$}}
\label{eqn:fluid MC1} \\
\text{s.t.} \quad
& d_{jt} \ \ = \ \ \exp\left( - \sum_{k \in \mathcal{K}_{j}} x_{ikt} \right) v_{jt}
		& \forall \ j \in \mathcal{J}, \ t = 0,1,\ldots,T \nonumber \\
& v_{jt} \ \ = \ \ \theta_{jt} + \sum_{i \in \mathcal{J}} \rho_{ij} (v_{it} - d_{it})
		& \forall \ j \in \mathcal{J}, \ t = 0,1,\ldots,T \nonumber \\
& \sum_{t = 0}^{T} \lambda_{t} \sum_{j \in \mathcal{J}} a_{rj} d_{jt} \ \ \le \ \ b_{r}
		& \forall \ r \in \mathcal{R} \nonumber \\
& x_{jkt} \ \ \geq \ \ \underline{x}_{jkt}
		& \forall \ j \in \mathcal{J}, \ k \in \underline{\mathcal{K}}_{jt}, 
		\ t = 0,1,\ldots,T \nonumber \\
& x_{jkt} \ \ \leq \ \ \overline{x}_{jkt}
		& \forall \ j \in \mathcal{J}, \ k \in \underline{\mathcal{K}}_{jt}, 
		\ t = 0,1,\ldots,T \nonumber
\end{align}
and its following convex conic relaxation:
\begin{subequations}
\begin{align}
\max_{v, \, d, \, u} \quad & \sum_{t = 0}^{T} \lambda_{t} \sum_{j \in \mathcal{J}} \left(\sum_{k \in \mathcal{K}_{jt}} \phi_{kt} u_{jkt} - \psi_{jt} d_{jt}\right)
\tag{\textsf{$\mathsf{FP_2^{MC}}$}}
\label{eqn:fluid MC2} \\
\text{s.t.} \quad
& v_{jt} \ \ = \ \ \theta_{jt} + \sum_{i \in \mathcal{J}} \rho_{ij} (v_{it} - d_{it})
		& \forall \ j \in \mathcal{J}, \ t = 0,1,\ldots,T \\
& \sum_{t = 0}^{T} \lambda_{t} \sum_{j \in \mathcal{J}} a_{rj} d_{jt} \ \ \le \ \ b_{r}
		& \forall \ r \in \mathcal{R}
		\label{eqn:fluid MC2 resource} \\
& \left(v_{jt}, d_{jt}, \sum_{k \in \mathcal{K}_{jt}} u_{jkt}\right) \ \ \in \ \ \mathcal{K}_{\exp}
		& \forall \ j \in \mathcal{J}, \ t = 0,1,\ldots,T
		\label{eqn:fluid MC2 cone} \\
& \left(d_{jt}, v_{jt}, - \sum_{k \in \mathcal{K}_{jt}} \overline{x}_{jkt} v_{jt}\right) 
		\ \ \in \ \ \mathcal{K}_{\exp}
		& \forall \ j \in \overline{\mathcal{J}}_{t}, \ t = 0,1,\ldots,T
		\label{eqn:fluid MC2 cone dummy} \\
& u_{jkt} \ \ \geq \ \ \underline{x}_{jkt} d_{jt}
		& \forall \ j \in \mathcal{J}, \ k \in \underline{\mathcal{K}}_{jt}, \ t = 0,1,\ldots,T
		\label{eqn:fluid MC2 lower} \\
& u_{jkt} \ \ \leq \ \ \overline{x}_{jkt} d_{jt}
		& \forall \ j \in \mathcal{J}, \ k \in \overline{\mathcal{K}}_{jt}, \ t = 0,1,\ldots,T
		\label{eqn:fluid MC2 upper} 
\end{align}
\end{subequations}

Next we summarize the results for~\ref{eqn:fluid general} under the MC choice model.

\begin{theorem}
\label{thm:fpmc1=fpmc2}
\ref{eqn:fluid MC1} can be solved by solving \ref{eqn:fluid MC2} and taking into account the following:
\begin{enumerate}
\item
If \ref{eqn:fluid MC2} is infeasible, then \ref{eqn:fluid MC1} is infeasible.
\item
If \ref{eqn:fluid MC2} is feasible, then:
\begin{enumerate}
\item
\ref{eqn:fluid MC2} is bounded, and has an optimal solution. 
\item
\ref{eqn:fluid MC1} is bounded, and has an optimal solution. 
\item
Let $(v^*, d^*, u^*)$ be an optimal solution of \ref{eqn:fluid MC2}. Let $x_{jkt}^* = u_{jkt}^* / d_{jt}^*$ for every $j \in \mathcal{J}$, $t \in \{0,1,\ldots,T\}$, and $k \in \mathcal{K}_{jt}$ such that $d_{jt}^* > 0$, and let $x_{jkt}^* = \underline{x}_{jkt}^*$ otherwise. Then $(v^*, d^*, x^*)$ is an optimal solution of \ref{eqn:fluid MC1}. 
\end{enumerate}
\end{enumerate}
\end{theorem}

\begin{remark}
Recall that, with the three assumptions hold, any feasible $y$ corresponds to a unique solution of the system \ref{eqn:sys-mc}. Thus, if there is no resource constraint in \ref{eqn:fluid MC1}, then it must be feasible (assuming that the upper and lower bounds of attributes are proper, such that a feasible $y$ exists). Therefore, while \ref{eqn:static general} is a special case of \ref{eqn:fluid general}, the problem \ref{eqn:static general} under the MC model is always feasible, and has an optimal solution. 
\end{remark}

%------------------------------------------------------------------------

\subsection{\ref{eqn:static general} and \ref{eqn:fluid general} under the NL Model}

Next we consider problem \ref{eqn:fluid general} (and as a special case, \ref{eqn:static general}) in which the choice probabilities are given by a nested logit model with non-overlapping nests, and one single null alternative. 
Let $\mathcal{I}$ denotes the set of all nests. For each $i \in \mathcal{I}$, let $\mathcal{J}_{i} \in \mathcal{J}$ denote the set of products in nest $i$.
In this model: (1) $\mathcal{J}_{i} \neq \varnothing$ for each $i \in \mathcal{I}$, (2) $\mathcal{J} = \cup_{i \in \mathcal{I}} \mathcal{J}_{i}$, and (3) $\mathcal{J}_{i} \cap \mathcal{J}_{i'} = \varnothing$ for each pair of two different nests $i, i' \in \mathcal{I}$. 
Then the choice probabilities $\check{P}_{jt}(y_{t})$ are: 
\[
\check{P}_{jt}(y) \ \ = \ \ 
		\frac{\exp\left(\alpha_{jt} - \sum_{k \in \mathcal{K}_{jt}} \beta_{kt} y_{jkt}\right)
		\left( \sum_{j' \in \mathcal{J}_{i}} \exp\left( \alpha_{j't} 
		- \sum_{k \in \mathcal{K}_{j't}} \beta_{kt} y_{j'kt} \right)\right)^{\gamma_{it} - 1}}
		{1 + \sum_{i' \in \mathcal{I}} \left( \sum_{j' \in \mathcal{J}_{i'}} 
		\exp\left(\alpha_{j't} - \sum_{k \in \mathcal{K}_{j't}} \beta_{kt} y_{j'kt}
		\right)\right)^{\gamma_{i't}}}
\]
for each $i$ and $j$ such that $j \in \mathcal{J}_{i}$. 
Same as in the MNL model, $\alpha_{jt}$ are the product ``baseline attractiveness'' parameters, and $\beta_{kt} > 0$ are the attribute ``sensitivity'' parameters. 
In the nested logit model, $\gamma_{it} > 0$ denotes the ``nest dissimilarity'' parameters. Note that, when $\gamma_i > 1$ for some $i \in \mathcal{I}$, the model can be inconsistent with the random utility theory. See the following quote from \cite{train2009discrete}: 

\begin{displayquote}
The value of $\gamma_i$ must be within a particular range for the model to be consistent with utility-maximizing behavior. If $\gamma_i$ is between zero and one for all $i \in \mathcal{I}$, the model is consistent with utility maximization for all possible values of the explanatory variables. For $\gamma_i$ greater than one, the model is consistent with utility-maximizing behavior for some range of the explanatory variables but not for all values. \cite{kling1995empirical} and \cite{herriges1996testing} provide tests of consistency of nested logit with utility maximization when $\gamma_i > 1$; and \cite{train1987demand} and \cite{lee1999calling} provide examples of models for which $\gamma_i > 1$. A negative value of $\gamma_i$ is inconsistent with utility maximization and implies that improving the attributes of an alternative (such as lowering its price) can decrease the probability of the alternative being chosen. With positive $\gamma_i$, the nested logit approaches the “elimination by aspects” model of \cite{tversky1972elimination} as $\gamma_i \to 0$.
\end{displayquote}
In this part, we let $\mathcal{I}_{\leq} := \{i \in \mathcal{I} : \gamma_i \leq 1\}$, let $\mathcal{I}_{>} := \{i \in \mathcal{I} : \gamma_i > 1\}$, let $\mathcal{J}_{\leq} := \cup_{i \in \mathcal{I}_{\leq}} \mathcal{J}_{i}$, and let $\mathcal{J}_{>} := \cup_{i \in \mathcal{I}_{>}} \mathcal{J}_{i}$. Different types of nests will be discussed separately. In addition, similar to the previous parts, we let $x_{jkt} \defi \beta_{kt} y_{jkt} - \alpha_{jt} / K_{jt}$, $\underline{x}_{jkt} \defi \beta_{kt} \underline{y}_{jkt} - \alpha_{jt} / K_{jt}$, $\overline{x}_{jkt} \defi \beta_{kt} \overline{y}_{jkt} - \alpha_{jt} / K_{jt}$, $\phi_{kt} \defi \beta_{kt} \check{\phi}_{kt} > 0$, $\psi_{jt} \defi \check{\psi}_{jt} + (\alpha_{jt} / K_{jt}) \sum_{k \in \mathcal{K}_{jt}} \check{\phi}_{kt}$, for every $j \in \mathcal{J}$, $t = 0,1,\ldots,T$, and $k \in \mathcal{K}_{jt}$.

Let $x \defi (x_{jkt}, \, j \in \mathcal{J}, \, k \in \mathcal{K}_{jt}, \, t = 0,1,\ldots,T)$, $d \defi (d_{jt}, \, j \in \mathcal{J}, \, t = 0,1,\ldots,T)$ , $p \defi (p_{it}, \, i \in \mathcal{I}, \, t = 0,1,\ldots,T)$, and $p_{0} \defi (p_{0t}, \, t = 0,1,\ldots,T)$.
In addition, let $u \defi (u_{jkt}, \, j \in \mathcal{J}, \, k \in \mathcal{K}_{jt}, \, t = 0,1,\ldots,T)$, $u_{>} \defi (u_{jkt}, \, j \in \mathcal{J}_{>}, \, k \in \mathcal{K}_{jt}, \, t = 0,1,\ldots,T)$, $u_{\leq} \defi (u_{jkt}, \, j \in \mathcal{J}_{\leq}, \, k \in \mathcal{K}_{jt}, \, t = 0,1,\ldots,T)$, and $v_{\leq} \defi (v_{ikt}, \, i \in \mathcal{I}_{\leq}, \, k \in \mathcal{K}_{it}', \, t = 0,1,\ldots,T)$. 
We consider the optimization problem
{\small
\begin{align}
\max_{d, \, p, \, x} \quad & \sum_{t = 0}^{T} \lambda_{t} \sum_{j \in \mathcal{J}} 
		\left(\sum_{k \in \mathcal{K}_{jt}} \phi_{kt} x_{jkt} - \psi_{jt}\right) d_{jt}
\tag{\textsf{$\mathsf{FP_1^{NL}}$}}
\label{eqn:fluid NL1} \\
\text{s.t.} \quad
& d_{jt} \ \ = \ \ p_{it} \left(
		\frac{\exp\left(\alpha_{jt} - \sum_{k \in \mathcal{K}_{jt}} \beta_{kt} y_{jkt}\right)}
		{\sum_{j' \in \mathcal{J}_{i}} \exp\left(\alpha_{j't} 
		- \sum_{k \in \mathcal{K}_{j't}} \beta_{kt} y_{j'kt} \right)} \right)
& \forall \ i \in \mathcal{I}, \ j \in \mathcal{J}_{i}, \ t = 0,1,\ldots,T \nonumber \\
& p_{it} \ \ = \ \ 
		\frac{\left( \sum_{j' \in \mathcal{J}_{i}} \exp\left( \alpha_{j't} 
		- \sum_{k \in \mathcal{K}_{j't}} \beta_{kt} y_{j'kt} \right)\right)^{\gamma_{it}}}
		{1 + \sum_{i' \in \mathcal{I}} \left( \sum_{j' \in \mathcal{J}_{i'}} 
		\exp\left(\alpha_{j't} - \sum_{k \in \mathcal{K}_{j't}} \beta_{kt} y_{j'kt}
		\right)\right)^{\gamma_{it}}}
& \forall \ i \in \mathcal{I}, \ t = 0,1,\ldots,T \nonumber \\
& \sum_{t = 0}^{T} \lambda_{t} \sum_{j \in \mathcal{J}} a_{rj} d_{jt} \ \ \le \ \ b_{r}
		& \forall \ r \in \mathcal{R} \nonumber \\
& x_{jkt} \ \ \geq \ \ \underline{x}_{jkt}
& \forall \ j \in \mathcal{J}, \ k \in \underline{\mathcal{K}}_{jt}, \ t = 0,1,\ldots,T \nonumber \\
& x_{jkt} \ \ \leq \ \ \overline{x}_{jkt}
& \forall \ j \in \mathcal{J}, \ k \in \underline{\mathcal{K}}_{jt}, \ t = 0,1,\ldots,T \nonumber
\end{align}
}%
and the convex program:
{\small
\begin{align}
\max_{\substack{d, \ p, \ p_{0t}, \\ v_{\leq}, \ u_{>}}} \quad
& \sum_{t = 0}^{T} \lambda_{t} \sum_{i \in \mathcal{I}_{\leq}} 
		\sum_{k \in \mathcal{K}_{it}'} \phi_{kt} v_{ikt} \nonumber \\
		& + \sum_{t = 0}^{T} \lambda_{t} \sum_{j \in \mathcal{J}_{>}} 
		\sum_{k \in \mathcal{K}_{jt}} \phi_{kt} u_{jkt} \nonumber \\
		& - \sum_{t = 0}^{T} \lambda_{t} \sum_{j \in \mathcal{J}} \psi_{jt} d_{jt}
\tag{\textsf{$\mathsf{FP_2^{NL}}$}} 
\label{eqn:fluid NL2} \\
\text{s.t.} \quad
& \sum_{t = 0}^{T} \lambda_{t} \sum_{j \in \mathcal{J}} a_{rj} d_{jt} \ \ \le \ \ b_{r}
		& \forall \ r \in \mathcal{R} \nonumber \\
& \left(\frac{1}{\gamma_{it}} - 1\right)
		p_{it} \ln\left(\frac{p_{0t}}{p_{it}}\right)
		+ \sum_{j \in \mathcal{J}_i} d_{jt} \ln\left(\frac{p_{0t}}{d_{jt}}\right)
		\ \ \geq \ \
		\sum_{k \in \mathcal{K}_{it}'} v_{ikt}
		& \forall \ i \in \mathcal{I}_{\leq}, \ t = 0,1,\ldots,T 
		\nonumber \\
& \left(\frac{1}{\gamma_{it}} - 1\right)
		p_{0t} \ln\left(\frac{p_{0t}}{p_{it}}\right)
		+ p_{0t} \ln\left(\frac{p_{0t}}{d_{jt}}\right)
		\ \ \leq \ \
		p_{0t} \sum_{k \in \mathcal{K}_{jt}} \overline{x}_{jkt}^{[i]}
		& \forall \ i \in \mathcal{I}_{\leq}, \ j \in \overline{\mathcal{J}}_{i}, \ t = 0,1,\ldots,T 
		\nonumber \\
& \left(1 - \frac{1}{\gamma_{it}}\right)
		d_{jt} \ln\left(\frac{p_{it}}{d_{jt}}\right)
		+ \left(\frac{1}{\gamma_{it}}\right)
		d_{jt} \ln\left(\frac{p_{0t}}{d_{jt}}\right)
		\ \ \geq \ \
		\sum_{k \in \mathcal{K}_{jt}} u_{jkt}
		& \forall \ i \in \mathcal{I}_{>}, \ j \in \mathcal{J}_{i}, \ t = 0,1,\ldots,T 
		\nonumber \\
& p_{0t} + \sum_{i \in \mathcal{I}} p_{it} \ \ = \ \ 1
		& \forall \ t = 0,1,\ldots,T 
		\nonumber \\
& p_{it} \ \ = \ \ \sum_{j \in \mathcal{J}_{i}} d_{jt}
		& \forall \ i \in \mathcal{I}, \ t = 0,1,\ldots,T 
		\nonumber \\
& d \ \ > \ \ 0
		\nonumber \\
& v_{ikt} \ \ \geq \ \ \sum_{j \in \mathcal{J}_{i}} \underline{x}_{jkt} d_{jt}
& \forall \ i \in \mathcal{I}_{\leq}, \ k \in \underline{\mathcal{K}}_{it}', \ t = 0,1,\ldots,T 
		\nonumber \\
& v_{ikt} \ \ \leq \ \ \sum_{j \in \mathcal{J}_{i}} \overline{x}_{jkt} d_{jt}
& \forall \ i \in \mathcal{I}_{\leq}, \ k \in \underline{\mathcal{K}}_{it}', \ t = 0,1,\ldots,T 
		\nonumber \\
& u_{jkt} \ \ \geq \ \ \underline{x}_{jkt} d_{jt}
& \forall \ j \in \mathcal{J}_{>}, \ k \in \underline{\mathcal{K}}_{jt}, \ t = 0,1,\ldots,T 
		\nonumber \\
& u_{jkt} \ \ \leq \ \ \overline{x}_{jkt} d_{jt}
& \forall \ j \in \mathcal{J}_{>}, \ k \in \overline{\mathcal{K}}_{jt}, \ t = 0,1,\ldots,T
		\nonumber 
\end{align}
}%
where $\overline{\mathcal{K}}_{it}' \defi \cap_{j \in \mathcal{J}_{i}} \overline{\mathcal{K}}_{jt}$, $\underline{\mathcal{K}}_{it}' \defi \cap_{j \in \mathcal{J}_{i}} \underline{\mathcal{K}}_{jt}$ and $\mathcal{K}_{it}' \defi \cup_{j \in \mathcal{J}_{i}} \mathcal{K}_{jt}$. 
Then the following result holds: 

\begin{theorem}
\label{thm:fpnl1=fpnl2}
If the following conditions hold: 

\noindent Condition (i): for any feasible solution $(d, p, p_{0}, v_{\leq}, u_{>})$ to \ref{eqn:fluid NL2} at which
\begin{align*}
	& \left(\frac{1}{\gamma_{it}} - 1\right)
		p_{it} \ln\left(\frac{p_{0t}}{p_{it}}\right)
		+ \sum_{j \in \mathcal{J}_i} d_{jt} \ln\left(\frac{p_{0t}}{d_{jt}}\right)
		\ \ = \ \
		\sum_{k \in \mathcal{K}_{it}'} v_{ikt}
		& \forall \ i \in \mathcal{I}_{\leq}, \ t = 0,1,\ldots,T 
\end{align*}
the system
\begin{align*}
& \left(\frac{1}{\gamma_{it}} - 1\right)
		d_{jt} \ln\left(\frac{p_{0t}}{p_{it}}\right)
		+ d_{jt} \ln\left(\frac{p_{0t}}{d_{jt}}\right)
		\ \ = \ \
		\sum_{k \in \mathcal{K}_{jt}} u_{jkt} 
		& \forall \ i \in \mathcal{I}_{\leq}, \ j \in \mathcal{J}_{i}, \ t = 0,1,\ldots,T \\
& v_{ikt} \ \ = \ \ \sum_{j \in \mathcal{J}_{ikt}'} u_{jkt}
		& \forall \ i \in \mathcal{I}_{\leq}, \ k \in \underline{\mathcal{K}}_{it}', \ t = 0,1,\ldots,T \\
& u_{jkt} \ \ \geq \ \ \underline{x}_{jkt} d_{jt}
& \forall \ j \in \mathcal{J}_{\leq}, \ k \in \underline{\mathcal{K}}_{jt}, \ t = 0,1,\ldots,T 
		\nonumber \\
& u_{jkt} \ \ \leq \ \ \overline{x}_{jkt} d_{jt}
& \forall \ j \in \mathcal{J}_{\leq}, \ k \in \overline{\mathcal{K}}_{jt}, \ t = 0,1,\ldots,T
		\nonumber 
\end{align*}
has a solution $u_{\leq}$. Here $\mathcal{J}_{ikt}' \defi \{j \in \mathcal{J}_{i} \ : \ k \in \mathcal{K}_{jt} \}$. 

\noindent Condition (ii): $\overline{\mathcal{K}}_{jt} \subsetneq \mathcal{K}_{jt}$ for every $j \in \mathcal{J}_{>}$, i.e. $\sum_{k \in \mathcal{K}_{jt}} u_{jkt}$ is not upper bounded for any $j \in \mathcal{J}_{>}$. 

\noindent then \ref{eqn:fluid NL1} can be solved by solving \ref{eqn:fluid NL2} and taking into account the following possibilities:
\begin{enumerate}
\item
If \ref{eqn:fluid NL2} is infeasible, then \ref{eqn:fluid NL1} is infeasible.
\item
If \ref{eqn:fluid NL2} is feasible and unbounded, then \ref{eqn:fluid NL1} is feasible and unbounded. 
\item
If \ref{eqn:fluid NL2} is feasible and has an optimal solution $(d^*, p^*, p_{0}^*, v_{\leq}^*, u_{>}^*)$, then $(d^*, p^*, x^*)$ with $x^*$ given by $x_{jkt}^* = u_{jkt}^* / d_{jt}^*$ for every $j \in \mathcal{J}$, $t \in \{0,1,\ldots,T\}$, and $k \in \mathcal{K}_{jt}$, is an optimal solution for \ref{eqn:fluid NL1}. (Here $u_{>}^*$ is part of the optimal solution $(d^*, p^*, p_{0}^*, v_{\leq}^*, u_{>}^*)$ to \ref{eqn:fluid NL2}, while $u_{\leq}^*$ is obtained by solving the system in Condition (i), using $(d^*, p^*, p_{0}^*, v_{\leq}^*, u_{>}^*)$ as an input.)
\end{enumerate}
\end{theorem}

\begin{remark}
A simple sufficient condition for Condition (i) to hold in Theorem \ref{thm:fpmnl1=fpmnl2}, is that for any $i \in \mathcal{I}_{\leq}$, there is one unbounded attribute $k$ that is shared by all products in $\mathcal{J}_{i}$ (for example, the price attribute). 
That is, for any $i \in \mathcal{I}_{\leq}$, $j \in \mathcal{J}_{i}$ and $t = 0,1,\ldots,T$, we need $k \in \mathcal{K}_{jt}$ (which implies $\mathcal{J}_{i} = \mathcal{J}_{ikt}'$), $k \notin \overline{\mathcal{K}}_{jt}$, and $k \notin \underline{\mathcal{K}}_{jt}$. 
Then we can split all other attribute $k' \in \mathcal{J}_{i}$ for each $i \in \mathcal{I}_{\leq}$, $j \in \mathcal{J}_{i}$ and $t = 0,1,\ldots,T$, such that each of those attributes only apply to one product (i.e. $\mathcal{J}_{ik't}'$ are all singletons). 
Now, given any feasible solution $(d, p, p_{0}, v_{\leq}, u_{>})$ to \ref{eqn:fluid NL2} at which
\begin{align*}
	& \left(\frac{1}{\gamma_{it}} - 1\right)
		p_{it} \ln\left(\frac{p_{0t}}{p_{it}}\right)
		+ \sum_{j \in \mathcal{J}_i} d_{jt} \ln\left(\frac{p_{0t}}{d_{jt}}\right)
		\ \ = \ \
		\sum_{k \in \mathcal{K}_{it}'} v_{ikt}
		& \forall \ i \in \mathcal{I}_{\leq}, \ t = 0,1,\ldots,T 
\end{align*}
we can find a feasible solution to the system in Condition (i) as follow: Consider any $i \in \mathcal{I}_{\leq}$ and $j \in \mathcal{J}_{i}$. Let  
\begin{align*}
& u_{jk't} \ \ = \ \ v_{ik't}
		& \forall \ t = 0,1,\ldots,T 
\end{align*}
for any $k' \neq k$, and let 
\begin{align*}
& u_{jkt} \ \ = \ \ \left(\frac{1}{\gamma_{it}} - 1\right)
		d_{jt} \ln\left(\frac{p_{0t}}{p_{it}}\right)
		+ d_{jt} \ln\left(\frac{p_{0t}}{d_{jt}}\right) 
		- \sum_{k' \in \mathcal{K}_{jt} \setminus \{k\}} u_{jk't}
		& \forall \ t = 0,1,\ldots,T 
\end{align*}
Note that we automatically have
\begin{align*}
& v_{ikt} \ \ = \ \ \sum_{j \in \mathcal{J}_{i}} u_{jkt}
		& \forall \ t = 0,1,\ldots,T
\end{align*}
which implies that the $u_{\leq}$ constructed is a feasible solution to the system in Condition (i). 
\end{remark}

Similar to the case under the MNL model, the feasible set of~\ref{eqn:fluid NL2} may not be closed. Next we relax~\ref{eqn:fluid NL2} to make the feasible set closed, and then we show that optimal solutions are not affected by the relaxation. Consider the following convex conic program: 
{\small
\begin{subequations}
\begin{align}
\max_{\substack{d, \ p, \ p_{0}, \ v_{\leq}, \ u_{>} \\ 
		e, \ f_{\leq}, \ g_{>}, \ r_{\leq}, \ s_{>}}} \quad
& \sum_{t = 0}^{T} \lambda_{t} \sum_{i \in \mathcal{I}_{\leq}} 
		\sum_{k \in \mathcal{K}_{it}'} \phi_{kt} v_{ikt} \nonumber \\
		& + \sum_{t = 0}^{T} \lambda_{t} \sum_{j \in \mathcal{J}_{>}} 
		\sum_{k \in \mathcal{K}_{jt}} \phi_{kt} u_{jkt} \nonumber \\
		& - \sum_{t = 0}^{T} \lambda_{t} \sum_{j \in \mathcal{J}} \psi_{jt} d_{jt}
\tag{\textsf{$\mathsf{FP_3^{NL}}$}}
\label{eqn:fluid NL3} \\
& \sum_{t = 0}^{T} \lambda_{t} \sum_{j \in \mathcal{J}} a_{rj} d_{jt} \ \ \le \ \ b_{r}
		& \forall \ r \in \mathcal{R} \\
\text{s.t.} \quad
& \left(\frac{1}{\gamma_{it}} - 1\right) r_{it}
		+ \sum_{j \in \mathcal{J}_i} e_{jt}
		\ \ = \ \
		\sum_{k \in \mathcal{K}_{it}'} v_{ikt}
		& \forall \ i \in \mathcal{I}_{\leq}, \ t = 0,1,\ldots,T \\
& \left(\frac{1}{\gamma_{it}} - 1\right) s_{it} + f_{jt}
		\ \ = \ \
		- p_{0t} \sum_{k \in \mathcal{K}_{jt}} \overline{x}_{jkt}^{[i]}
		& \forall \ i \in \mathcal{I}_{\leq}, \ j \in \overline{\mathcal{J}}_{i}, \ t = 0,1,\ldots,T \\
& \left(1 - \frac{1}{\gamma_{it}}\right) g_{jt}
		+ \left(\frac{1}{\gamma_{it}}\right) e_{jt}
		\ \ = \ \
		\sum_{k \in \mathcal{K}_{jt}} u_{jkt}
		& \forall \ i \in \mathcal{I}_{>}, \ j \in \mathcal{J}_{i}, \ t = 0,1,\ldots,T \\
& \left(p_{0t}, d_{jt}, e_{jt}\right) \ \ \in \ \ \mathcal{K}_{\exp}
		& \forall \ i \in \mathcal{I}, \ j \in \mathcal{J}_{i}, \ t = 0,1,\ldots,T \\
& \left(d_{jt}, p_{0t}, f_{jt}\right) \ \ \in \ \ \mathcal{K}_{\exp}
		& \forall \ i \in \mathcal{I}_{\leq}, \ j \in \overline{\mathcal{J}}_{i}, \ t = 0,1,\ldots,T \\
& \left(p_{it}, d_{jt}, g_{jt}\right) \ \ \in \ \ \mathcal{K}_{\exp}
		& \forall \ i \in \mathcal{I}_{>}, \ j \in \mathcal{J}_{i}, \ t = 0,1,\ldots,T \\
& \left(p_{it}, p_{0t}, s_{it}\right) \ \ \in \ \ \mathcal{K}_{\exp}
		& \forall \ i \in \mathcal{I}_{\leq}, \ t = 0,1,\ldots,T \\
& \left(p_{0t}, p_{it}, r_{it}\right) \ \ \in \ \ \mathcal{K}_{\exp}
		& \forall \ i \in \mathcal{I}_{\leq}, \ t = 0,1,\ldots,T \\
& p_{0t} + \sum_{i \in \mathcal{I}} p_{it} \ \ = \ \ 1
		& \forall \ t = 0,1,\ldots,T \\
& p_{it} \ \ = \ \ \sum_{j \in \mathcal{J}_{i}} d_{jt}
		& \forall \ i \in \mathcal{I}, \ t = 0,1,\ldots,T \\
& v_{ikt} \ \ \geq \ \ \sum_{j \in \mathcal{J}_{i}} \underline{x}_{jkt} d_{jt}
& \forall \ i \in \mathcal{I}_{\leq}, \ k \in \underline{\mathcal{K}}_{it}', \ t = 0,1,\ldots,T
\label{eqn:fluid NL3 v lower} \\
& v_{ikt} \ \ \leq \ \ \sum_{j \in \mathcal{J}_{i}} \overline{x}_{jkt} d_{jt}
& \forall \ i \in \mathcal{I}_{\leq}, \ k \in \overline{\mathcal{K}}_{it}', \ t = 0,1,\ldots,T
\label{eqn:fluid NL3 v upper} \\
& u_{jkt} \ \ \geq \ \ \underline{x}_{jkt} d_{jt}
& \forall \ j \in \mathcal{J}_{>}, \ k \in \underline{\mathcal{K}}_{jt}, \ t = 0,1,\ldots,T
\label{eqn:fluid NL3 u lower} \\
& u_{jkt} \ \ \leq \ \ \overline{x}_{jkt} d_{jt}
& \forall \ j \in \mathcal{J}_{>}, \ k \in \overline{\mathcal{K}}_{jt}, \ t = 0,1,\ldots,T
\label{eqn:fluid NL3 u upper} 
\end{align}
\end{subequations}
}%
Then the following result holds: 
\begin{theorem}
\label{thm:fpnl2=fpnl3}
Assume that Condition (i) in Theorem \ref{thm:fpnl1=fpnl2} holds.

\noindent Then, \ref{eqn:fluid NL2} can be solved by solving \ref{eqn:fluid NL3} and taking into account the following:
\begin{enumerate}
\item
If \ref{eqn:fluid NL3} is infeasible, then \ref{eqn:fluid NL2} is infeasible.
\item
If \ref{eqn:fluid NL3} is feasible, then the following possibilities hold:
\begin{enumerate}
\item
If \ref{eqn:fluid NL3} is unbounded (which happens if and only if for some $j \in \mathcal{J}$, $t \in \{0, 1, \ldots, T\}$, $k_{1} \in \mathcal{K}_{jt} \setminus \overline{\mathcal{K}}_{jt}$, and $k_{2} \in \mathcal{K}_{jt} \setminus \underline{\mathcal{K}}_{jt}$, it holds that $\phi_{k_{1},t} > \phi_{k_{2},t}$), then the following possibilities hold:
\begin{enumerate}
\item
If \ref{eqn:fluid NL3} does not have a feasible solution $(d', p', p_{0}', v_{\leq}', u_{>}', e', f_{\leq}', g_{>}', r_{\leq}', s_{>}')$ with $d' > 0$, then \ref{eqn:fluid NL2} is infeasible.
\item
If \ref{eqn:fluid NL3} has a feasible solution $(d', p', p_{0}', v_{\leq}', u_{>}', e', f_{\leq}', g_{>}', r_{\leq}', s_{>}')$ with $d' > 0$, then \ref{eqn:fluid NL2} is unbounded.
\end{enumerate}
\item
If \ref{eqn:fluid NL3} has an optimal solution $(d^*, p^*, p_{0}^*, v_{\leq}^*, u_{>}^*, e^*, f_{\leq}^*, g_{>}^*, r_{\leq}^*, s_{>}^*)$ (which happens if and only if for each $j \in \mathcal{J}$, $t \in \{0, 1, \ldots, T\}$, $k_{1} \in \mathcal{K}_{jt} \setminus \overline{\mathcal{K}}_{jt}$, and $k_{2} \in \mathcal{K}_{jt} \setminus \underline{\mathcal{K}}_{jt}$, it holds that $\phi_{k_{1},t} \le \phi_{k_{2},t}$), then the following possibilities hold:
\begin{enumerate}
\item
If $d^* > 0$, then $(d^*, p^*, p_{0}^*, v_{\leq}^*, u_{>}^*)$ is an optimal solution for \ref{eqn:fluid NL2}. 
\item
If $d_{jt}^* = 0$ for some $j \in \mathcal{J}$ and $t \in \{0,1,\ldots,T\}$, then \ref{eqn:fluid NL2} is infeasible.
\end{enumerate}
\end{enumerate}
\end{enumerate}
\end{theorem}

%------------------------------------------------------------------------
%	Section
%------------------------------------------------------------------------

\section{Remarks}
\label{sec:remarks}

%------------------------------------------------------------------------

\emph{\underline{Examples where \ref{eqn:fluid MNL1} is infeasible, while \ref{eqn:fluid MNL2} is feasible}:}

\begin{example}
\label{ex:FP1 infeasible, FP2 feasible}
In this example \ref{eqn:fluid MNL1} is infeasible, while \ref{eqn:fluid MNL2} has an optimal solution.
Consider the following single-period instance of \ref{eqn:fluid MNL1} with two products indexed $\{1,2\}$.
Each product~$i$ has one attribute also indexed~$i$.
The time index is omitted.
\begin{align}
\max_{d, \, x} \quad & \big\{x_{1,1} d_{1} + x_{2,2} d_{2}\big\}
\tag{$\mathsf{FP^{Example}_{1}}$}
\label{eqn:FP EX1} \\
\text{s.t.} \quad
& d_{1} + d_{2} \ \ \leq \ \ \frac{1}{2} \nonumber \\
& d_{1} \ \ = \ \ \frac{\exp(- x_{1,1})}{1 + \exp(- x_{1,1}) + \exp(- x_{2,2})} \nonumber \\
& d_{2} \ \ = \ \ \frac{\exp(- x_{2,2})}{1 + \exp(- x_{1,1}) + \exp(- x_{2,2})} \nonumber \\
& x_{1,1} \ \ \leq \ \ 0 \nonumber
\end{align}
For any $x$ feasible for \ref{eqn:FP EX1}, it holds that
\begin{align*}
d_{1} + d_{2}
\ \ & = \ \ \frac{\exp(- x_{1,1}) + \exp(- x_{2,2})}{1 + \exp(- x_{1,1}) + \exp(- x_{2,2})} \\
& > \ \ \frac{\exp(- x_{1,1})}{1 + \exp(- x_{1,1})}
\ \ \geq \ \ \frac{\exp(0)}{1 + \exp(0)} \ \ = \ \ \frac{1}{2}
\end{align*}
Thus, \ref{eqn:FP EX1} is infeasible.

The corresponding instance of \ref{eqn:fluid MNL2} is
\begin{align}
\max_{d, \, d_{0}, \, u} \quad & \big\{u_{1,1} + u_{2,2}\big\}
\tag{$\mathsf{FP^{Example}_{2}}$}
\label{eqn:FP EX2} \\
\text{s.t.} \quad
& d_{1} + d_{2} \ \ \leq \ \ \frac{1}{2} \nonumber \\
& \left(d_{0}, d_{1}, u_{1,1}\right) \ \ \in \ \ \mathcal{K}_{\exp} \nonumber \\
& \left(d_{0}, d_{2}, u_{2,2}\right) \ \ \in \ \ \mathcal{K}_{\exp} \nonumber \\
& \left(d_{1}, d_{0}, 0\right) \ \ \in \ \ \mathcal{K}_{\exp} \nonumber \\
& d_{0} + d_{1} + d_{2} \ \ = \ \ 1 \nonumber \\
& u_{1,1} \ \ \leq \ \ 0 \nonumber
\end{align}
The only feasible solution of \ref{eqn:FP EX2} is $d_{0} = 1/2$, $d_{1} = 1/2$, $d_{2} = 0$, $u_{1,1} = 0$, $u_{2,2} = 0$.
Thus, \ref{eqn:FP EX2} has an optimal solution.
\end{example}

\begin{example}
\label{ex:FP1 infeasible, FP2 unbounded}
In this example \ref{eqn:fluid MNL1} is infeasible, while \ref{eqn:fluid MNL2} is unbounded.
Consider the following single-period instance of \ref{eqn:fluid MNL1} with two products indexed $\{1,2\}$.
Product~$1$ has one attribute indexed~$1$, and product~$2$ has two attributes indexed~$1,2$.
The time index is omitted.
\begin{align}
\max_{d, \, x} \quad & x_{1,1} d_{1} + \left(x_{2,1} + 2 x_{2,2}\right) d_{2}
\tag{$\mathsf{FP^{Example}_{3}}$}
\label{eqn:FP EX3} \\
\text{s.t.} \quad
& d_{1} + d_{2} \ \ \leq \ \ \frac{1}{2} \nonumber \\
& d_{1} \ \ = \ \ \frac{\exp(- x_{1,1})}{1 + \exp(- x_{1,1}) + \exp(- x_{2,1} - x_{2,2})} \nonumber \\
& d_{2} \ \ = \ \ \frac{\exp(- x_{2,1} - x_{2,2})}{1 + \exp(- x_{1,1}) + \exp(- x_{2,1} - x_{2,2})} \nonumber \\
& x_{1,1} \ \ \leq \ \ 0 \nonumber
\end{align}
For any $x$ feasible for \ref{eqn:FP EX3}, it holds that
\begin{align*}
d_{1} + d_{2}
\ \ & = \ \ \frac{\exp(- x_{1,1}) + \exp(- x_{2,1} - x_{2,2})}{1 + \exp(- x_{1,1}) + \exp(- x_{2,1} - x_{2,2})} \\
& > \ \ \frac{\exp(- x_{1,1})}{1 + \exp(- x_{1,1})}
\ \ \geq \ \ \frac{\exp(0)}{1 + \exp(0)} \ \ = \ \ \frac{1}{2}
\end{align*}
Thus, \ref{eqn:FP EX3} is infeasible.

The corresponding instance of \ref{eqn:fluid MNL2} is
\begin{align}
\max_{d, \, d_{0}, \, u} \quad & \big\{u_{1,1} + u_{2,1} + 2 u_{2,2}\big\}
\tag{$\mathsf{FP^{Example}_{4}}$}
\label{eqn:FP EX4} \\
\text{s.t.} \quad
& d_{1} + d_{2} \ \ \leq \ \ \frac{1}{2} \nonumber \\
& \left(d_{0}, d_{1}, u_{1,1}\right) \ \ \in \ \ \mathcal{K}_{\exp} \nonumber \\
& \left(d_{0}, d_{2}, u_{2,1} + u_{2,2}\right) \ \ \in \ \ \mathcal{K}_{\exp} \nonumber \\
& \left(d_{1}, d_{0}, 0\right) \ \ \in \ \ \mathcal{K}_{\exp} \nonumber \\
& d_{0} + d_{1} + d_{2} \ \ = \ \ 1 \nonumber \\
& u_{1,1} \ \ \leq \ \ 0 \nonumber
\end{align}
Consider the sequence of solutions $\{(d_{1}^{n}, \, d_{2}^{n}, \, d_{0}^{n}, \, u_{1,1}^{n}, \, u_{2,1}^{n}, \, u_{2,2}^{n}) = (1 / 2, 0, 1 / 2, 0, -n, n)\}_{n=1}^{\infty}$.
The solutions are feasible for \ref{eqn:FP EX4} for all~$n$, and the objective value of $(d_{1}^{n}, \, d_{2}^{n}, \, d_{0}^{n}, \, u_{1,1}^{n}, \, u_{2,1}^{n}, \, u_{2,2}^{n})$ is $n$.
Thus, \ref{eqn:FP EX4} is unbounded.
\end{example}

%------------------------------------------------------------------------

\emph{\underline{The MC model subsumes the MNL model when attributes are lower bounded}:} Assume that all attributes are lower bounded. Consider a generic MNL model
\begin{align*}
\check{P}_{jt}(y) \ \ = \ \ 
		\frac{\exp\left(\alpha_{jt} - \sum_{k \in \mathcal{K}_{jt}} \beta_{kt} y_{jkt}\right)}
		{1 + \sum_{j' \in \mathcal{J}} 
		\exp\left(\alpha_{j't} - \sum_{k \in \mathcal{K}_{j't}} \beta_{kt} y_{j'kt}\right)}
\end{align*}
where
\begin{align*}
\underline{y}_{jkt} \ \ \leq \ \ y_{jkt}
    &\quad , \quad \forall \ j \in \mathcal{J}, \ k \in \mathcal{K}_{jt}
\end{align*}
We let
\begin{align*}
\theta_{jt} \ \ &= \ \ \frac{1}{\exp\left( \sum_{k \in \mathcal{K}_{jt}} 
		\beta_{kt} \underline{y}_{jkt} - \alpha_{jt} \right) + |\mathcal{J}|}
    	\quad , \quad j \in \mathcal{J} \\
\rho_{ijt} \ \ &= \ \ \frac{1}{\exp\left( \sum_{k \in \mathcal{K}_{jt}} 
		\beta_{kt} \underline{y}_{jkt} - \alpha_{jt} \right) + |\mathcal{J}|}
    	\quad , \quad j \in \mathcal{J} 
\end{align*}
let
\begin{align*}
\check{Q}_{jt}(y_{t}) \ \ &= \ \ \exp\left( \sum_{k \in \mathcal{K}_{jt}} \beta_{kt} 
		\left( \underline{y}_{jkt} - y_{jk} \right) \right)
    	\quad , \quad j \in \mathcal{J} 
\end{align*}
and let
\begin{align*}
\check{V}_{jt}(y_{t}) \ \ &= \ \ \frac{1}
		{\exp\left( \sum_{k \in \mathcal{K}_{jt}} 
		\beta_{kt} \underline{y}_{jkt} - \alpha_{jt} \right) + 
		\sum_{j' \in \mathcal{J}} \exp\left( \sum_{k \in \mathcal{K}_{j't}} \beta_{kt} 
		\left( \underline{y}_{j'kt} - y_{j'k} \right) \right)}
    	\quad , \quad j \in \mathcal{J} 
\end{align*}
Now, we have
\begin{align*}
\check{P}_{jt}(y) = \check{Q}_{jt}(y_{t}) \check{V}_{jt}(y_{t})
    	\quad , \quad j \in \mathcal{J} 
\end{align*}
Meanwhile, since
\begin{align*}
& 1 + \sum_{i \in \mathcal{J}_{t}} \left[1 - \check{Q}_{it}(y_{t})\right]\check{V}_{it}(y_{t}) \\
& \ \ = \ \ 1 + \sum_{i \in \mathcal{J}_{t}} 
		\frac{1 - \exp\left( \sum_{k \in \mathcal{K}_{i}} \beta_{kt} 
		\left( \underline{y}_{ikt} - y_{ik} \right) \right)}
		{\exp\left( \sum_{k \in \mathcal{K}_{it}} 
		\beta_{kt} \underline{y}_{ikt} - \alpha_{it} \right) + 
		\sum_{i' \in \mathcal{J}} \exp\left( \sum_{k \in \mathcal{K}_{i'}} \beta_{kt} 
		\left( \underline{y}_{i'kt} - y_{i'k} \right) \right)} \\
& \ \ = \ \ 1 + \sum_{i \in \mathcal{J}_{t}} 
		\frac{\exp\left( \alpha_{it} - 
				\sum_{k \in \mathcal{K}_{it}} \beta_{kt} \underline{y}_{ikt}\right) - 
		\exp\left( \alpha_{i't} - 
				\sum_{k \in \mathcal{K}_{it}} \beta_{kt} y_{ik} \right)}
		{1 + \sum_{i' \in \mathcal{J}} \exp\left( \alpha_{it} - 
				\sum_{k \in \mathcal{K}_{i't}} \beta_{kt} y_{i'k} \right)} \\
& \ \ = \ \ 
		\frac{1 + \sum_{i \in \mathcal{J}_{t}} \exp\left( \alpha_{it} - 
				\sum_{k \in \mathcal{K}_{it}} \beta_{kt} \underline{y}_{ikt}\right)}
		{1 + \sum_{i' \in \mathcal{J}_{t}} \exp\left( \alpha_{it} - 
				\sum_{k \in \mathcal{K}_{i't}} \beta_{kt} y_{i'k} \right)} \\
& \ \ = \ \ 
		\frac{\exp\left( \sum_{k \in \mathcal{K}_{jt}} 
		\beta_{kt} \underline{y}_{jkt} - \alpha_{jt} \right) + |\mathcal{J}|}
		{\exp\left( \sum_{k \in \mathcal{K}_{jt}} 
		\beta_{kt} \underline{y}_{jkt} - \alpha_{jt} \right) + 
		\sum_{j' \in \mathcal{J}} \exp\left( \sum_{k \in \mathcal{K}_{j't}} \beta_{kt} 
		\left( \underline{y}_{j'kt} - y_{j'k} \right) \right)}
\end{align*}
we know that $\check{V}(y)$ is a solution to the system \ref{eqn:sys-mc}. 

It left to prove that the $\check{V}(y)$ is the unique solution to \ref{eqn:sys-mc}. Indeed. it is clear that $\theta_{t} \geq 0$, $\sum_{j \in \mathcal{J}} \theta_{jt} > 0$, $\rho_{t} \geq 0$ and $\sum_{n=1}^{\infty} \rho_{t}^{n} < \infty$. In addition, since $\underline{y}_{jkt} \leq y_{jkt}$ for any $j \in \mathcal{J}$ and $k \in \mathcal{K}_{jt}$, we have $\check{Q}_{j}(y) \in [0,1]$. Thus, all necessary assumptions hold, and the solution is unique.

%------------------------------------------------------------------------
%	Section
%------------------------------------------------------------------------

\section{Conclusion}
\label{sec:conclusion}

In this paper, we formulate the profit maximization problems \ref{eqn:static general} and \ref{eqn:fluid general} with multiple product attributes, and provide their convex conic reformulations under three choice models. When there is only one attribute associating with each product, \ref{eqn:static general} and \ref{eqn:fluid general} reduces to the classical pricing problems and revenue management problems with resource constraints. Thus, our paper reproduces the tractability results of pricing problems under the MNL model (\cite{song2007demand}, \cite{dong2009dynamic}), the MC model (\cite{dong2019pricing}), and the NL model (\cite{li2011pricing}). For the NL model, our approach also covers the case when the nest dissimilarity parameters $\gamma$ is larger than $1$. Tractability result of this scenario was given in \cite{gallego2014constrained}. However, their method only words for unconstrained problems, while our approach allows different type of constraints, such as resource constraints and upper and lower bound of attributes. 

While we kept our discussion minimum, our approach can be extended to work with multiple types of objective functions and constraints. For example, revenue maximization can be viewed as a special case of profit maximization with zero costs, while market share maximization  can be viewed as a special case of profit maximization with zero profit margins. In general, any objective function is allowed as long as (1) it can be transformed into a concave function in the convex conic reformulation, and (2) the nonlinear inequality constraints hold tight at optimality. Similarly, we can work with different types of constraints, as long as there is a good convex transformation. Another possible direction of extension is allowing utility and / or profit margin to be nonlinear functions of attributes. One example is to replace an attribute $x_{j}$ by $x_{j}^{+} - x_{j}^{-}$, where $x_{j}^{+}, x_{j}^{-} \geq 0$. Then $x_{j}^{+}$ and $x_{j}^{-}$ can have different sensitivity parameters $\beta_{j}^{+}$ and $\beta_{j}^{-}$. (Note that $\beta_{j}^{+}$ and $\beta_{j}^{-}$ need to satisfy certain conditions, such that one cannot benefit from increasing $x_{j}^{+}$ and $x_{j}^{-}$ simultaneously.) Thus, the utility is a piecewise linear function of $x_{j}$.

\newpage

% ------------------------------------------------------------
% REFERENCE
% ------------------------------------------------------------

\bibliography{library}
\bibliographystyle{plainnat}

 \begin{appendix}

\newpage

%------------------------------------------------------------------------
%	Section
%------------------------------------------------------------------------

\section{Proofs}

%------------------------------------------------------------------------

\subsection{Proof of Lemma~\ref{lem:spmnl5 optimal}}

{\color{purple}
\begin{manuallemma}{\ref{lem:spmnl5 optimal}}
\ref{eqn:static MNL5} has an optimal solution if and only if for each $j \in \mathcal{J}$, and for all $k_{1} \in \mathcal{K}_{j} \setminus \overline{\mathcal{K}}_{j}$, $k_{2} \in \mathcal{K}_{j} \setminus \underline{\mathcal{K}}_{j}$, it holds that $\phi_{k_{1}} \le \phi_{k_{2}}$; otherwise, \ref{eqn:static MNL5} is unbounded.
\end{manuallemma}
}

\begin{proof}
First, note that \ref{eqn:static MNL5} is feasible.
For example, consider $(d^{0}, d^{0}_{0}, u^{0})$ given as follows:
For all $j \in \mathcal{J}$ and $k \in \overline{\mathcal{K}}_{j}$, let $x^{0}_{jk} = \overline{x}_{jk}$.
For all $j \in \mathcal{J}$ and $k \in \underline{\mathcal{K}}_{j} \setminus \overline{\mathcal{K}}_{j}$, let $x^{0}_{jk} = \underline{x}_{jk}$.
For all $j \in \mathcal{J}$ and $k \in \mathcal{K}_{j} \setminus \left(\underline{\mathcal{K}}_{j} \cup \overline{\mathcal{K}}_{j}\right)$, let $x^{0}_{jk} = 0$.
Next, let $d^{0}_{0}$ be given by~\eqref{eqn:d0}, let $d^{0}_{j}$ be given by~\eqref{eqn:spmnl1 dj} for all $j \in \mathcal{J}$, and let $u^{0}_{jk} = d^{0}_{j} x^{0}_{jk}$ for all $j \in \mathcal{J}$ and $k \in \mathcal{K}_{j}$.
Then $(d^{0}, d^{0}_{0}, u^{0})$ is feasible for \ref{eqn:static MNL5}.

Suppose that there is a $j \in \mathcal{J}$, a $k_{1} \in \mathcal{K}_{j} \setminus \overline{\mathcal{K}}_{j}$, and a $k_{2} \in \mathcal{K}_{j} \setminus \underline{\mathcal{K}}_{j}$, such that $\phi_{k_{1}} > \phi_{k_{2}}$.
Then consider a sequence $\{(d^{n}, d^{n}_{0}, u^{n})\}_{n=0}^{\infty}$ of feasible solutions with $u^{n}_{jk_{1}} = u^{0}_{jk_{1}} + n$, $u^{n}_{jk_{2}} = u^{0}_{jk_{2}} - n$, $u^{n}_{jk} = u^{0}_{jk}$ for all $k \notin \{k_{1}, k_{2}\}$, $u^{n}_{j'k} = u^{0}_{j'k}$ for all $j' \neq j$ and $k \in \mathcal{K}_{j'}$, $d^{n}_{j'} = d^{0}_{j'}$ for all $j' \in \mathcal{J}$, and $d^{n}_{0} = d^{0}_{0}$.
Then $(d^{n}, d^{n}_{0}, u^{n})$ is feasible for \ref{eqn:static MNL5} for all~$n$, and the objective value of $(d^{n+1}, d^{n+1}_{0}, u^{n+1})$ exceeds the objective value of $(d^{n}, d^{n}_{0}, u^{n})$ by $\phi_{k_{1}} - \phi_{k_{2}}$ for all~$n$, and thus the objective value of $(d^{n}, d^{n}_{0}, u^{n})$ increases without bound as $n \to \infty$.

Next suppose that for each $j \in \mathcal{J}$, and for all $k_{1} \in \mathcal{K}_{j} \setminus \overline{\mathcal{K}}_{j}$, $k_{2} \in \mathcal{K}_{j} \setminus \underline{\mathcal{K}}_{j}$, it holds that $\phi_{k_{1}} \le \phi_{k_{2}}$.
Let
\begin{align*}
\mathcal{K}_{\exp}^* \ \ & \defi \ \ \mbox{closure}\big\{(y_{1}, y_{2}, y_{3}) \; : \; y_{1} \ge -y_{3} \exp(y_{2}/y_{3} - 1), \; y_{1} > 0, \; y_{3} < 0\big\} \\& = \ \ \big\{(y_{1}, y_{2}, y_{3}) \; : \; y_{1} \ge -y_{3} \exp(y_{2}/y_{3} - 1), \; y_{1} > 0, \; y_{3} < 0\big\} \cup \big\{(y_{1}, y_{2}, 0) \; : \; y_{1} \ge 0, \; y_{2} \ge 0\big\}
\end{align*}
denote the dual cone of $\mathcal{K}_{\exp}$.
Let $\pi_{j} \in \mathcal{K}_{\exp}^*$ denote the dual variables corresponding to~\eqref{mp5-mnl}, let $\varpi_{j} \in \mathcal{K}_{\exp}^*$ denote the dual variables corresponding to~\eqref{mp5-total upper bound}, let $\underline{\nu}_{jk} \ge 0$ denote the dual variables corresponding to~\eqref{eqn:spmnl5-lower-bound}, let $\overline{\nu}_{jk} \ge 0$ denote the dual variables corresponding to~\eqref{eqn:spmnl5-upper-bound}, and let $\gamma$ denote the dual variable corresponding to~\eqref{mp5-totalmarket}.
Then the dual of \ref{eqn:static MNL5} is
{\small
\begin{subequations}
\begin{align}
\min_{\pi, \, \varpi, \, \underline{\nu}, \, \overline{\nu}, \, \gamma} \quad & \gamma
\tag{$\mathsf{SD^{MNL}_{5}}$}
\label{eqn:static MNL5D} \\
\text{s.t.} \quad
& \pi_{j3} \ \ = \ \ -\phi_{k}
& \forall \ j \in \mathcal{J}, \ k \in \mathcal{K}_{j} \setminus \left(\underline{\mathcal{K}}_{j} \cup \overline{\mathcal{K}}_{j}\right) \label{md5-u1} \\
& \pi_{j3} + \underline{\nu}_{jk} \ \ = \ \ -\phi_{k}
& \forall \ j \in \mathcal{J}, \ k \in \underline{\mathcal{K}}_{j} \setminus \overline{\mathcal{K}}_{j} \label{md5-u2} \\
& \pi_{j3} - \overline{\nu}_{jk} \ \ = \ \ -\phi_{k}
& \forall \ j \in \mathcal{J}, \ k \in \overline{\mathcal{K}}_{j} \setminus \underline{\mathcal{K}}_{j} \label{md5-u3} \\
& \pi_{j3} + \underline{\nu}_{jk} - \overline{\nu}_{jk} \ \ = \ \ -\phi_{k}
& \forall \ j \in \mathcal{J}, \ k \in \underline{\mathcal{K}}_{j} \cap \overline{\mathcal{K}}_{j} \label{md5-u4} \\
& \pi_{j2} + \varpi_{j1} - \sum_{k \in \underline{\mathcal{K}}_{j}} \underline{\nu}_{jk} \underline{x}_{jk} + \sum_{k \in \overline{\mathcal{K}}_{j}} \overline{\nu}_{jk} \overline{x}_{jk} - \gamma \ \ = \ \ \psi_{j}
& \forall \ j \in \overline{\mathcal{J}} \label{md5-dj1} \\
& \pi_{j2} - \sum_{k \in \underline{\mathcal{K}}_{j}} \underline{\nu}_{jk} \underline{x}_{jk} + \sum_{k \in \overline{\mathcal{K}}_{j}} \overline{\nu}_{jk} \overline{x}_{jk} - \gamma \ \ = \ \ \psi_{j}
& \forall \ j \in \mathcal{J} \setminus \overline{\mathcal{J}} \label{md5-dj2} \\
& \sum_{j \in \mathcal{J}} \pi_{j1} + \sum_{j \in \overline{\mathcal{J}}} \varpi_{j2} - \sum_{j \in \overline{\mathcal{J}}} \varpi_{j3} \sum_{k \in \mathcal{K}_{j}} \overline{x}_{jk} - \gamma \ \ = \ \ 0 \label{md5-d0} \\
& \pi_{j} \ \ \in \ \ \mathcal{K}_{\exp}^*
& \forall \ j \in \mathcal{J} \label{md5-j1} \\
& \varpi_{j} \ \ \in \ \ \mathcal{K}_{\exp}^*
& \forall \ j \in \overline{\mathcal{J}} \label{md5-j2} \\
& \underline{\nu}_{jk} \ \ \ge \ \ 0
& \forall \ j \in \mathcal{J}, \ k \in \underline{\mathcal{K}}_{j} \label{md5-lower bound sign} \\
& \overline{\nu}_{jk} \ \ \ge \ \ 0
& \forall \ j \in \mathcal{J}, \ k \in \overline{\mathcal{K}}_{j} \label{md5-upper bound sign}
\end{align}
\end{subequations}
}%

Next we show that there is a feasible solution $(\pi, \varpi, \underline{\nu}, \overline{\nu}, \gamma)$ of~\ref{eqn:static MNL5D} such that $\pi_{j} \in \mbox{int}\left(\mathcal{K}_{\exp}^*\right)$ for all~$j \in \mathcal{J}$ and $\varpi_{j} \in \mbox{int}\left(\mathcal{K}_{\exp}^*\right)$ for all~$j \in \overline{\mathcal{J}}$.
Choose $\varpi_{j2} = 0$, $\varpi_{j3} = -1$, and $\varpi_{j1} = \exp(-1) + 1$ for all~$j \in \overline{\mathcal{J}}$.
Thus $\varpi_{j} \in \mbox{int}\left(\mathcal{K}_{\exp}^*\right)$ for all~$j \in \overline{\mathcal{J}}$.
For each $j \in \mathcal{J}$, choose $\pi_{j3}$, $\underline{\nu}_{jk}$ for all $k \in \underline{\mathcal{K}}_{j}$, and $\overline{\nu}_{jk}$ for all $k \in \overline{\mathcal{K}}_{j}$, by considering the following cases: 
\begin{enumerate}
\item
If $\mathcal{K}_{j} \setminus \underline{\mathcal{K}}_{j} \neq \varnothing$, then choose $\pi_{j3} = - \min\left\{\phi_{k} \, : \, k \in \mathcal{K}_{j} \setminus \underline{\mathcal{K}}_{j}\right\}$.
Since $\phi_{k} > 0$ for all~$k$, it follows that $\pi_{j3} < 0$.
Then, for each $k \in \mathcal{K}_{j}$, consider the following four cases.
\begin{enumerate}
\item
Suppose $k \in \mathcal{K}_{j} \setminus \left(\underline{\mathcal{K}}_{j} \cup \overline{\mathcal{K}}_{j}\right)$.
Since $k \in \mathcal{K}_{j} \setminus \underline{\mathcal{K}}_{j}$, it follows that $- \pi_{j3} \le \phi_{k}$.
Since $k \in \mathcal{K}_{j} \setminus \overline{\mathcal{K}}_{j}$ and $\phi_{k_{1}} \le \phi_{k_{2}}$ for all $k_{1} \in \mathcal{K}_{j} \setminus \overline{\mathcal{K}}_{j}$, $k_{2} \in \mathcal{K}_{j} \setminus \underline{\mathcal{K}}_{j}$, it follows that $\phi_{k} \le - \pi_{j3}$.
Thus, $\phi_{k} = - \pi_{j3}$, and hence constraint~\eqref{md5-u1} is satisfied.
\item
Suppose $k \in \underline{\mathcal{K}}_{j} \setminus \overline{\mathcal{K}}_{j}$.
Then choose $\underline{\nu}_{jk} = - \pi_{j3} - \phi_{k}$.
Hence constraint~\eqref{md5-u2} is satisfied.
Since $k \in \mathcal{K}_{j} \setminus \overline{\mathcal{K}}_{j}$ and $\phi_{k_{1}} \le \phi_{k_{2}}$ for all $k_{1} \in \mathcal{K}_{j} \setminus \overline{\mathcal{K}}_{j}$, $k_{2} \in \mathcal{K}_{j} \setminus \underline{\mathcal{K}}_{j}$, it follows that $\phi_{k} \le - \pi_{j3}$.
Thus $\underline{\nu}_{jk} = - \pi_{j3} - \phi_{k} \ge 0$, and hence constraint~\eqref{md5-lower bound sign} is satisfied.
\item
Suppose $k \in \overline{\mathcal{K}}_{j} \setminus \underline{\mathcal{K}}_{j}$.
Then choose $\overline{\nu}_{jk} = \pi_{j3} + \phi_{k}$.
Hence constraint~\eqref{md5-u3} is satisfied.
Since $k \in \mathcal{K}_{j} \setminus \underline{\mathcal{K}}_{j}$, it follows that $- \pi_{j3} \le \phi_{k}$.
Thus $\overline{\nu}_{jk} = \pi_{j3} + \phi_{k} \ge 0$, and hence constraint~\eqref{md5-upper bound sign} is satisfied.
\item
Suppose $k \in \underline{\mathcal{K}}_{j} \cap \overline{\mathcal{K}}_{j}$.
Then choose $\underline{\nu}_{jk} = \max\left\{0, \, - \pi_{j3} - \phi_{k}\right\}$ and $\overline{\nu}_{jk} = \max\left\{0, \, \pi_{j3} + \phi_{k}\right\}$.
Hence constraints~\eqref{md5-u4}, \eqref{md5-lower bound sign} and~\eqref{md5-upper bound sign} are satisfied.
\end{enumerate}
\item
If $\mathcal{K}_{j} \setminus \underline{\mathcal{K}}_{j} = \varnothing$ (that is, $\mathcal{K}_{j} = \underline{\mathcal{K}}_{j}$), then choose $\pi_{j3} = - \max\left\{\phi_{k} \, : \, k \in \mathcal{K}_{j}\right\}$.
Since $\phi_{k} > 0$ for all~$k$, it follows that $\pi_{j3} < 0$.
Then, for each $k \in \mathcal{K}_{j}$, consider the following two cases (only these two cases are possible).
\begin{enumerate}
\item
Suppose $k \in \underline{\mathcal{K}}_{j} \setminus \overline{\mathcal{K}}_{j}$.
Then choose $\underline{\nu}_{jk} = - \pi_{j3} - \phi_{k}$.
Hence constraint~\eqref{md5-u2} is satisfied.
Since $\phi_{k} \le - \pi_{j3}$, it follows that $\underline{\nu}_{jk} = - \pi_{j3} - \phi_{k} \ge 0$, and hence constraint~\eqref{md5-lower bound sign} is satisfied.
\item
Suppose $k \in \underline{\mathcal{K}}_{j} \cap \overline{\mathcal{K}}_{j}$.
Then choose $\underline{\nu}_{jk} = \max\left\{0, \, - \pi_{j3} - \phi_{k}\right\}$ and $\overline{\nu}_{jk} = \max\left\{0, \, \pi_{j3} + \phi_{k}\right\}$.
Hence constraints~\eqref{md5-u4}, \eqref{md5-lower bound sign} and~\eqref{md5-upper bound sign} are satisfied.
\end{enumerate}
\end{enumerate}
Thus $\pi_{j3}$, $\underline{\nu}_{jk}$, and $\overline{\nu}_{jk}$ have been chosen for all~$j$ and~$k$.

Next we use~\eqref{md5-dj1} and~\eqref{md5-dj2} to express $\pi_{j2}$ as a function of $\gamma$ for each $j \in \mathcal{J}$.
Then we choose $\pi_{j1}$ as a function of $\gamma$ such that $\pi_{j1} = - \pi_{j3} \exp(\pi_{j2} / \pi_{j3} - 1) + 1$.
Then it will follow that $\pi_{j} \in \mbox{int}\left(\mathcal{K}_{\exp}^*\right)$ for all~$j$.
Finally we will show that there exists a $\gamma$ such that~\eqref{md5-d0} is satisfied, and thus all the constraints of~\ref{eqn:static MNL5D} will be satisfied.

If $j \in \overline{\mathcal{J}}$, then let $e_{j} = - \varpi_{j1} + \sum_{k \in \underline{\mathcal{K}}_{j}} \underline{\nu}_{jk} \underline{x}_{jk} - \sum_{k \in \overline{\mathcal{K}}_{j}} \overline{\nu}_{jk} \overline{x}_{jk} + \psi_{j}$, and if $j \in \mathcal{J} \setminus \overline{\mathcal{J}}$, then let $e_{j} = \sum_{k \in \underline{\mathcal{K}}_{j}} \underline{\nu}_{jk} \underline{x}_{jk} - \sum_{k \in \overline{\mathcal{K}}_{j}} \overline{\nu}_{jk} \overline{x}_{jk} + \psi_{j}$.
For all $j \in \mathcal{J}$, let $\pi_{j2} = \pi_{j2}(\gamma) = \gamma + e_{j}$, and let $\pi_{j1} = \pi_{j1}(\gamma) = - \pi_{j3} \exp\left(\pi_{j2}(\gamma) / \pi_{j3} - 1\right) + 1 = - \pi_{j3} \exp\left(\left[\gamma + e_{j}\right] / \pi_{j3} - 1\right) + 1$.
It follows that for any $\gamma$, \eqref{md5-dj1} is satisfied for all $j \in \overline{\mathcal{J}}$ and~\eqref{md5-dj2} is satisfied for all $j \in \mathcal{J} \setminus \overline{\mathcal{J}}$.
It also follows that for any $\gamma$, $\pi_{j} \in \mbox{int}\left(\mathcal{K}_{\exp}^*\right)$ for all~$j$.
For constraint~\eqref{md5-d0} to be satisfied, $\gamma$ has to satisfy
$\sum_{j \in \mathcal{J}} \pi_{j1}(\gamma) + \sum_{j \in \overline{\mathcal{J}}} \varpi_{j2} - \sum_{j \in \overline{\mathcal{J}}} \varpi_{j3} \sum_{k \in \mathcal{K}_{j}} \overline{x}_{jk} - \gamma = 0$, that is,
\[
h(\gamma) \ \ \defi \ \ \sum_{j \in \mathcal{J}} \left\{- \pi_{j3} \exp\left(\frac{\gamma + e_{j}}{\pi_{j3}} - 1\right) + 1\right\} + \sum_{j \in \overline{\mathcal{J}}} \varpi_{j2} - \sum_{j \in \overline{\mathcal{J}}} \varpi_{j3} \sum_{k \in \mathcal{K}_{j}} \overline{x}_{jk} - \gamma \ \ = \ \ 0
\]
Since $\pi_{j3} < 0$ for all~$j$, it follows that $h$ is decreasing, $h(\gamma) \to \infty$ as $\gamma \to -\infty$, and $h(\gamma) \to -\infty$ as $\gamma \to \infty$.
Thus it follows from the intermediate value theorem that there is a (unique) $\gamma^*$ such that $h(\gamma^*) = 0$, and hence $(\gamma^*, \, \pi_{j1}(\gamma^*), j \in \mathcal{J})$ satisfies~\eqref{md5-d0}.
The resulting solution $(\pi(\gamma^*), \varpi, \underline{\nu}, \overline{\nu}, \gamma^*)$ is feasible for~\ref{eqn:static MNL5D} and satisfies $\pi_{j}(\gamma^*) \in \mbox{int}\left(\mathcal{K}_{\exp}^*\right)$ for all~$j \in \mathcal{J}$ and $\varpi_{j} \in \mbox{int}\left(\mathcal{K}_{\exp}^*\right)$ for all~$j \in \overline{\mathcal{J}}$.
Therefore it follows from the conic duality theorem (Theorem~1.4.2 in \cite{ben2001lectures}) that \ref{eqn:static MNL5} has an optimal solution.
\end{proof}

\vspace{5mm}

\newpage

%------------------------------------------------------------------------

\subsection{Proof of Lemma~\ref{lem:spmnl4=spmnl5}}

{\color{purple}
\begin{manuallemma}{\ref{lem:spmnl4=spmnl5}}
Every optimal solution for~\ref{eqn:static MNL5} is also optimal for~\ref{eqn:static MNL4}. 
\end{manuallemma}
}

\begin{proof}
Consider any $(d, d_{0}, u)$ feasible for~\ref{eqn:static MNL4}.
Then $(d, d_{0}, u)$ is also feasible for~\ref{eqn:static MNL5} and has the same objective value.
Next, consider any $(d, d_{0}, u)$ feasible for~\ref{eqn:static MNL5}.
Note that if $d > 0$ and $d_{0} > 0$, then $(d, d_{0}, u)$ is feasible for~\ref{eqn:static MNL4} and has the same objective value.
We show by contradiction that $d_{0} > 0$, and that $d_{j} > 0$ for any $j \in \overline{\mathcal{J}}$.
Suppose that $d_{0} = 0$.
Then it follows from constraint~\eqref{mp5-mnl} that $d_{j} = 0$ for all~$j$, which would violate constraint~\eqref{mp5-totalmarket}.
Similarly, suppose that $d_{j} = 0$ for some~$j \in \overline{\mathcal{J}}$, then it follows from constraint~\eqref{mp5-total upper bound} that $d_{0} = 0$, which as shown above cannot be feasible for~\ref{eqn:static MNL5}.
However, $d_{j} = 0$ for $j \in \mathcal{J} \setminus \overline{\mathcal{J}}$ can be feasible for~\ref{eqn:static MNL5}.

Next we show that if $(d^*, d_{0}^*, u^*)$ is an optimal solution for \ref{eqn:static MNL5}, then it is also optimal for \ref{eqn:static MNL4}.
It remains to show that $d_{j}^* > 0$ for every $j \in \mathcal{J} \setminus \overline{\mathcal{J}}$.
%As shown above, we then have $d_{0}^* > 0, d^* > 0$, which means $(d^*, d_{0}^*, u^*)$ is feasible to \ref{eqn:static MNL4}, and has the same objective value.
%Since any feasible solution to \ref{eqn:static MNL4} is also feasible to \ref{eqn:static MNL5}, and has the same objective value, we can conclude that $(d^*, d_{0}^*, u^*)$ is an optimal solution to \ref{eqn:static MNL4}.
%We show that $d_{j}^* > 0$ for every $j \in \mathcal{J} \setminus \overline{\mathcal{J}}$ by contradiction.
Suppose that $d_{i}^* = 0$ for some $i \in \mathcal{J} \setminus \overline{\mathcal{J}}$.
%Since $i \in \mathcal{J} \setminus \overline{\mathcal{J}}$, we have $\sum_{k \in \mathcal{K}_{i}} u_{ik}^* = 0$.
%(Otherwise, there must be a better feasible solution.)
Let $\underline{\mathcal{K}}_{i}^{+} \defi \{k \in \underline{\mathcal{K}}_{i} \, : \, \underline{x}_{ik} > 0\}$ and $\overline{\mathcal{K}}_{i}^{-} \defi \{k \in \overline{\mathcal{K}}_{i} \, : \, \overline{x}_{ik} < 0\}$.
Consider any $\delta \in \big(0, \; \min\big\{1, \; \min\{u_{ik}^* / \underline{x}_{ik} \, : \, k \in \underline{\mathcal{K}}_{i}^{+}, u_{ik}^* > 0\}, \; \min\{u_{ik}^* / \overline{x}_{ik} \, : \, k \in \overline{\mathcal{K}}_{i}^{-}, u_{ik}^* < 0\}\big\}\big)$.
Consider the solution $(d', d_{0}', u')$ for~\ref{eqn:static MNL5} given by
{\small
\begin{align*}
d'_{i} \ \ & = \ \ \delta \\
d_{0}' \ \ & = \ \ d_{0}^* \left(1 - \delta\right) \\
d'_{j} \ \ & = \ \ d_{j}^* \left(1 - \delta\right) \hspace{110mm} \mbox{for all } j \in \mathcal{J} \setminus \{i\} \\
u'_{il} \ \ & = \ \ u_{il}^* \; - \sum_{\left\{k \in \underline{\mathcal{K}}_{i}^{+} \setminus \{l\} \, : \, u_{ik}^* = 0\right\}} \underline{x}_{ik} \delta \; - \sum_{\left\{k \in \overline{\mathcal{K}}_{i}^{-} \, : \, u_{ik}^* = 0\right\}} \overline{x}_{ik} \delta \; - \; \delta \ln\left(\frac{\delta}{d_{0}^* \left(1 - \delta\right)}\right) \hspace{15mm} \mbox{for some } l \in \mathcal{K}_{i} \setminus \overline{\mathcal{K}}_{i} \\
u'_{ik} \ \ & = \ \ u_{ik}^* + \underline{x}_{ik} \delta \ \ = \ \ \underline{x}_{ik} \delta \hspace{75mm} \mbox{for all } k \in \underline{\mathcal{K}}_{i}^{+} \setminus \{l\} \, : \, u_{ik}^* = 0 \\
u'_{ik} \ \ & = \ \ u_{ik}^* + \overline{x}_{ik} \delta \ \ = \ \ \overline{x}_{ik} \delta \hspace{83mm} \mbox{for all } k \in \overline{\mathcal{K}}_{i}^{-} \, : \, u_{ik}^* = 0 \\
u'_{ik} \ \ & = \ \ u_{ik}^* \hspace{50mm} \mbox{for all } k \in \mathcal{K}_{i} \setminus \big(\{l\} \cup \left\{k \in \underline{\mathcal{K}}_{i}^{+} \, : \, u_{ik}^* = 0\right\} \cup \big\{k \in \overline{\mathcal{K}}_{i}^{-} \, : \, u_{ik}^* = 0\big\}\big) \\
u'_{jk} \ \ & = \ \ u_{jk}^* \left(1 - \delta\right) \hspace{95mm} \mbox{for all } j \in \mathcal{J} \setminus \{i\}, \ k \in \mathcal{K}_{j}
\end{align*}
}
Note that, since $i \in \mathcal{J} \setminus \overline{\mathcal{J}}$, it holds that $\mathcal{K}_{i} \setminus \overline{\mathcal{K}}_{i} \neq \varnothing$, and thus an attribute $l \in \mathcal{K}_{i} \setminus \overline{\mathcal{K}}_{i}$ exists.

Next we show that $(d', d_{0}', u')$ is feasible for~\ref{eqn:static MNL5}.
\begin{itemize}
\item
Note that for all $j \in \mathcal{J} \setminus \{i\}$, it holds that $d_{j}' / d_{0}' =  d_{j}^* / d_{0}^*$, and thus $(d', d_{0}', u')$ satisfies~\eqref{mp5-total upper bound}.
\item
Next we show that $(d', d_{0}', u')$ satisfies~\eqref{mp5-mnl}.
\begin{itemize}
\item
By~\eqref{mp5-mnl}, $\sum_{k \in \mathcal{K}_{i}} u_{ik}^* \leq 0$, and thus
\begin{align*}
\sum_{k \in \mathcal{K}_{i}} u'_{ik}
\ \ & = \ \ u_{il}^* \; - \sum_{\left\{k \in \underline{\mathcal{K}}_{i}^{+} \setminus \{l\} \, : \, u_{ik}^* = 0\right\}} \underline{x}_{ik} \delta \; - \sum_{\left\{k \in \overline{\mathcal{K}}_{i}^{-} \, : \, u_{ik}^* = 0\right\}} \overline{x}_{ik} \delta \; - \; \delta \ln\left(\frac{\delta}{d_{0}^* \left(1 - \delta\right)}\right) \\
& \qquad + \sum_{\left\{k \in \underline{\mathcal{K}}_{i}^{+} \setminus \{l\} \, : \, u_{ik}^* = 0\right\}} \big(u_{ik}^* + \underline{x}_{ik} \delta\big) \; + \sum_{\left\{k \in \overline{\mathcal{K}}_{i}^{-} \, : \, u_{ik}^* = 0\right\}} \big(u_{ik}^* + \overline{x}_{ik} \delta\big) \\
& \qquad + \sum_{k \in \mathcal{K}_{i} \setminus \big(\{l\} \cup \left\{k \in \underline{\mathcal{K}}_{i}^{+} \, : \, u_{ik}^* = 0\right\} \cup \big\{k \in \overline{\mathcal{K}}_{i}^{-} \, : \, u_{ik}^* = 0\big\}\big)} u_{ik}^* \\
& = \ \ \sum_{k \in \mathcal{K}_{i}} u_{ik}^* - \delta \ln\left(\frac{\delta}{d_{0}^* \left(1 - \delta\right)}\right)
\ \ \le \ \ - d'_{i} \ln\left(\frac{d'_{i}}{d_{0}'}\right)
\end{align*}
and hence $(d', d_{0}', u')$ satisfies~\eqref{mp5-mnl} for $j = i$.
\item
For $j \in \mathcal{J} \setminus \{i\}$, if $d_{j}^* > 0$ then it holds that
\[
\sum_{k \in \mathcal{K}_{j}} u'_{jk}
\ \ = \ \ \left(1 - \delta\right) \sum_{k \in \mathcal{K}_{j}} u_{jk}^*
\ \ \le \ \ - \left(1 - \delta\right) d_{j}^* \ln\left(\frac{d_{j}^*}{d_{0}^*}\right)
\ \ = \ \ - d'_{j} \ln\left(\frac{d'_{j}}{d_{0}'}\right)
\]
and if $d_{j}^* = 0$ then it holds that $d'_{j} = 0$ and
\[
\sum_{k \in \mathcal{K}_{j}} u'_{jk}
\ \ = \ \ \left(1 - \delta\right) \sum_{k \in \mathcal{K}_{j}} u_{jk}^*
\ \ \le \ \ 0
\]
and hence $(d', d_{0}', u')$ satisfies~\eqref{mp5-mnl} for all $j \in \mathcal{J}$.
\end{itemize}
\item
Next, we show that $(d', d_{0}', u')$ satisfies~\eqref{eqn:spmnl5-lower-bound}.
\begin{itemize}
\item
Suppose that $l \in \underline{\mathcal{K}}_{i}$.
Let $f(\delta) \defi - \sum_{\left\{k \in \underline{\mathcal{K}}_{i}^{+} \setminus \{l\} \, : \, u_{ik}^* = 0\right\}} \underline{x}_{ik} \delta \, - \sum_{\left\{k \in \overline{\mathcal{K}}_{i}^{-} \, : \, u_{ik}^* = 0\right\}} \overline{x}_{ik} \delta \, - \, \delta \ln\left(\delta / \left[d_{0}^* \left(1 - \delta\right)\right]\right) - \underline{x}_{il} \delta$.
Note that $f'(\delta) = - \sum_{\left\{k \in \underline{\mathcal{K}}_{i}^{+} \setminus \{l\} \, : \, u_{ik}^* = 0\right\}} \underline{x}_{ik} \, - \sum_{\left\{k \in \overline{\mathcal{K}}_{i}^{-} \, : \, u_{ik}^* = 0\right\}} \overline{x}_{ik} \, - \, \ln\left(\delta / \left[d_{0}^* \left(1 - \delta\right)\right]\right) \, - \, 1 \, - \, \delta / \left(1 - \delta\right) - \underline{x}_{il} \to \infty$ as $\delta \downarrow 0$.
Thus, $f'(s) > 0$ for all $s > 0$ sufficiently small, hence $\left[u'_{il} - \underline{x}_{il} d'_{i}\right] - \left[u_{il}^* - \underline{x}_{il} d_{i}^*\right] = \int_{0}^{\delta} f'(s) \diff s > 0$, and therefore $u'_{il} - \underline{x}_{il} d'_{i} > u_{il}^* - \underline{x}_{il} d_{i}^* \ge 0$ for all $\delta > 0$ sufficiently small.
\item
Consider any $k \in \underline{\mathcal{K}}_{i}^{+} \setminus \{l\}$.
If $u_{ik}^* = 0$, then $u'_{ik} = \underline{x}_{ik} \delta = \underline{x}_{ik} d'_{i}$.
Otherwise, $u_{ik}^* > 0$, then $u'_{ik} = u_{ik}^* > \underline{x}_{ik} \delta = \underline{x}_{ik} d'_{i}$.
\item
Consider any $k \in \underline{\mathcal{K}}_{i} \setminus \left(\underline{\mathcal{K}}_{i}^{+} \cup \{l\}\right)$.
Then $u'_{ik} = u_{ik}^* \ge 0 \ge \underline{x}_{ik} \delta = \underline{x}_{ik} d'_{i}$.
\item
Next, consider any $j \in \mathcal{J} \setminus \{i\}$ and any $k \in \underline{\mathcal{K}}_{j}$.
Then $u'_{jk} = u_{jk}^* \left(1 - \delta\right) \ge \underline{x}_{jk} d_{j}^* \left(1 - \delta\right) = \underline{x}_{jk} d'_{j}$.
Thus, $(d', d_{0}', u')$ satisfies~\eqref{eqn:spmnl5-lower-bound} for all $j \in \mathcal{J}$ and $k \in \underline{\mathcal{K}}_{j}$.
\end{itemize}
\item
Next, we show that $(d', d_{0}', u')$ satisfies~\eqref{eqn:spmnl5-upper-bound}.
\begin{itemize}
\item
Consider any $k \in \overline{\mathcal{K}}_{i}^{-}$.
If $u_{ik}^* = 0$, then $u'_{ik} = \overline{x}_{ik} \delta = \underline{x}_{ik} d'_{i}$.
Otherwise, $u_{ik}^* < 0$, then $u'_{ik} = u_{ik}^* < \overline{x}_{ik} \delta = \overline{x}_{ik} d'_{i}$.
\item
Consider any $k \in \overline{\mathcal{K}}_{i} \setminus \overline{\mathcal{K}}_{i}^{-}$.
Then $u'_{ik} = u_{ik}^* \le 0 \le \overline{x}_{ik} \delta = \overline{x}_{ik} d'_{i}$.
\item
Consider any $j \in \mathcal{J} \setminus \{i\}$ and any $k \in \overline{\mathcal{K}}_{j}$.
Then $u'_{jk} = u_{jk}^* \left(1 - \delta\right) \le \overline{x}_{jk} d_{j}^* \left(1 - \delta\right) = \overline{x}_{jk} d'_{j}$.
Thus, $(d', d_{0}', u')$ satisfies~\eqref{eqn:spmnl5-upper-bound} for all $j \in \mathcal{J}$ and $k \in \overline{\mathcal{K}}_{j}$.
\end{itemize}
\item
Since $d_{i}^* = 0$, it follows from~\eqref{mp5-totalmarket} that $d_{0}^* + \sum_{j \in \mathcal{J} \setminus \{i\}} d_{j}^* = 1$, and thus
\[
d_{0}' + \sum_{j \in \mathcal{J}} d'_{j}
\ \ = \ \ d_{0}^* \left(1 - \delta\right) + \delta + \sum_{j \in \mathcal{J} \setminus \{i\}} d_{j}^* \left(1 - \delta\right)
\ \ = \ \ \left(d_{0}^* + \sum_{j \in \mathcal{J} \setminus \{i\}} d_{j}^*\right) \left(1 - \delta\right) + \delta
\ \ = \ \ 1 - \delta + \delta \ \ = \ \ 1
\]
Hence $(d', d_{0}', u')$ satisfies~\eqref{mp5-totalmarket}, and therefore $(d', d_{0}', u')$ is feasible for~\ref{eqn:static MNL5}.
\end{itemize}
Next we show that the objective value of~\ref{eqn:static MNL5} at $(d', d_{0}', u')$ is greater than the objective value at $(d^*, d_{0}^*, u^*)$.
The change in objective value from $(d^*, d_{0}^*, u^*)$ to $(d', d_{0}', u')$ is
\begin{align*}
& \sum_{j \in \mathcal{J}} \left(\sum_{k \in \mathcal{K}_{j}} \phi_{k} u'_{jk} - \psi_{j} d'_{j}\right) - \sum_{j \in \mathcal{J}} \left(\sum_{k \in \mathcal{K}_{j}} \phi_{k} u_{jk}^* - \psi_{j} d_{j}^*\right) \\
= \ \ & \phi_{l} \left[u_{il}^* \; - \sum_{\left\{k \in \underline{\mathcal{K}}_{i}^{+} \setminus \{l\} \, : \, u_{ik}^* = 0\right\}} \underline{x}_{ik} \delta \; - \sum_{\big\{k \in \overline{\mathcal{K}}_{i}^{-} \, : \, u_{ik}^* = 0\big\}} \overline{x}_{ik} \delta \; - \; \delta \ln\left(\frac{\delta}{d_{0}^* \left(1 - \delta\right)}\right)\right] \\
& + \sum_{\left\{k \in \underline{\mathcal{K}}_{i}^{+} \setminus \{l\} \, : \, u_{ik}^* = 0\right\}} \phi_{k} \left[u_{ik}^* + \underline{x}_{ik} \delta\right] \; + \sum_{\big\{k \in \overline{\mathcal{K}}_{i}^{-} \, : \, u_{ik}^* = 0\big\}} \phi_{k} \left[u_{ik}^* + \overline{x}_{ik} \delta\right] \\
& + \sum_{k \in \mathcal{K}_{i} \setminus \big(\{l\} \cup \big\{k \in \underline{\mathcal{K}}_{i}^{+} \, : \, u_{ik}^* = 0\big\} \cup \big\{k \in \overline{\mathcal{K}}_{i}^{-} \, : \, u_{ik}^* = 0\big\}\big)} \phi_{k} u_{ik}^* \; - \; \psi_{i} \delta \\
& + \sum_{j \in \mathcal{J} \setminus \{i\}} \left(\sum_{k \in \mathcal{K}_{j}} \phi_{k} u_{jk}^* \left(1 - \delta\right) - \psi_{j} d_{j}^* \left(1 - \delta\right)\right) - \sum_{j \in \mathcal{J}} \left(\sum_{k \in \mathcal{K}_{j}} \phi_{k} u_{jk}^* - \psi_{j} d_{j}^*\right) \\
= \ \ & \int_{0}^{\delta} g'(s) \diff s
\end{align*}
where $g : (0,1) \mapsto \mathbb{R}$ is given by
\begin{align*}
g(s) \ \ \defi \ \ & \phi_{l} \left[u_{il}^* \; - \sum_{\left\{k \in \underline{\mathcal{K}}_{i}^{+} \setminus \{l\} \, : \, u_{ik}^* = 0\right\}} \underline{x}_{ik} s \; - \sum_{\big\{k \in \overline{\mathcal{K}}_{i}^{-} \, : \, u_{ik}^* = 0\big\}} \overline{x}_{ik} s \; - \; s \ln\left(\frac{s}{d_{0}^* \left(1 - s\right)}\right)\right] \\
& + \sum_{\left\{k \in \underline{\mathcal{K}}_{i}^{+} \setminus \{l\} \, : \, u_{ik}^* = 0\right\}} \phi_{k} \left[u_{ik}^* + \underline{x}_{ik} s\right] \; + \sum_{\big\{k \in \overline{\mathcal{K}}_{i}^{-} \, : \, u_{ik}^* = 0\big\}} \phi_{k} \left[u_{ik}^* + \overline{x}_{ik} s\right] \\
& + \sum_{k \in \mathcal{K}_{i} \setminus \big(\{l\} \cup \big\{k \in \underline{\mathcal{K}}_{i}^{+} \, : \, u_{ik}^* = 0\big\} \cup \big\{k \in \overline{\mathcal{K}}_{i}^{-} \, : \, u_{ik}^* = 0\big\}\big)} \phi_{k} u_{ik}^* \; - \; \psi_{i} s \\
& + \sum_{j \in \mathcal{J} \setminus \{i\}} \left(\sum_{k \in \mathcal{K}_{j}} \phi_{k} u_{jk}^* \left(1 - s\right) - \psi_{j} d_{j}^* \left(1 - s\right)\right)
\end{align*}
Note that
{\small
\begin{align*}
g'(s) & \ \ = \ \ \phi_{l} \left[- \sum_{\left\{k \in \underline{\mathcal{K}}_{i}^{+} \setminus \{l\} \, : \, u_{ik}^* = 0\right\}} \underline{x}_{ik} \; - \sum_{\big\{k \in \overline{\mathcal{K}}_{i}^{-} \, : \, u_{ik}^* = 0\big\}} \overline{x}_{ik} \; - \; \ln\left(\frac{s}{d_{0}^* \left(1 - s\right)}\right) \; - \; 1 \; - \; \frac{s}{1 - s}\right] \\
& + \sum_{\left\{k \in \underline{\mathcal{K}}_{i}^{+} \setminus \{l\} \, : \, u_{ik}^* = 0\right\}} \phi_{k} \underline{x}_{ik} \; + \sum_{\big\{k \in \overline{\mathcal{K}}_{i}^{-} \, : \, u_{ik}^* = 0\big\}} \phi_{k} \overline{x}_{ik} \; - \; \psi_{i} \; - \sum_{j \in \mathcal{J} \setminus \{i\}} \left(\sum_{k \in \mathcal{K}_{j}} \phi_{k} u_{jk}^* - \psi_{j} d_{j}^*\right)
\ \ \to \ \ \infty \ \ \mbox{ as } \ \ s \downarrow 0
\end{align*}
}%
Thus $g'(s) > 0$ for all $s > 0$ sufficiently small, and therefore for $\delta > 0$ sufficiently small the objective value of~\ref{eqn:static MNL5} at $(d', d_{0}', u')$ is greater than the objective value at $(d^*, d_{0}^*, u^*)$.
This contradiction establishes that $d_{j}^* > 0$ for all $j \in \mathcal{J} \setminus \overline{\mathcal{J}}$.
In conclusion, it has been shown that $(d^*, d_{0}^*, u^*)$ is feasible for and therefore also optimal for~\ref{eqn:static MNL4}.
\end{proof}

\vspace{5mm}

\newpage

%------------------------------------------------------------------------

\subsection{Proof of Lemma~\ref{lem:fpmnl2 optimal}}

{\color{purple}
\begin{manuallemma}{\ref{lem:fpmnl2 optimal}}
If \ref{eqn:fluid MNL2} is feasible, then \ref{eqn:fluid MNL2} has an optimal solution if and only if for each $j \in \mathcal{J}$, $t \in \{0, 1, \ldots, T\}$, $k_{1} \in \mathcal{K}_{jt} \setminus \overline{\mathcal{K}}_{jt}$, and $k_{2} \in \mathcal{K}_{jt} \setminus \underline{\mathcal{K}}_{jt}$, it holds that $\phi_{k_{1},t} \le \phi_{k_{2},t}$.
\end{manuallemma}
}

\begin{proof}
Suppose that there is a $j \in \mathcal{J}$, a $t \in \{0, 1, \ldots, T\}$, a $k_{1} \in \mathcal{K}_{jt} \setminus \overline{\mathcal{K}}_{jt}$, and a $k_{2} \in \mathcal{K}_{jt} \setminus \underline{\mathcal{K}}_{jt}$, such that $\phi_{k_{1},t} > \phi_{k_{2},t}$.
Consider any feasible solution $(d^{0}, d^{0}_{0}, u^{0})$ and a sequence $\{(d^{n}, d^{n}_{0}, u^{n})\}_{n=0}^{\infty}$ of feasible solutions with $u^{n}_{j,k_{1},t} = u^{0}_{j,k_{1},t} + n$, $u^{n}_{j,k_{2},t} = u^{0}_{j,k_{2},t} - n$, $u^{n}_{jkt} = u^{0}_{jkt}$ for all $k \notin \{k_{1}, k_{2}\}$, $u^{n}_{j'kt'} = u^{0}_{j'kt'}$ for all $j' \neq j$ or $t' \neq t$, and $k \in \mathcal{K}_{j't'}$, $d^{n}_{j't'} = d^{0}_{j't'}$ for all $j' \in \mathcal{J}$ and $t' \in \{0, 1, \ldots, T\}$, and $d^{n}_{0t'} = d^{0}_{0t'}$ for all $t' \in \{0, 1, \ldots, T\}$.
Then $(d^{n}, d^{n}_{0}, u^{n})$ is feasible for \ref{eqn:fluid MNL2} for all~$n$, and the objective value of $(d^{n+1}, d^{n+1}_{0}, u^{n+1})$ exceeds the objective value of $(d^{n}, d^{n}_{0}, u^{n})$ by $\phi_{k_{1},t} - \phi_{k_{2},t}$ for all~$n$, and thus the objective value of $(d^{n}, d^{n}_{0}, u^{n})$ increases without bound as $n \to \infty$.

Suppose that, for each $j \in \mathcal{J}$, $t \in \{0, 1, \ldots, T\}$, and $k_{1} \in \mathcal{K}_{jt} \setminus \overline{\mathcal{K}}_{jt}$, $k_{2} \in \mathcal{K}_{jt} \setminus \underline{\mathcal{K}}_{jt}$, it holds that $\phi_{k_{1},t} \le \phi_{k_{2},t}$.
Let problem~\textsf{$\mathsf{FP_3^{MNL}}$} be the same as problem~\ref{eqn:fluid MNL2} with the resource constraint~\eqref{eqn:fluid MNL2 resource} removed.
Thus, \textsf{$\mathsf{FP_3^{MNL}}$} is a collection of separate problems, one problem for each $t \in \{0, 1, \ldots, T\}$, similar to~\ref{eqn:static MNL5}.
By the proof of Lemma~\ref{lem:spmnl5 optimal}, \textsf{$\mathsf{FP_3^{MNL}}$} is feasible, the dual of \textsf{$\mathsf{FP_3^{MNL}}$} is bounded, and has a feasible solution in the interior of the dual cones.

Since \ref{eqn:fluid MNL2} is feasible, the dual of~\ref{eqn:fluid MNL2} is bounded.
Also, since \ref{eqn:fluid MNL2} is the same as \textsf{$\mathsf{FP_3^{MNL}}$} with additional linear inequality constraints, the dual of~\ref{eqn:fluid MNL2} has additional signed variables beyond the dual of \textsf{$\mathsf{FP_3^{MNL}}$}.
Then a feasible solution of the dual of \textsf{$\mathsf{FP_3^{MNL}}$} in the interior of the dual cones, combined with these additional signed variables set to $0$, gives a feasible solution of the dual of \ref{eqn:fluid MNL2} in the interior of the dual cones.
Then it follows from the conic duality theorem (Theorem~1.4.2 in \cite{ben2001lectures}) that \ref{eqn:fluid MNL2} has an optimal solution.
\end{proof}

\vspace{5mm}

\newpage
%------------------------------------------------------------------------

\subsection{Proof of Lemma~\ref{lem:fpmnl1=fpmnl2}}

{\color{purple}
\begin{manuallemma}{\ref{lem:fpmnl1=fpmnl2}}
Suppose that \ref{eqn:fluid MNL2} has a feasible solution $(d', d_{0}', u')$ with $d' > 0$.
Then, any optimal solution $(d^*, d_{0}^*, u^*)$ for \ref{eqn:fluid MNL2} satisfies $d^* > 0$, $d_{0}^* > 0$.
\end{manuallemma}
}

\begin{proof}
First, suppose that $d_{0t}^* = 0$ for some $t \in \{0, 1, \ldots, T\}$.
Then it follows from~\eqref{eqn:fluid MNL2 cone} that $d_{jt}^* = 0$ for all $j \in \mathcal{J}$, which violates~\eqref{eqn:fluid MNL2 unit}.
Thus $d_{0t}^* > 0$ for all $t \in \{0, 1, \ldots, T\}$.

Next, suppose that $d_{jt}^* = 0$ for some $t \in \{0, 1, \ldots, T\}$ and some $j \in \overline{\mathcal{J}}_{t}$.
Then it follows from~\eqref{eqn:fluid MNL2 cone dummy} that $d_{0t}^* = 0$.
It was shown above that $d_{0t}^* > 0$ for all $t \in \{0, 1, \ldots, T\}$, and thus $d_{jt}^* > 0$ for all $t \in \{0, 1, \ldots, T\}$ and $j \in \overline{\mathcal{J}}_{t}$.

Next, suppose that $d_{is}^* = 0$ for some $s \in \{0, 1, \ldots, T\}$ and some $i \in \mathcal{J} \setminus \overline{\mathcal{J}}_{s}$.
Note that it follows from $\phi_{ks} > 0$ for all $k \in \mathcal{K}_{is}$, $i \in \mathcal{J} \setminus \overline{\mathcal{J}}_{s}$, and~\eqref{eqn:fluid MNL2 cone}, that $\sum_{k \in \mathcal{K}_{is}} u_{iks}^* = 0$.
Let $(d', d_{0}', u')$ be a feasible solution of \ref{eqn:fluid MNL2} with $d' > 0$.
By the same argument used above, it holds that $d_{0}' > 0$.
Consider any $\delta \in (0,1)$.
Since the feasible set of \ref{eqn:fluid MNL2} is convex, it follows that $(1 - \delta) (d^*, d_{0}^*, u^*) + \delta (d', d_{0}', u')$ is feasible for \ref{eqn:fluid MNL2}.
Consider the feasible solution $(\hat{d}, \hat{d}_{0}, \hat{u})$ for \ref{eqn:fluid MNL2} constructed to be the same as $(1 - \delta) (d^*, d_{0}^*, u^*) + \delta (d', d_{0}', u')$, except for a single variable $u_{ils}$, as follows:
Since $i \in \mathcal{J} \setminus \overline{\mathcal{J}}_{s}$, it holds that $\mathcal{K}_{is} \setminus \overline{\mathcal{K}}_{is} \neq \varnothing$.
Choose any $l \in \mathcal{K}_{is} \setminus \overline{\mathcal{K}}_{is}$.
Let $\hat{d}_{jt} = (1 - \delta) d_{jt}^* + \delta d_{jt}'$ for all $j \in \mathcal{J}$ and $t \in \{0,1,\ldots,T\}$, $\hat{d}_{0t} = (1 - \delta) d_{0t}^* + \delta d_{0t}' > 0$ for all $t \in \{0,1,\ldots,T\}$, and $\hat{u}_{jkt} = (1 - \delta) u_{jkt}^* + \delta u_{jkt}'$ for all $j \in \mathcal{J}$, $t \in \{0,1,\ldots,T\}$, and $k \in \mathcal{K}_{jt}$ such that $(j,k,t) \neq (i,l,s)$.
Note that $\hat{d}_{is} = (1 - \delta) d_{is}^* + \delta d_{is}' = \delta d_{is}' > 0$.
Finally, let
\begin{align*}
\hat{u}_{ils} \ \ = \ \ \hat{d}_{is} \ln\left(\frac{\hat{d}_{0s}}{\hat{d}_{is}}\right) - \sum_{k \in \mathcal{K}_{is} \setminus \{l\}} \hat{u}_{iks}
\ \ = \ \ \delta d_{is}' \ln\left(\frac{(1 - \delta) d_{0s}^* + \delta d_{0s}'}{\delta d_{is}'}\right) - \sum_{k \in \mathcal{K}_{is} \setminus \{l\}} \big[(1 - \delta) u_{iks}^* + \delta u_{iks}'\big]
\end{align*}
Note that, by choice of $\hat{u}_{ils}$, \eqref{eqn:fluid MNL2 cone} is satisfied for $i \in \mathcal{J}$ and $s \in \{0, 1, \ldots, T\}$.
If $l \in \underline{\mathcal{K}}_{is}$, then note that, since $(1 - \delta) (d^*, d_{0}^*, u^*) + \delta (d', d_{0}', u')$ is feasible for \ref{eqn:fluid MNL2}, it follows from~\eqref{eqn:fluid MNL2 lower} and~\eqref{eqn:fluid MNL2 cone} that
\begin{align*}
\underline{x}_{ils} \hat{d}_{is} \ \ & = \ \ \underline{x}_{ils} \big[(1 - \delta) d_{is}^* + \delta d_{is}'\big]
\ \ \le \ \ (1 - \delta) u_{ils}^* + \delta u_{ils}' \\
& \le \ \ \big[(1 - \delta) d_{is}^* + \delta d_{is}'\big] \ln\left(\frac{(1 - \delta) d_{0s}^* + \delta d_{0s}'}{(1 - \delta) d_{is}^* + \delta d_{is}'}\right) - \sum_{k \in \mathcal{K}_{is} \setminus \{l\}} \big[(1 - \delta) u_{iks}^* + \delta u_{iks}'\big] \\
& = \ \ \hat{d}_{is} \ln\left(\frac{\hat{d}_{0s}}{\hat{d}_{is}}\right) - \sum_{k \in \mathcal{K}_{is} \setminus \{l\}} \hat{u}_{iks}
\ \ = \ \ \hat{u}_{ils}
\end{align*}
and thus $\hat{u}_{ils}$ satisfies~\eqref{eqn:fluid MNL2 lower}.
In addition, since $l \in \mathcal{K}_{is} \setminus \overline{\mathcal{K}}_{is}$, the variable $u_{ils}$ is not upper bounded.
Hence $(\hat{d}, \hat{d}_{0}, \hat{u})$ is feasible for \ref{eqn:fluid MNL2}.

The difference between the objective values of~\ref{eqn:fluid MNL2} at $(\hat{d}, \hat{d}_{0}, \hat{u})$ and $(d^*, d_{0}^*, u^*)$ is
\begin{align*}
& \sum_{t = 0}^{T} \lambda_{t} \sum_{j \in \mathcal{J}} \left(\sum_{k \in \mathcal{K}_{jt}} \phi_{kt} \hat{u}_{jkt} - \psi_{jt} \hat{d}_{jt}\right)
- \sum_{t = 0}^{T} \lambda_{t} \sum_{j \in \mathcal{J}} \left(\sum_{k \in \mathcal{K}_{jt}} \phi_{kt} u^*_{jkt} - \psi_{jt} d^*_{jt}\right) \\
& = \ \ \sum_{t = 0}^{T} \lambda_{t} \sum_{j \in \mathcal{J}} \left(\sum_{k \in \mathcal{K}_{jt}} \phi_{kt} \big[(1 - \delta) u_{jkt}^* + \delta u_{jkt}'\big] - \psi_{jt} \big[(1 - \delta) d_{jt}^* + \delta d_{jt}'\big]\right) \\
& \qquad + \lambda_{s} \phi_{ls} \left\{\delta d_{is}' \ln\left(\frac{(1 - \delta) d_{0s}^* + \delta d_{0s}'}{\delta d_{is}'}\right) - \sum_{k \in \mathcal{K}_{is} \setminus \{l\}} \big[(1 - \delta) u_{iks}^* + \delta u_{iks}'\big] - \big[(1 - \delta) u_{ils}^* + \delta u_{ils}'\big]\right\} \\
& \qquad - \sum_{t = 0}^{T} \lambda_{t} \sum_{j \in \mathcal{J}} \left(\sum_{k \in \mathcal{K}_{jt}} \phi_{kt} u^*_{jkt} - \psi_{jt} d^*_{jt}\right) \\
& = \ \ \sum_{t = 0}^{T} \lambda_{t} \sum_{j \in \mathcal{J}} \left(\sum_{k \in \mathcal{K}_{jt}} \phi_{kt} u_{jkt}^* - \psi_{jt} d_{jt}^*\right) + \delta \sum_{t = 0}^{T} \lambda_{t} \sum_{j \in \mathcal{J}} \left(\sum_{k \in \mathcal{K}_{jt}} \phi_{kt} \big[u_{jkt}' - u_{jkt}^*\big] - \psi_{jt} \big[d_{jt}' - d_{jt}^*\big]\right) \\
& \qquad + \lambda_{s} \phi_{ls} \left\{\delta d_{is}' \ln\left(\frac{(1 - \delta) d_{0s}^* + \delta d_{0s}'}{\delta d_{is}'}\right) - \sum_{k \in \mathcal{K}_{is}} u_{iks}^* - \delta \sum_{k \in \mathcal{K}_{is}} \big[u_{iks}' - u_{iks}^*\big]\right\} \\
& \qquad - \sum_{t = 0}^{T} \lambda_{t} \sum_{j \in \mathcal{J}} \left(\sum_{k \in \mathcal{K}_{jt}} \phi_{kt} u^*_{jkt} - \psi_{jt} d^*_{jt}\right) \\
& = \ \ \delta \sum_{t = 0}^{T} \lambda_{t} \sum_{j \in \mathcal{J}} \left(\sum_{k \in \mathcal{K}_{jt}} \phi_{kt} \big[u_{jkt}' - u_{jkt}^*\big] - \psi_{jt} \big[d_{jt}' - d_{jt}^*\big]\right) \\
& \qquad + \lambda_{s} \phi_{ls} \left\{\delta d_{is}' \ln\left(\frac{(1 - \delta) d_{0s}^* + \delta d_{0s}'}{\delta d_{is}'}\right) - \delta \sum_{k \in \mathcal{K}_{is}} \big[u_{iks}' - u_{iks}^*\big]\right\} \\
& = \ \ \int_{0}^{\delta} h'(v) \diff v
\end{align*}
where the third equality follows from the result that $\sum_{k \in \mathcal{K}_{is}} u_{iks}^* = 0$, and $h : (0,1) \mapsto \mathbb{R}$ is given by
\begin{align*}
h(v) \ \ = \ \ & v \sum_{t = 0}^{T} \lambda_{t} \sum_{j \in \mathcal{J}} \left(\sum_{k \in \mathcal{K}_{jt}} \phi_{kt} \big[u_{jkt}' - u_{jkt}^*\big] - \psi_{jt} \big[d_{jt}' - d_{jt}^*\big]\right) \\
& + \lambda_{s} \phi_{ls} \left\{v d_{is}' \ln\left(\frac{(1 - v) d_{0s}^* + v d_{0s}'}{v d_{is}'}\right) - v \sum_{k \in \mathcal{K}_{is}} \big[u_{iks}' - u_{iks}^*\big]\right\}
\end{align*}
Note that
{\small
\begin{align*}
h'(v) \ \ = \ \ & \sum_{t = 0}^{T} \lambda_{t} \sum_{j \in \mathcal{J}} \left(\sum_{k \in \mathcal{K}_{jt}} \phi_{kt} \big[u_{jkt}' - u_{jkt}^*\big] - \psi_{jt} \big[d_{jt}' - d_{jt}^*\big]\right) \\
& + \lambda_{s} \phi_{ls} \left\{d_{is}' \ln\left(\frac{(1 - v) d_{0s}^* + v d_{0s}'}{v d_{is}'}\right) - \frac{d_{is}' d_{0s}^*}{(1 - v) d_{0s}^* + v d_{0s}'} - \sum_{k \in \mathcal{K}_{is}} \big[u_{iks}' - u_{iks}^*\big]\right\}
\ \ \to \ \ \infty \ \ \mbox{ as } \ \ v \downarrow 0
\end{align*}
}%
Thus $h'(v) > 0$ for all $v > 0$ sufficiently small, and therefore for $\delta > 0$ sufficiently small the objective value of~\ref{eqn:fluid MNL2} at $(\hat{d}, \hat{d}_{0}, \hat{u})$ is greater than the objective value at $(d^*, d_{0}^*, u^*)$.
This contradiction establishes that $d_{jt}^* > 0$ for all $t \in \{0, 1, \ldots, T\}$ and $j \in \mathcal{J} \setminus \overline{\mathcal{J}}_{t}$.
In conclusion, it has been shown that $d^* > 0$, $d_{0}^* > 0$ for any optimal solution $(d^*, d_{0}^*, u^*)$ of~\ref{eqn:fluid MNL2}.
\end{proof}

\vspace{5mm}

\newpage
%------------------------------------------------------------------------

\subsection{Proof of Theorem~\ref{thm:fpmc1=fpmc2}}

{\color{purple}
\begin{manualtheorem}{\ref{thm:fpmc1=fpmc2}}
\ref{eqn:fluid MC1} can be solved by solving \ref{eqn:fluid MC2} and taking into account the following:
\begin{enumerate}
\item
If \ref{eqn:fluid MC2} is infeasible, then \ref{eqn:fluid MC1} is infeasible.
\item
If \ref{eqn:fluid MC2} is feasible, then:
\begin{enumerate}
\item
\ref{eqn:fluid MC2} is bounded, and has an optimal solution. 
\item
\ref{eqn:fluid MC1} is bounded, and has an optimal solution. 
\item
Let $(v^*, d^*, u^*)$ be an optimal solution of \ref{eqn:fluid MC2}. Let $x_{jkt}^* = u_{jkt}^* / d_{jt}^*$ for every $j \in \mathcal{J}$, $t \in \{0,1,\ldots,T\}$, and $k \in \mathcal{K}_{jt}$ such that $d_{jt}^* > 0$, and let $x_{jkt}^* = \underline{x}_{jkt}^*$ otherwise. Then $(v^*, d^*, x^*)$ is an optimal solution of \ref{eqn:fluid MC1}. 
\end{enumerate}
\end{enumerate}
\end{manualtheorem}
}

\begin{proof} First, consider any $(v,d,u)$ that is feasible to \ref{eqn:fluid MC2}. Consider any $j \in \mathcal{J}$ and $t \in \{0,1,\ldots,T\}$ such that $d_{jt}>0$ (which implies $v_{jt}>0$). By constraint (\ref{eqn:fluid MC2 cone}), we know
\begin{align*}
d_{jt}\ln\left(\frac{v_{jt}}{d_{jt}}\right) 
		\ \ \geq \ \ \sum_{k \in \mathcal{K}_{jt}} u_{jkt}
\end{align*}
Suppose that
\begin{align*}
d_{jt}\ln\left(\frac{v_{jt}}{d_{jt}}\right) 
		\ \ > \ \ \sum_{k \in \mathcal{K}_{jt}} u_{jkt}
\end{align*}
If $j \in \overline{\mathcal{J}}_{t}$, then by constraint (\ref{eqn:fluid MC2 cone dummy}), we know
\begin{align*}
\sum_{k \in \mathcal{K}_{jt}} \overline{x}_{jkt} d_{jt} 
		\ \ \geq \ \ d_{jt} \ln\left(\frac{v_{jt}}{d_{jt}}\right)
		\ \ \geq \ \ \sum_{k \in \mathcal{K}_{jt}} u_{jkt}
\end{align*}
 otherwise, $\sum_{k \in \mathcal{K}_{jt}} u_{jkt}$ is not upper bounded by any constraint other than (\ref{eqn:fluid MC2 cone}). In either case, we can increase $u_{jkt}$ for some $k$ while let $v_{jt}$ and $d_{jt}$ remain unchanged, and get a feasible solution with a better objective value than $(v,d,u)$. Thus, an optimal solution for \ref{eqn:fluid MC2} satisfies 
 \begin{align*}
d_{jt}\ln\left(\frac{v_{jt}}{d_{jt}}\right) \ \ = \ \ \sum_{k \in \mathcal{K}_{jt}} u_{jkt}
\end{align*}
 for any $j \in \mathcal{J}$ and $t \in \{0,1,\ldots,T\}$ such that $d_{jt}>0$. 
 
Second, consider any $(v,d,x)$ that is feasible to \ref{eqn:fluid MC1}. Let $u_{jkt} = x_{jkt} d_{jt}$ for every $j \in \mathcal{J}$, $t \in \{0,1,\ldots,T\}$, and $k \in \mathcal{K}_{jt}$. Then the solution $(v,d,x)$ is feasible to \ref{eqn:fluid MC2}, with the same objective value as $(v,d,x)$ in \ref{eqn:fluid MC1}. On the other hand, Consider any $(v,d,u)$ that is feasible to \ref{eqn:fluid MC2}. Then (by the previous paragraph) there is a solution $(v,d,u')$ that satisfies 
 \begin{align*}
d_{jt}\ln\left(\frac{v_{jt}}{d_{jt}}\right) \ \ = \ \ \sum_{k \in \mathcal{K}_{jt}} u_{jkt}'
\end{align*}
 for any $j \in \mathcal{J}$ and $t \in \{0,1,\ldots,T\}$ such that $d_{jt}>0$, whose objective value is at least as good as $(v,d,u)$. Let $x_{jkt} = u_{jkt}' / d_{jt}$ for every $j \in \mathcal{J}$, $t \in \{0,1,\ldots,T\}$, and $k \in \mathcal{K}_{jt}$ such that $d_{jt} > 0$, and let $x_{jkt} = \underline{x}_{jkt}$ otherwise. Then $(v, d, x)$ is a feasible solution of \ref{eqn:fluid MC1}, whose objective value is at least as good as $(v,d,u)$ in \ref{eqn:fluid MC2}. 
 
 Thus, we know that: 
 \begin{enumerate}
\item
\ref{eqn:fluid MC2} is infeasible if and only if \ref{eqn:fluid MC1} is infeasible.
\item
\ref{eqn:fluid MC2} is feasible and unbounded if and only if \ref{eqn:fluid MC1} is feasible and unbounded.
\item
\ref{eqn:fluid MC2} is feasible and bounded if and only if \ref{eqn:fluid MC1} is feasible and bounded.
\item
\ref{eqn:fluid MC2} has an optimal solution if and only if \ref{eqn:fluid MC1} has an optimal solution.
\item
Let $(v^*, d^*, u^*)$ be an optimal solution of \ref{eqn:fluid MC2}. Then
\begin{align*}
d_{jt}^*\ln\left(\frac{v_{jt}^*}{d_{jt}^*}\right) = \sum_{k \in \mathcal{K}_{jt}} u_{jkt}^*
\end{align*}
 for any $j \in \mathcal{J}$ and $t \in \{0,1,\ldots,T\}$ such that $d_{jt}^*>0$. 
\item
Let $(v^*, d^*, u^*)$ be an optimal solution of \ref{eqn:fluid MC2}. Let $x_{jkt}^* = u_{jkt}^* / d_{jt}^*$ for every $j \in \mathcal{J}$, $t \in \{0,1,\ldots,T\}$, and $k \in \mathcal{K}_{jt}$ such that $d_{jt}^* > 0$, and let $x_{jkt}^* = \underline{x}_{jkt}^*$ otherwise. Then $(v^*, d^*, x^*)$ is an optimal solution of \ref{eqn:fluid MC1}. 
\end{enumerate}

It left to show that an optimal solution exists for \ref{eqn:fluid MC2} if \ref{eqn:fluid MC2} is feasible. Indeed. Let problem~\textsf{$\mathsf{FP_3^{MC}}$} be the same as problem~\ref{eqn:fluid MC2} with the resource constraints removed.
Thus, \textsf{$\mathsf{FP_3^{MC}}$} is a collection of separate problems, one problem for each $t \in \{0, 1, \ldots, T\}$. Each of the separate problem have the following formulation: 
\begin{subequations}
\begin{align}
\max_{v, \, d, \, u} \quad & \sum_{j \in \mathcal{J}} \left(\sum_{k \in \mathcal{K}_{j}} \phi_{k} u_{jk} - \psi_{j} d_{j}\right)
\tag{\textsf{$\mathsf{SP_1^{MC}}$}}
\label{eqn:static MC1} \\
\text{s.t.} \quad
& v_{j} \ \ = \ \ \theta_{j} + \sum_{i \in \mathcal{J}} \rho_{ij} (v_{i} - d_{i})
		& \forall \ j \in \mathcal{J} \\
& \left(v_{j}, d_{j}, \sum_{k \in \mathcal{K}_{j}} u_{jk}\right) \ \ \in \ \ \mathcal{K}_{\exp}
		& \forall \ j \in \mathcal{J}
		\label{eqn:static MC1 cone} \\
& \left(d_{j}, v_{j}, - \sum_{k \in \mathcal{K}_{j}} \overline{x}_{jk} v_{j}\right) 
		\ \ \in \ \ \mathcal{K}_{\exp}
		& \forall \ j \in \overline{\mathcal{J}}_{}
		\label{eqn:static MC1 cone dummy} \\
& u_{jk} \ \ \geq \ \ \underline{x}_{jk} d_{j}
		& \forall \ j \in \mathcal{J}, \ k \in \underline{\mathcal{K}}_{j}
		\label{eqn:static MC1 lower} \\
& u_{jk} \ \ \leq \ \ \overline{x}_{jk} d_{j}
		& \forall \ j \in \mathcal{J}, \ k \in \overline{\mathcal{K}}_{j}
		\label{eqn:static MC1 upper} 
\end{align}
\end{subequations}
and is feasible. The dual problem of \ref{eqn:static MC1} can be written as: 
{\small
\begin{subequations}
\begin{align}
\min_{\pi, \, \varpi, \, \overline{\nu}, \, \underline{\nu}, \, \eta} \quad & \sum_{j \in \mathcal{J}} \theta_{j} \eta_{j} 
\tag{\textsf{$\mathsf{SD_1^{MC}}$}} 
\label{eqn:static MC dual1} \\
\text{s.t.} \quad
& \pi_{j1} + \varpi_{j2} - \varpi_{j3} \sum_{k \in \overline{\mathcal{K}}_j} \overline{x}_{jk}
		- \eta_j + \sum_{i \in \mathcal{J}} \rho_{ji} \eta_{i} \ \ = \ \ 0
		& \forall \ j \in \overline{\mathcal{J}}
		\label{eqn:static MC dual1 balance1} \\
& \pi_{j1} - \eta_j + \sum_{i \in \mathcal{J}} \rho_{ji} \eta_{i} \ \ = \ \ 0
		& \forall \ j \in \mathcal{J} \setminus \overline{\mathcal{J}}
		\label{eqn:static MC dual1 balance2} \\
& \pi_{j2} + \varpi_{j1}
		+ \sum_{k \in \overline{\mathcal{K}}_j} \overline{\nu}_{jk} \overline{x}_{jk} 
		- \sum_{k \in \mathcal{K}_j} \underline{\nu}_{jk} \underline{x}_{jk} 
		- \sum_{i \in \mathcal{J}} \rho_{ji} \eta_{i} \ \ = \ \ \psi_{j}
		& \forall \ j \in \overline{\mathcal{J}} \\
& \pi_{j2} 
		+ \sum_{k \in \overline{\mathcal{K}}_j} \overline{\nu}_{jk} \overline{x}_{jk} 
		- \sum_{k \in \mathcal{K}_j} \underline{\nu}_{jk} \underline{x}_{jk} 
		- \sum_{i \in \mathcal{J}} \rho_{ji} \eta_{i} \ \ = \ \ \psi_{j}
		& \forall \ j \in \mathcal{J} \setminus \overline{\mathcal{J}} \\
& \pi_{j3} + \underline{\nu}_{jk} - \overline{\nu}_{jk} \ \ = \ \ - \phi_{k}
		& \forall \ j \in \mathcal{J}, \ k \in \overline{\mathcal{K}}_{j} \\
& \pi_{j3} + \underline{\nu}_{jk} \ \ = \ \ - \phi_{k}
		& \forall \ j \in \mathcal{J}, \ k \in \mathcal{K}_{j} \setminus \overline{\mathcal{K}}_{j} \\
& \pi_{j} \ \ \in \ \ \mathcal{K}_{\exp}^*
		& \forall \ j \in \mathcal{J} \\
& \varpi_{j} \ \ \in \ \ \mathcal{K}_{\exp}^*
		& \forall \ j \in \overline{\mathcal{J}} \\
& \overline{\nu}_{jk} \ \ \ge \ \ 0
		& \forall \ j \in \mathcal{J}, \ k \in \overline{\mathcal{K}}_{j}\\
& \underline{\nu}_{jk} \ \ \ge \ \ 0
		& \forall \ j \in \mathcal{J}, \ k \in \mathcal{K}_{j} 
\end{align}
\end{subequations}
}%
The refined conic duality theorem (see \cite{ben2001lectures}) states that \ref{eqn:static MC1} has an optimal solution if its dual problem \ref{eqn:static MC dual1} is bounded below and is ``strictly feasible'' (i.e. it has a feasible solution in the interior of the cones). Since \ref{eqn:static MC1} is feasible, we know that \ref{eqn:static MC dual1} is bounded below. Next we show that \ref{eqn:static MC dual1} has a feasible solution such that the conic variables are all in the interior of $\mathcal{K}_{\text{exp}}^*$. First, let
\begin{align*}
\overline{\nu}_{jk} \ \ &= \ \ 0 
		& \forall \ j \in \mathcal{J}, \ k \in \overline{\mathcal{K}}_{j} \\
\underline{\nu}_{jk} \ \ &= \ \ \max_{k}\{\phi_{k}\} - \phi_{k} \ \ \geq \ \ 0 
		& \forall \ j \in \mathcal{J}, \ k \in \mathcal{K}_{j} \\
\pi_{j1} \ \ &= \ \ 1 + \max_{k}\{\phi_{k}\} \cdot 
		\exp\left(-\frac{\psi_{j}}{\max_{k}\{\phi_{k}\}} - 1\right) \ \ > \ \ 0 
		& \forall \ j \in \mathcal{J} \\
\pi_{j2} \ \ &= \ \ \psi_{j} + \sum_{k \in \mathcal{K}_j} \underline{\nu}_{jk} \underline{x}_{jk} 
		+ \sum_{i \in \mathcal{J}} \rho_{ji} \eta_{i} \ \ \geq \ \ \psi_{j}
		& \forall \ j \in \mathcal{J} \\
\pi_{j3} \ \ &= \ \ -\max_{k}\{\phi_{k}\} \ \ < \ \ 0 
		& \forall \ j \in \mathcal{J} \\
\varpi \ \ &= \ \ 0 
\end{align*}
where $\eta = (1 - \rho)^{-1} \pi_{\cdot 1} \geq 0$ is a solution to constraints (\ref{eqn:static MC dual1 balance1}) and (\ref{eqn:static MC dual1 balance2}). (Recall that $I - \rho$ needs to be non-singular and have nonnegative entries for the MC model). Clearly, the solution we constructed is feasible. Meanwhile, for each $j \in \mathcal{J}$, we have
\begin{align*}
\pi_{j1} \ \ &= \ \ 1 + \max_{k}\{\phi_{k}\} \cdot 
		\exp\left(-\frac{\psi_{j}}{\max_{k}\{\phi_{k}\}} - 1\right) \\
&\geq \ \ 1 + \max_{k}\{\phi_{k}\} \cdot 
		\exp\left(-\frac{\pi_{j2}}{\max_{k}\{\phi_{k}\}} - 1\right) \\
&= \ \ 1 - \pi_{j3} \cdot 
		\exp\left(\frac{\pi_{j2}}{\pi_{j3}} - 1\right) \\
&> \ \ - \pi_{j3} \cdot 
		\exp\left(\frac{\pi_{j2}}{\pi_{j3}} - 1\right) 
\end{align*}
Thus, $\pi_{j} >_{\mathcal{K}_{\text{exp}}^*} 0$ for all $j \in \mathcal{J}$. Then consider each $j \in \overline{\mathcal{J}}$. From the current solution, we increase $\varpi_{j1}$ and decrease $\varpi_{j3}$ by some small enough $\epsilon > 0$. To stay feasible, we need to decrease $\pi_{j1}$ by $\epsilon \sum_{k \in \overline{\mathcal{K}}_j} \overline{x}_{jk}$ and decrease $\pi_{j2}$ by $\epsilon$. Let $(\pi', \varpi', \overline{\nu}', \underline{\nu}', \eta')$ be the updated solution, which is also feasible. For any $\epsilon > 0$, we have $\varpi_{j}' >_{\mathcal{K}_{\text{exp}}^*} 0$. Meanwhile, if $\epsilon$ is small enough, then $\pi_{j}'$ will stay in the interior of $\mathcal{K}_{\text{exp}}^*$ (as it is close enough to $\pi_{j}$). Thus, \ref{eqn:static MC1} has an optimal solution. 

Now, since \ref{eqn:fluid MC2} is feasible, the dual of~\ref{eqn:fluid MC2} is bounded.
Also, since \ref{eqn:fluid MC2} is the same as \textsf{$\mathsf{FP_3^{MC}}$} with additional linear inequality constraints, the dual of~\ref{eqn:fluid MC2} has additional signed variables beyond the dual of \textsf{$\mathsf{FP_3^{MC}}$}.
Then a feasible solution of the dual of \textsf{$\mathsf{FP_3^{MC}}$} in the interior of the dual cones, combined with these additional signed variables set to $0$, gives a feasible solution of the dual of \ref{eqn:fluid MC2} in the interior of the dual cones.
Then it follows from the conic duality theorem (Theorem~1.4.2 in \cite{ben2001lectures}) that \ref{eqn:fluid MC2} has an optimal solution.
\end{proof}

\newpage
%------------------------------------------------------------------------

\subsection{Proof of Theorem~\ref{thm:fpnl1=fpnl2}}

{\color{purple}
\begin{manualtheorem}{\ref{thm:fpnl1=fpnl2}}
If the following conditions hold: 

\noindent Condition (i): for any feasible solution $(d, p, p_{0}, v_{\leq}, u_{>})$ to \ref{eqn:fluid NL2} at which
\begin{align*}
	& \left(\frac{1}{\gamma_{it}} - 1\right)
		p_{it} \ln\left(\frac{p_{0t}}{p_{it}}\right)
		+ \sum_{j \in \mathcal{J}_i} d_{jt} \ln\left(\frac{p_{0t}}{d_{jt}}\right)
		\ \ = \ \
		\sum_{k \in \mathcal{K}_{it}'} v_{ikt}
		& \forall \ i \in \mathcal{I}_{\leq}, \ t = 0,1,\ldots,T 
\end{align*}
the system
\begin{align*}
& \left(\frac{1}{\gamma_{it}} - 1\right)
		d_{jt} \ln\left(\frac{p_{0t}}{p_{it}}\right)
		+ d_{jt} \ln\left(\frac{p_{0t}}{d_{jt}}\right)
		\ \ = \ \
		\sum_{k \in \mathcal{K}_{jt}} u_{jkt} 
		& \forall \ i \in \mathcal{I}_{\leq}, \ j \in \mathcal{J}_{i}, \ t = 0,1,\ldots,T \\
& v_{ikt} \ \ = \ \ \sum_{j \in \mathcal{J}_{ikt}'} u_{jkt}
		& \forall \ i \in \mathcal{I}_{\leq}, \ k \in \underline{\mathcal{K}}_{it}', \ t = 0,1,\ldots,T \\
& u_{jkt} \ \ \geq \ \ \underline{x}_{jkt} d_{jt}
& \forall \ j \in \mathcal{J}_{\leq}, \ k \in \underline{\mathcal{K}}_{jt}, \ t = 0,1,\ldots,T 
		\nonumber \\
& u_{jkt} \ \ \leq \ \ \overline{x}_{jkt} d_{jt}
& \forall \ j \in \mathcal{J}_{\leq}, \ k \in \overline{\mathcal{K}}_{jt}, \ t = 0,1,\ldots,T
		\nonumber 
\end{align*}
has a solution $u_{\leq}$. Here $\mathcal{J}_{ikt}' \defi \{j \in \mathcal{J}_{i} \ : \ k \in \mathcal{K}_{jt} \}$. 

\noindent Condition (ii): $\overline{\mathcal{K}}_{jt} \subsetneq \mathcal{K}_{jt}$ for every $j \in \mathcal{J}_{>}$, i.e. $\sum_{k \in \mathcal{K}_{jt}} u_{jkt}$ is not upper bounded for any $j \in \mathcal{J}_{>}$. 

\noindent then \ref{eqn:fluid NL1} can be solved by solving \ref{eqn:fluid NL2} and taking into account the following possibilities:
\begin{enumerate}
\item
If \ref{eqn:fluid NL2} is infeasible, then \ref{eqn:fluid NL1} is infeasible.
\item
If \ref{eqn:fluid NL2} is feasible and unbounded, then \ref{eqn:fluid NL1} is feasible and unbounded. 
\item
If \ref{eqn:fluid NL2} is feasible and has an optimal solution $(d^*, p^*, p_{0}^*, v_{\leq}^*, u_{>}^*)$, then $(d^*, p^*, x^*)$ with $x^*$ given by $x_{jkt}^* = u_{jkt}^* / d_{jt}^*$ for every $j \in \mathcal{J}$, $t \in \{0,1,\ldots,T\}$, and $k \in \mathcal{K}_{jt}$, is an optimal solution for \ref{eqn:fluid NL1}. (Here $u_{>}^*$ is part of the optimal solution $(d^*, p^*, p_{0}^*, v_{\leq}^*, u_{>}^*)$ to \ref{eqn:fluid NL2}, while $u_{\leq}^*$ is obtained by solving the system in Condition (i), using $(d^*, p^*, p_{0}^*, v_{\leq}^*, u_{>}^*)$ as an input.)
\end{enumerate}
\end{manualtheorem}
}

\begin{proof} 
First, consider any feasible solution $(d,p,x)$ to \ref{eqn:fluid NL1}. Let
\begin{align*}
p_{0} \ \ &= \ \ 1 - \sum_{i \in \mathcal{I}} p_{i} \\
v_{ikt} \ \ &= \ \ \sum_{j \in \mathcal{J}_{ikt}'} x_{jkt} d_{jt} 
		& \forall \ i \in \mathcal{I}_{\leq}, \ k \in \underline{\mathcal{K}}_{it}', \ t = 0,1,\ldots,T \\
u_{jkt} \ \ &= \ \ x_{jkt} d_{jt} 
		& \forall \ j \in \mathcal{J}_{>}, \ k \in \underline{\mathcal{K}}_{jt}, \ t = 0,1,\ldots,T 
\end{align*}
Then $(d, p, p_{0}, v_{\leq}, u_{>})$ is a feasible solution to \ref{eqn:fluid NL2} with the same objective value as $(d,p,x)$ in \ref{eqn:fluid NL1}. 
Second, consider any feasible solution $(d, p, p_{0}, v_{\leq}, u_{>})$ to \ref{eqn:fluid NL2}. 
Then we can find a solution $(d, p, p_{0}, v_{\leq}', u_{>}')$ (by increased some components in $v_{\leq}$ and $u_{>}$) whose objective value is at least as good as $(d, p, p_{0}, v_{\leq}, u_{>})$, such that 
\begin{align*}
& \left(\frac{1}{\gamma_{it}} - 1\right)
		p_{it} \ln\left(\frac{p_{0t}}{p_{it}}\right)
		+ \sum_{j \in \mathcal{J}_i} d_{jt} \ln\left(\frac{p_{0t}}{d_{jt}}\right)
		\ \ = \ \
		\sum_{k \in \mathcal{K}_{it}'} v_{ikt}'
		& \forall \ i \in \mathcal{I}_{\leq}, \ t = 0,1,\ldots,T  \\
& \left(1 - \frac{1}{\gamma_{it}}\right)
		d_{jt} \ln\left(\frac{p_{it}}{d_{jt}}\right)
		+ \left(\frac{1}{\gamma_{it}}\right)
		d_{jt} \ln\left(\frac{p_{0t}}{d_{jt}}\right)
		\ \ = \ \
		\sum_{k \in \mathcal{K}_{jt}} u_{jkt}'
		& \forall \ i \in \mathcal{I}_{>}, \ j \in \mathcal{J}_{i}, \ t = 0,1,\ldots,T 
\end{align*}
(The argument above can be proved the same way as in the proof of Theorem \ref{thm:fpmc1=fpmc2}.) Let $u_{\leq}'$ be a solution of the system in Condition (i). Let 
\begin{align*}
x_{jkt} \ \ &= \ \ 
\begin{cases} 
u_{jkt}' / d_{jt}  & \mbox{if } d_{jt} > 0 \\
\mbox{anything that satisfies the upper and lower bound (if any)} & \mbox{if } d_{jt} = 0
\end{cases}
\end{align*}
for each $j \in \mathcal{J}_{>}$, $k \in \underline{\mathcal{K}}_{jt}$ and $t = 0,1,\ldots,T$. The $(d,p,x)$ is a feasible solution to \ref{eqn:fluid NL1} with the same objective value as $(d, p, p_{0}, v_{\leq}', u_{>}')$ in \ref{eqn:fluid NL2}, which is at least as good as the objective value of $(d, p, p_{0}, v_{\leq}, u_{>})$. Thus, we know that: 
\begin{enumerate}
\item
If \ref{eqn:fluid NL2} is infeasible, then \ref{eqn:fluid NL1} is infeasible.
\item
If \ref{eqn:fluid NL2} is feasible and unbounded, then \ref{eqn:fluid NL1} is feasible and unbounded. 
\item
If \ref{eqn:fluid NL2} is feasible and has an optimal solution $(d^*, p^*, p_{0}^*, v_{\leq}^*, u_{>}^*)$, then $(d^*, p^*, x^*)$ with $x^*$ given by $x_{jkt}^* = u_{jkt}^* / d_{jt}^*$ for every $j \in \mathcal{J}$, $t \in \{0,1,\ldots,T\}$, and $k \in \mathcal{K}_{jt}$, is an optimal solution for \ref{eqn:fluid NL1}. (Here $u_{>}^*$ is part of the optimal solution $(d^*, p^*, p_{0}^*, v_{\leq}^*, u_{>}^*)$ to \ref{eqn:fluid NL2}, while $u_{\leq}^*$ is obtained by solving the system in Condition (i), using $(d^*, p^*, p_{0}^*, v_{\leq}^*, u_{>}^*)$ as an input.)
\end{enumerate}
\end{proof}

\newpage
%------------------------------------------------------------------------

\subsection{Proof of Theorem~\ref{thm:fpnl2=fpnl3}}

{\color{purple}
\begin{manualtheorem}{\ref{thm:fpnl2=fpnl3}}
Assume that Condition (i) in Theorem \ref{thm:fpnl1=fpnl2} holds.

\noindent Then, \ref{eqn:fluid NL2} can be solved by solving \ref{eqn:fluid NL3} and taking into account the following:
\begin{enumerate}
\item
If \ref{eqn:fluid NL3} is infeasible, then \ref{eqn:fluid NL2} is infeasible.
\item
If \ref{eqn:fluid NL3} is feasible, then the following possibilities hold:
\begin{enumerate}
\item
If \ref{eqn:fluid NL3} is unbounded (which happens if and only if for some $j \in \mathcal{J}$, $t \in \{0, 1, \ldots, T\}$, $k_{1} \in \mathcal{K}_{jt} \setminus \overline{\mathcal{K}}_{jt}$, and $k_{2} \in \mathcal{K}_{jt} \setminus \underline{\mathcal{K}}_{jt}$, it holds that $\phi_{k_{1},t} > \phi_{k_{2},t}$), then the following possibilities hold:
\begin{enumerate}
\item
If \ref{eqn:fluid NL3} does not have a feasible solution $(d', p', p_{0}', v_{\leq}', u_{>}', e', f_{\leq}', g_{>}', r_{\leq}', s_{>}')$ with $d' > 0$, then \ref{eqn:fluid NL2} is infeasible.
\item
If \ref{eqn:fluid NL3} has a feasible solution $(d', p', p_{0}', v_{\leq}', u_{>}', e', f_{\leq}', g_{>}', r_{\leq}', s_{>}')$ with $d' > 0$, then \ref{eqn:fluid NL2} is unbounded.
\end{enumerate}
\item
If \ref{eqn:fluid NL3} has an optimal solution $(d^*, p^*, p_{0}^*, v_{\leq}^*, u_{>}^*, e^*, f_{\leq}^*, g_{>}^*, r_{\leq}^*, s_{>}^*)$ (which happens if and only if for each $j \in \mathcal{J}$, $t \in \{0, 1, \ldots, T\}$, $k_{1} \in \mathcal{K}_{jt} \setminus \overline{\mathcal{K}}_{jt}$, and $k_{2} \in \mathcal{K}_{jt} \setminus \underline{\mathcal{K}}_{jt}$, it holds that $\phi_{k_{1},t} \le \phi_{k_{2},t}$), then the following possibilities hold:
\begin{enumerate}
\item
If $d^* > 0$, then $(d^*, p^*, p_{0}^*, v_{\leq}^*, u_{>}^*)$ is an optimal solution for \ref{eqn:fluid NL2}. 
\item
If $d_{jt}^* = 0$ for some $j \in \mathcal{J}$ and $t \in \{0,1,\ldots,T\}$, then \ref{eqn:fluid NL2} is infeasible.
\end{enumerate}
\end{enumerate}
\end{enumerate}

\end{manualtheorem}
}

\begin{proof} 
Same as Theorem \ref{thm:fpmnl1=fpmnl2}, Theorem \ref{thm:fpnl2=fpnl3} can be proved by showing two facts: 
\begin{enumerate}
\item
If \ref{eqn:fluid NL3} is feasible, then \ref{eqn:fluid NL3} has an optimal solution if and only if for each $j \in \mathcal{J}$, $t \in \{0, 1, \ldots, T\}$, $k_{1} \in \mathcal{K}_{jt} \setminus \overline{\mathcal{K}}_{jt}$, and $k_{2} \in \mathcal{K}_{jt} \setminus \underline{\mathcal{K}}_{jt}$, it holds that $\phi_{k_{1},t} \le \phi_{k_{2},t}$.
\item
Suppose that \ref{eqn:fluid NL3} has a feasible solution $(d', d_{0}', u')$ with $d' > 0$.
Then, any optimal solution $(d^*, p^*, p_{0}^*, v_{\leq}^*, u_{>}^*, e^*, f_{\leq}^*, g_{>}^*, r_{\leq}^*, s_{>}^*)$ for \ref{eqn:fluid NL3} satisfies $d^* > 0$, $d_{0}^* > 0$.
\end{enumerate}

The second result can be shown in the same way as in the proof of Lemma~\ref{lem:fpmnl1=fpmnl2}. 
%As a sketch, if some $d_{jt}^* = 0$, then we can find a point on the line connecting $(d^*, p^*, p_{0}^*, v_{\leq}^*, u_{>}^*, e^*, f_{\leq}^*, g_{>}^*, r_{\leq}^*, s_{>}^*)$ and $(d', p', p_{0}', v_{\leq}', u_{>}', e', f_{\leq}', g_{>}', r_{\leq}', s_{>}')$ with $d' > 0$, and then increase some component of $v_{\leq}$ or $u_{>}$ to get a better solution than $(d^*, p^*, p_{0}^*, v_{\leq}^*, u_{>}^*, e^*, f_{\leq}^*, g_{>}^*, r_{\leq}^*, s_{>}^*)$. 
%The key idea we use here is that the derivative of $a \ln(b/a)$ goes to infinity as $a \to 0$, whenever $b>0$. 
%Thus, any $d_{jt}^* = 0$ will give us a chance to improve the objective value with an infinity large derivative, whenever there is a feasible direction (e.g. a solution with $d>0$ exists). 
On the other hand, if there is a $j \in \mathcal{J}$, a $t \in \{0, 1, \ldots, T\}$, a $k_{1} \in \mathcal{K}_{jt} \setminus \overline{\mathcal{K}}_{jt}$, and a $k_{2} \in \mathcal{K}_{jt} \setminus \underline{\mathcal{K}}_{jt}$, such that $\phi_{k_{1},t} > \phi_{k_{2},t}$, then \ref{eqn:fluid NL3} is unbounded when feasible. This can be shown by constructing a sequence of feasible solutions with objective values going to infinity, same as in the proof of Lemma~\ref{lem:fpmnl2 optimal}. 

It left to show that if \ref{eqn:fluid NL3} is feasible, then \ref{eqn:fluid NL3} has an optimal solution when for each $j \in \mathcal{J}$, $t \in \{0, 1, \ldots, T\}$, $k_{1} \in \mathcal{K}_{jt} \setminus \overline{\mathcal{K}}_{jt}$, and $k_{2} \in \mathcal{K}_{jt} \setminus \underline{\mathcal{K}}_{jt}$, it holds that $\phi_{k_{1},t} \le \phi_{k_{2},t}$. Let problem~\textsf{$\mathsf{FP_4^{NL}}$} be the same as problem~\ref{eqn:fluid NL3} with the resource constraints removed.
Thus, \textsf{$\mathsf{FP_4^{NL}}$} is a collection of separate problems, one problem for each $t \in \{0, 1, \ldots, T\}$. Each of the separate problem have the following formulation: 
{\small
\begin{subequations}
\begin{align}
\max_{\substack{d, \ p, \ p_{0}, \ v_{\leq}, \ u_{>} \\ 
		e, \ f_{\leq}, \ g_{>}, \ r_{\leq}, \ s_{>}}} \quad
& \sum_{i \in \mathcal{I}_{\leq}} 
		\sum_{k \in \mathcal{K}_{i}'} \phi_{k} v_{ik} 
		 + \sum_{j \in \mathcal{J}_{>}} \sum_{k \in \mathcal{K}_{j}} \phi_{k} u_{jk} 
		 - \sum_{j \in \mathcal{J}} \psi_{j} d_{j}
\tag{\textsf{$\mathsf{SP_1^{NL}}$}}
\label{eqn:static NL1} \\
\text{s.t.} \quad
& \left(\frac{1}{\gamma_{i}} - 1\right) r_{i}
		+ \sum_{j \in \mathcal{J}_i} e_{j}
		\ \ = \ \
		\sum_{k \in \mathcal{K}_{i}'} v_{ik}
		& \forall \ i \in \mathcal{I}_{\leq} \ \ ((\sigma_{\leq})_{i}) \\
& \left(\frac{1}{\gamma_{i}} - 1\right) s_{i} + f_{j}
		\ \ = \ \
		- p_{0} \sum_{k \in \mathcal{K}_{j}} \overline{x}_{jk}
		& \forall \ i \in \mathcal{I}_{\leq}, \ j \in \overline{\mathcal{J}}_{i} \ \ ((\kappa_{\leq})_{j}) \\
& \left(1 - \frac{1}{\gamma_{i}}\right) g_{j}
		+ \left(\frac{1}{\gamma_{i}}\right) e_{j}
		\ \ = \ \
		\sum_{k \in \mathcal{K}_{j}} u_{jk}
		& \forall \ i \in \mathcal{I}_{>}, \ j \in \mathcal{J}_{i} \ \ ((\kappa_{>})_{j}) \\
& \left(p_{0}, d_{j}, e_{j}\right) \ \ \in \ \ \mathcal{K}_{\exp}
		& \forall \ i \in \mathcal{I}, \ j \in \mathcal{J}_{i} \ \ (\pi_{j}) \\
& \left(p_{i}, d_{j}, g_{j}\right) \ \ \in \ \ \mathcal{K}_{\exp}
		& \forall \ i \in \mathcal{I}_{>}, \ j \in \mathcal{J}_{i} \ \ ((\varpi_{>})_{j}) \\
& \left(d_{j}, p_{0}, f_{j}\right) \ \ \in \ \ \mathcal{K}_{\exp}
		& \forall \ i \in \mathcal{I}_{\leq}, \ j \in \overline{\mathcal{J}}_{i} 
		\ \ ((\varpi_{\leq})_{j}) \\
& \left(p_{i}, p_{0}, s_{i}\right) \ \ \in \ \ \mathcal{K}_{\exp}
		& \forall \ i \in \mathcal{I}_{\leq} \ \ ((\omega_{\leq})_{i}) \\
& \left(p_{0}, p_{i}, r_{i}\right) \ \ \in \ \ \mathcal{K}_{\exp}
		& \forall \ i \in \mathcal{I}_{\leq} \ \ ((\tau_{\leq})_{i}) \\
& p_{0} + \sum_{i \in \mathcal{I}} p_{i} \ \ = \ \ 1 
		& \ \ (\eta_{0}) \\
& p_{i} \ \ = \ \ \sum_{j \in \mathcal{J}_{i}} d_{j}
		& \forall \ i \in \mathcal{I} \ \ (\eta_{j}) \\
& v_{ik} \ \ \geq \ \ \sum_{j \in \mathcal{J}_{i}} \underline{x}_{jk} d_{j}
& \forall \ i \in \mathcal{I}_{\leq}, \ k \in \underline{\mathcal{K}}_{i}'
\label{eqn:static NL1 v lower} \ \ ((\underline{\mu}_{\leq})_{ik}) \\
& v_{ik} \ \ \leq \ \ \sum_{j \in \mathcal{J}_{i}} \overline{x}_{jk} d_{j}
& \forall \ i \in \mathcal{I}_{\leq}, \ k \in \overline{\mathcal{K}}_{i}'
\label{eqn:static NL1 v upper} \ \ ((\overline{\mu}_{\leq})_{ik}) \\
& u_{jk} \ \ \geq \ \ \underline{x}_{jk} d_{j}
& \forall \ j \in \mathcal{J}_{>}, \ k \in \underline{\mathcal{K}}_{j}
\label{eqn:static NL1 u lower} \ \ ((\underline{\nu}_{>})_{jk}) \\
& u_{jk} \ \ \leq \ \ \overline{x}_{jk} d_{j}
& \forall \ j \in \mathcal{J}_{>}, \ k \in \overline{\mathcal{K}}_{j} 
\label{eqn:static NL1 u upper} \ \ ((\overline{\nu}_{>})_{jk})
\end{align}
\end{subequations}
}%
and is feasible. The dual problem of \ref{eqn:static NL1} can be written as: 
{\small
\begin{subequations}
\begin{align}
\min_{} \quad & \eta_{0} 
\tag{\textsf{$\mathsf{SD_1^{NL}}$}} 
\label{eqn:static NL dual1} \\
\text{s.t.} \quad
& \sum_{i \in \mathcal{I}} \sum_{j \in \mathcal{J}_{i}} \pi_{j1} 
		+ \sum_{i \in \mathcal{I}_{\leq}} \sum_{j \in \overline{\mathcal{J}}_{i}} \varpi_{j2} 
		+ \sum_{i \in \mathcal{I}_{\leq}} \omega_{i2} \nonumber \\
		& + \sum_{i \in \mathcal{I}_{\leq}} \tau_{i1}
		+ \sum_{i \in \mathcal{I}_{\leq}} \sum_{j \in \overline{\mathcal{J}}_{i}} 
				\left( \sum_{k \in \mathcal{K}_{j}} \overline{x}_{jk} \right) \kappa_{j}
		\ \ = \ \ \eta_{0} 
		& \ \ (p_{0}) 
		\nonumber \\
& \sum_{j \in \mathcal{J}_{i}} \varpi_{j1} + \eta_{i}
		\ \ = \ \ \eta_{0} 
		& i \in \mathcal{I}_{>} \ \ (p_{i}) 
		\nonumber \\
& \omega_{i1} + \tau_{i2} + \eta_{i}
		\ \ = \ \ \eta_{0} 
		& i \in \mathcal{I}_{\leq} \ \ (p_{i}) 
		\nonumber \\
& \pi_{j2} + \varpi_{j2} - \eta_{i} 
		- \sum_{k \in \underline{\mathcal{K}}_{j}} \underline{\nu}_{jk} \underline{x}_{jk} 
		+ \sum_{k \in \overline{\mathcal{K}}_{j}} \overline{\nu}_{jk} \overline{x}_{jk}
		\ \ = \ \ \psi_{j}
		& i \in \mathcal{I}_{>}, \ j \in \mathcal{J}_{i} \ \ (d_{j}) 
		\nonumber \\
& \pi_{j2} + \varpi_{j1} - \eta_{i} 
		- \sum_{k \in \underline{\mathcal{K}}_{i}'} \underline{\mu}_{jk} \underline{x}_{jk} 
		+ \sum_{k \in \overline{\mathcal{K}}_{i}'} \overline{\mu}_{jk} \overline{x}_{jk}
		\ \ = \ \ \psi_{j}
		& i \in \mathcal{I}_{\leq}, \ j \in \overline{\mathcal{J}}_{i} \ \ (d_{j}) 
		\nonumber \\
& \pi_{j2} - \eta_{i} 
		- \sum_{k \in \underline{\mathcal{K}}_{i}'} \underline{\mu}_{jk} \underline{x}_{jk} 
		+ \sum_{k \in \overline{\mathcal{K}}_{i}'} \overline{\mu}_{jk} \overline{x}_{jk}
		\ \ = \ \ \psi_{j}
		& i \in \mathcal{I}_{\leq}, \ j \in \mathcal{J} \setminus \overline{\mathcal{J}}_{i} \ \ (d_{j}) 
		\nonumber \\
& -\kappa_{j} \ \ = \ \ -\phi_{k}
		& i \in \mathcal{I}_{>}, \ j \in \mathcal{J}_{i} , 
		\ k \in \mathcal{K}_{j} \setminus 
		(\underline{\mathcal{K}}_{j} \cap \overline{\mathcal{K}}_{j}) \ \ (u_{jk}) 
		\nonumber \\
& -\kappa_{j} + \underline{\nu}_{jk} \ \ = \ \ -\phi_{k}
		& i \in \mathcal{I}_{>}, \ j \in \mathcal{J}_{i} , 
		\ k \in \underline{\mathcal{K}}_{j} \setminus \overline{\mathcal{K}}_{j} \ \ (u_{jk}) 
		\nonumber \\
& -\kappa_{j} - \overline{\nu}_{jk} \ \ = \ \ -\phi_{k}
		& i \in \mathcal{I}_{>}, \ j \in \mathcal{J}_{i} , 
		\ k \in \overline{\mathcal{K}}_{j} \setminus \underline{\mathcal{K}}_{j} \ \ (u_{jk}) 
		\nonumber \\
& -\kappa_{j} + \underline{\nu}_{jk} - \overline{\nu}_{jk} \ \ = \ \ -\phi_{k}
		& i \in \mathcal{I}_{>}, \ j \in \mathcal{J}_{i} , 
		\ k \in \underline{\mathcal{K}}_{j} \cap \overline{\mathcal{K}}_{j} \ \ (u_{jk}) 
		\nonumber \\
& -\sigma_{i} \ \ = \ \ -\phi_{k}
		& i \in \mathcal{I}_{\leq}, 
		\ k \in \mathcal{K}_{i}' \setminus 
		(\underline{\mathcal{K}}_{i}' \cap \overline{\mathcal{K}}_{i}') \ \ (v_{ik}) 
		\nonumber \\
& -\sigma_{i} + \underline{\mu}_{ik} \ \ = \ \ -\phi_{k}
		& i \in \mathcal{I}_{\leq}, 
		\ k \in \underline{\mathcal{K}}_{i}' \setminus \overline{\mathcal{K}}_{i}' \ \ (v_{ik}) 
		\nonumber \\
& -\sigma_{i} - \overline{\mu}_{ik} \ \ = \ \ -\phi_{k}
		& i \in \mathcal{I}_{\leq}, 
		\ k \in \overline{\mathcal{K}}_{i}' \setminus \underline{\mathcal{K}}_{i}' \ \ (v_{ik}) 
		\nonumber \\
& -\sigma_{i} + \underline{\mu}_{ik} - \overline{\mu}_{ik} \ \ = \ \ -\phi_{k}
		& i \in \mathcal{I}_{\leq}, 
		\ k \in \underline{\mathcal{K}}_{i}' \cap \overline{\mathcal{K}}_{i}' \ \ (v_{ik}) 
		\nonumber \\
& \pi_{j3} + \left(\frac{1}{\gamma_{i}}\right) \kappa_{j} \ \ = \ \ 0
		& i \in \mathcal{I}_{>}, \ j \in \mathcal{J}_{i} \ \ (e_{j}) 
		\nonumber \\
& \pi_{j3} + \sigma_{i} \ \ = \ \ 0
		& i \in \mathcal{I}_{\leq}, \ j \in \mathcal{J}_{i} \ \ (e_{j}) 
		\nonumber \\
& \varpi_{j3} + \kappa_{j} \ \ = \ \ 0
		& i \in \mathcal{I}_{\leq}, \ j \in \overline{\mathcal{J}}_{i} \ \ (f_{j}) 
		\nonumber \\
& \varpi_{j3} + \left(1 - \frac{1}{\gamma_{i}}\right) \kappa_{j} \ \ = \ \ 0
		& i \in \mathcal{I}_{>}, \ j \in \mathcal{J}_{i} \ \ (g_{j}) 
		\nonumber \\
& \tau_{i3} + \left(\frac{1}{\gamma_{i}} - 1\right) \sigma_{i} \ \ = \ \ 0
		& i \in \mathcal{I}_{\leq} \ \ (r_{i}) 
		\nonumber \\
& \omega_{i3} + \left(\frac{1}{\gamma_{i}} - 1\right) 
		\sum_{j \in \overline{\mathcal{J}}_{i}} \kappa_{j} \ \ = \ \ 0
		& i \in \mathcal{I}_{\leq} \ \ (s_{i}) 
		\nonumber \\
& \pi_{j} \ \ \in \ \ \mathcal{K}_{\exp}^*
		& \forall \ i \in \mathcal{I}, \ j \in \mathcal{J}_{i} 
		\nonumber \\
& (\varpi_{>})_{j} \ \ \in \ \ \mathcal{K}_{\exp}^*
		& \forall \ i \in \mathcal{I}_{>}, \ j \in \mathcal{J}_{i} 
		\nonumber \\
& (\varpi_{\leq})_{j} \ \ \in \ \ \mathcal{K}_{\exp}^*
		& \forall \ i \in \mathcal{I}_{\leq}, \ j \in \overline{\mathcal{J}}_{i} 
		\nonumber \\
& (\omega_{\leq})_{i} \ \ \in \ \ \mathcal{K}_{\exp}^*
		& \forall \ i \in \mathcal{I}_{\leq} 
		\nonumber \\
& (\tau_{\leq})_{i} \ \ \in \ \ \mathcal{K}_{\exp}^*
		& \forall \ i \in \mathcal{I}_{\leq} 
		\nonumber \\
& (\underline{\mu}_{\leq})_{ik} \ \ \geq \ \ 0
		& \forall \ i \in \mathcal{I}_{\leq}, \ k \in \underline{\mathcal{K}}_{i}' 
		\nonumber \\
& (\overline{\mu}_{\leq})_{ik} \ \ \geq \ \ 0
		& \forall \ i \in \mathcal{I}_{\leq}, \ k \in \overline{\mathcal{K}}_{i}' 
		\nonumber \\
& (\underline{\nu}_{>})_{jk} \ \ \geq \ \ 0
		& \forall \ j \in \mathcal{J}_{>}, \ k \in \underline{\mathcal{K}}_{j} 
		\nonumber \\
& (\overline{\nu}_{>})_{jk} \ \ \geq \ \ 0
		& \forall \ j \in \mathcal{J}_{>}, \ k \in \overline{\mathcal{K}}_{j} 
		\nonumber 
\end{align}
\end{subequations}
}%

The refined conic duality theorem (see \cite{ben2001lectures}) states that \ref{eqn:static NL1} has an optimal solution if its dual problem \ref{eqn:static NL dual1} is bounded below and is ``strictly feasible'' (i.e. it has a feasible solution in the interior of the cones). Since \ref{eqn:static NL1} is feasible, we know that \ref{eqn:static NL dual1} is bounded below. Next we show that \ref{eqn:static NL dual1} has a feasible solution such that the conic variables are all in the interior of $\mathcal{K}_{\text{exp}}^*$. As a preparation, note that given Condition (i) in Theorem \ref{thm:fpnl1=fpnl2}, the following two arguments are equivalent (either one applies the other): 
\begin{enumerate}
\item
For any $j \in \mathcal{J}$, $\phi_{k_{1}} \le \phi_{k_{2}}$ for all $k_{1} \in \mathcal{K}_{j} \setminus \overline{\mathcal{K}}_{j}$, $k_{2} \in \mathcal{K}_{j} \setminus \underline{\mathcal{K}}_{j}$. 
\item
For any $i \in \mathcal{I}$, $\phi_{k_{1}} \le \phi_{k_{2}}$ for all $k_{1} \in \mathcal{K}_{i}' \setminus \overline{\mathcal{K}}_{i}'$, $k_{2} \in \mathcal{K}_{i}' \setminus \underline{\mathcal{K}}_{i}'$. 
\end{enumerate}

For each $i \in \mathcal{I}_{>}$, $j \in \mathcal{J}$, choose $\kappa_{j}$, $\underline{\nu}_{jk}$ for all $k \in \underline{\mathcal{K}}_{j}$, and $\overline{\nu}_{jk}$ for all $k \in \overline{\mathcal{K}}_{j}$, by considering the following cases: 
\begin{enumerate}
\item
If $\mathcal{K}_{j} \setminus \underline{\mathcal{K}}_{j} \neq \varnothing$, then choose $\kappa_{j} = - \min\left\{\phi_{k} \, : \, k \in \mathcal{K}_{j} \setminus \underline{\mathcal{K}}_{j}\right\}$.
Since $\phi_{k} > 0$ for all~$k$, it follows that $\kappa_{j} < 0$.
Then, for each $k \in \mathcal{K}_{j}$, consider the following four cases.
\begin{enumerate}
\item
Suppose $k \in \mathcal{K}_{j} \setminus \left(\underline{\mathcal{K}}_{j} \cup \overline{\mathcal{K}}_{j}\right)$.
Since $k \in \mathcal{K}_{j} \setminus \underline{\mathcal{K}}_{j}$, it follows that $\kappa_{j} \le \phi_{k}$.
Since $k \in \mathcal{K}_{j} \setminus \overline{\mathcal{K}}_{j}$ and $\phi_{k_{1}} \le \phi_{k_{2}}$ for all $k_{1} \in \mathcal{K}_{j} \setminus \overline{\mathcal{K}}_{j}$, $k_{2} \in \mathcal{K}_{j} \setminus \underline{\mathcal{K}}_{j}$, it follows that $\phi_{k} \le \kappa_{j}$.
Thus, $\phi_{k} = \kappa_{j}$.
\item
Suppose $k \in \underline{\mathcal{K}}_{j} \setminus \overline{\mathcal{K}}_{j}$.
Then choose $\underline{\nu}_{jk} = \kappa_{j} - \phi_{k}$.
Since $k \in \mathcal{K}_{j} \setminus \overline{\mathcal{K}}_{j}$ and $\phi_{k_{1}} \le \phi_{k_{2}}$ for all $k_{1} \in \mathcal{K}_{j} \setminus \overline{\mathcal{K}}_{j}$, $k_{2} \in \mathcal{K}_{j} \setminus \underline{\mathcal{K}}_{j}$, it follows that $\phi_{k} \le \kappa_{j}$.
Thus $\underline{\nu}_{jk} = \kappa_{j} - \phi_{k} \ge 0$.
\item
Suppose $k \in \overline{\mathcal{K}}_{j} \setminus \underline{\mathcal{K}}_{j}$.
Then choose $\overline{\nu}_{jk} = \kappa_{j} + \phi_{k}$.
Since $k \in \mathcal{K}_{j} \setminus \underline{\mathcal{K}}_{j}$, it follows that $\kappa_{j} \le \phi_{k}$.
Thus $\overline{\nu}_{jk} = \kappa_{j} + \phi_{k} \ge 0$.
\item
Suppose $k \in \underline{\mathcal{K}}_{j} \cap \overline{\mathcal{K}}_{j}$.
Then choose $\underline{\nu}_{jk} = \max\left\{0, \, \kappa_{j} - \phi_{k}\right\}$ and $\overline{\nu}_{jk} = \max\left\{0, \, \kappa_{j} + \phi_{k}\right\}$.
\end{enumerate}
\item
If $\mathcal{K}_{j} \setminus \underline{\mathcal{K}}_{j} = \varnothing$ (that is, $\mathcal{K}_{j} = \underline{\mathcal{K}}_{j}$), then choose $\kappa_{j} = - \max\left\{\phi_{k} \, : \, k \in \mathcal{K}_{j}\right\}$.
Since $\phi_{k} > 0$ for all~$k$, it follows that $\kappa_{j} < 0$.
Then, for each $k \in \mathcal{K}_{j}$, consider the following two cases (only these two cases are possible).
\begin{enumerate}
\item
Suppose $k \in \underline{\mathcal{K}}_{j} \setminus \overline{\mathcal{K}}_{j}$.
Then choose $\underline{\nu}_{jk} = \kappa_{j} - \phi_{k}$.
Since $\phi_{k} \le \kappa_{j}$, it follows that $\underline{\nu}_{jk} = \kappa_{j} - \phi_{k} \ge 0$. 
\item
Suppose $k \in \underline{\mathcal{K}}_{j} \cap \overline{\mathcal{K}}_{j}$.
Then choose $\underline{\nu}_{jk} = \max\left\{0, \, \kappa_{j} - \phi_{k}\right\}$ and $\overline{\nu}_{jk} = \max\left\{0, \, \kappa_{j} + \phi_{k}\right\}$.
\end{enumerate}
\end{enumerate}

For each $i \in \mathcal{I}_{\leq}$, choose $\sigma_{i}$, $\underline{\mu}_{ik}$ for all $k \in \underline{\mathcal{K}}_{i}'$, and $\overline{\mu}_{ik}$ for all $k \in \overline{\mathcal{K}}_{i}'$, by considering the following cases: 
\begin{enumerate}
\item
If $\mathcal{K}_{i}' \setminus \underline{\mathcal{K}}_{i}' \neq \varnothing$, then choose $\sigma_{i} = - \min\left\{\phi_{k} \, : \, k \in \mathcal{K}_{i}' \setminus \underline{\mathcal{K}}_{i}'\right\}$.
Since $\phi_{k} > 0$ for all~$k$, it follows that $\sigma_{i} < 0$.
Then, for each $k \in \mathcal{K}_{i}'$, consider the following four cases.
\begin{enumerate}
\item
Suppose $k \in \mathcal{K}_{i}' \setminus \left(\underline{\mathcal{K}}_{i}' \cup \overline{\mathcal{K}}_{i}'\right)$.
Since $k \in \mathcal{K}_{i}' \setminus \underline{\mathcal{K}}_{i}'$, it follows that $\sigma_{i} \le \phi_{k}$.
Since $k \in \mathcal{K}_{i}' \setminus \overline{\mathcal{K}}_{i}'$ and $\phi_{k_{1}} \le \phi_{k_{2}}$ for all $k_{1} \in \mathcal{K}_{i}' \setminus \overline{\mathcal{K}}_{i}'$, $k_{2} \in \mathcal{K}_{i}' \setminus \underline{\mathcal{K}}_{i}'$, it follows that $\phi_{k} \le \sigma_{i}$.
Thus, $\phi_{k} = \sigma_{i}$.
\item
Suppose $k \in \underline{\mathcal{K}}_{i}' \setminus \overline{\mathcal{K}}_{i}'$.
Then choose $\underline{\mu}_{ik} = \sigma_{i} - \phi_{k}$.
Since $k \in \mathcal{K}_{i}' \setminus \overline{\mathcal{K}}_{i}'$ and $\phi_{k_{1}} \le \phi_{k_{2}}$ for all $k_{1} \in \mathcal{K}_{i}' \setminus \overline{\mathcal{K}}_{i}'$, $k_{2} \in \mathcal{K}_{i}' \setminus \underline{\mathcal{K}}_{i}'$, it follows that $\phi_{k} \le \sigma_{i}$.
Thus $\underline{\mu}_{ik} = \sigma_{i} - \phi_{k} \ge 0$.
\item
Suppose $k \in \overline{\mathcal{K}}_{i}' \setminus \underline{\mathcal{K}}_{i}'$.
Then choose $\overline{\mu}_{ik} = \sigma_{i} + \phi_{k}$.
Since $k \in \mathcal{K}_{i}' \setminus \underline{\mathcal{K}}_{i}'$, it follows that $\sigma_{i} \le \phi_{k}$.
Thus $\overline{\mu}_{ik} = \sigma_{i} + \phi_{k} \ge 0$.
\item
Suppose $k \in \underline{\mathcal{K}}_{i}' \cap \overline{\mathcal{K}}_{i}'$.
Then choose $\underline{\mu}_{ik} = \max\left\{0, \, \sigma_{i} - \phi_{k}\right\}$ and $\overline{\mu}_{ik} = \max\left\{0, \, \sigma_{i} + \phi_{k}\right\}$.
\end{enumerate}
\item
If $\mathcal{K}_{i}' \setminus \underline{\mathcal{K}}_{i}' = \varnothing$ (that is, $\mathcal{K}_{i}' = \underline{\mathcal{K}}_{i}'$), then choose $\sigma_{i} = - \max\left\{\phi_{k} \, : \, k \in \mathcal{K}_{i}'\right\}$.
Since $\phi_{k} > 0$ for all~$k$, it follows that $\sigma_{i} < 0$.
Then, for each $k \in \mathcal{K}_{i}'$, consider the following two cases (only these two cases are possible).
\begin{enumerate}
\item
Suppose $k \in \underline{\mathcal{K}}_{i}' \setminus \overline{\mathcal{K}}_{i}'$.
Then choose $\underline{\mu}_{ik} = \sigma_{i} - \phi_{k}$.
Since $\phi_{k} \le \sigma_{i}$, it follows that $\underline{\mu}_{ik} = \sigma_{i} - \phi_{k} \ge 0$. 
\item
Suppose $k \in \underline{\mathcal{K}}_{i}' \cap \overline{\mathcal{K}}_{i}'$.
Then choose $\underline{\mu}_{ik} = \max\left\{0, \, \sigma_{i} - \phi_{k}\right\}$ and $\overline{\mu}_{ik} = \max\left\{0, \, \sigma_{i} + \phi_{k}\right\}$.
\end{enumerate}
\end{enumerate}

In addition, let
\begin{align*}
& \pi_{j3} \ \ = \ \ - \left(\frac{1}{\gamma_{i}}\right) \kappa_{j} \ \ < \ \ 0
		& i \in \mathcal{I}_{>}, \ j \in \mathcal{J}_{i}
		\nonumber \\
& \pi_{j3} \ \ = \ \ - \sigma_{i} \ \ < \ \ 0
		& i \in \mathcal{I}_{\leq}, \ j \in \mathcal{J}_{i}
		\nonumber \\
& \varpi_{j3} \ \ = \ \ - \varpi_{j1} \ \ = \ \ - \left(1 - \frac{1}{\gamma_{i}}\right) \kappa_{j} \ \ < \ \ 0
		& i \in \mathcal{I}_{>}, \ j \in \mathcal{J}_{i}
		\nonumber \\
& \varpi_{j3} \ \ = \ \ - \varpi_{j1} \ \ = \ \ - \kappa_{j} \ \ < \ \ 0
		& i \in \mathcal{I}_{\leq}, \ j \in \overline{\mathcal{J}}_{i}
		\nonumber \\
& \omega_{i3} \ \ = \ \ - \omega_{i1} \ \ = \ \ - \left(\frac{1}{\gamma_{i}} - 1\right)
		\sum_{j \in \overline{\mathcal{J}}_{i}} \kappa_{j} \ \ < \ \ 0
		& i \in \mathcal{I}_{\leq}
		\nonumber \\
& \tau_{i3} \ \ = \ \ - \left(\frac{1}{\gamma_{i}} - 1\right) \sigma_{i} \ \ < \ \ 0
		& i \in \mathcal{I}_{\leq}
		\nonumber 
\end{align*}
and let
\begin{align*}
& \pi_{j1} \ \ = \ \ C
		& i \in \mathcal{I}, \ j \in \mathcal{J}_{i}
		\nonumber \\
& \tau_{i1} \ \ = \ \ C
		& i \in \mathcal{I}_{\leq}
		\nonumber \\
& \varpi_{j2} \ \ = \ \ 0
		& i \in \mathcal{I}_{>}, \ j \in \mathcal{J}_{i}
		\nonumber \\
& \varpi_{j2} \ \ = \ \ 0
		& i \in \mathcal{I}_{\leq}, \ j \in \overline{\mathcal{J}}_{i}
		\nonumber \\
& \omega_{i2} \ \ = \ \ 0
		& i \in \mathcal{I}_{\leq}
		\nonumber \\
& \eta_{0} \ \ = \ \ 
		\sum_{i \in \mathcal{I}} \sum_{j \in \mathcal{J}_{i}} \pi_{j1} 
		+ \sum_{i \in \mathcal{I}_{\leq}} \tau_{i1}
		+ \sum_{i \in \mathcal{I}_{\leq}} \sum_{j \in \overline{\mathcal{J}}_{i}} 
				\left( \sum_{k \in \mathcal{K}_{j}} \overline{x}_{jk} \right) \kappa_{j}
		\nonumber \\
& \eta_{i} \ \ = \ \ \eta_{0} - \sum_{j \in \mathcal{J}_{i}} \varpi_{j1}
		& i \in \mathcal{I}_{>}
		\nonumber \\
& \tau_{i2} \ \ = \ \ \eta_{0} - \eta_{i} - \omega_{i1} \ \ = \ \ \sum_{j \in \mathcal{J}_{i}} \varpi_{j1} 
		& i \in \mathcal{I}_{\leq} \ \ (p_{i}) 
		\nonumber \\
& \pi_{j2} \ \ = \ \ \eta_{i} + \psi_{j} - \varpi_{j2} 
		+ \sum_{k \in \underline{\mathcal{K}}_{j}} \underline{\nu}_{jk} \underline{x}_{jk} 
		- \sum_{k \in \overline{\mathcal{K}}_{j}} \overline{\nu}_{jk} \overline{x}_{jk}
		& i \in \mathcal{I}_{>}, \ j \in \mathcal{J}_{i} 
		\nonumber \\
& \pi_{j2} \ \ = \ \ \eta_{i} + \psi_{j} - \varpi_{j1} 
		+ \sum_{k \in \underline{\mathcal{K}}_{i}'} \underline{\mu}_{jk} \underline{x}_{jk} 
		- \sum_{k \in \overline{\mathcal{K}}_{i}'} \overline{\mu}_{jk} \overline{x}_{jk}
		& i \in \mathcal{I}_{\leq}, \ j \in \overline{\mathcal{J}}_{i} 
		\nonumber \\
& \pi_{j2} \ \ = \ \ \eta_{i} + \psi_{j} 
		+ \sum_{k \in \underline{\mathcal{K}}_{i}'} \underline{\mu}_{jk} \underline{x}_{jk} 
		- \sum_{k \in \overline{\mathcal{K}}_{i}'} \overline{\mu}_{jk} \overline{x}_{jk}
		& i \in \mathcal{I}_{\leq}, \ j \in \mathcal{J} \setminus \overline{\mathcal{J}}_{i}
		\nonumber 
\end{align*}
where $C$ is a positive number. Thus, all equality and inequality constraints in \ref{eqn:static NL dual1} hold. Meanwhile, we can make $C$ large enough, such that
\begin{align*}
& \eta_{0} \ \ > \ \ 0 \\
& \eta_{i} \ \ > \ \ 0 
		& i \in \mathcal{I}_{>} \\
& \pi_{j1} \ \ \geq \ \ - \pi_{j3} \ \ > \ \ 0 \ \ ; \ \ \pi_{j2} \ \ \geq \ \ 0 
		& \forall \ i \in \mathcal{I}, \ j \in \mathcal{J}_{i} 
		\nonumber \\
& \varpi_{j1} \ \ = \ \ - \varpi_{j3} \ \ > \ \ 0 \ \ ; \ \ \varpi_{j2} \ \ = \ \ 0 
		& \forall \ i \in \mathcal{I}_{>}, \ j \in \mathcal{J}_{i} 
		\nonumber \\
& \varpi_{j1} \ \ = \ \ - \varpi_{j3} \ \ > \ \ 0 \ \ ; \ \ \varpi_{j2} \ \ = \ \ 0 
		& \forall \ i \in \mathcal{I}_{\leq}, \ j \in \overline{\mathcal{J}}_{i} 
		\nonumber \\
& \omega_{j1} \ \ = \ \ - \omega_{j3} \ \ > \ \ 0 \ \ ; \ \ \omega_{j2} \ \ = \ \ 0 
		& \forall \ i \in \mathcal{I}_{\leq} 
		\nonumber 
& \tau_{j1} \ \ > \ \ - \tau_{j3} \exp\left( \tau_{j2}/\tau_{j3} - 1 \right)
		& \forall \ i \in \mathcal{I}_{\leq} 
		\nonumber 
\end{align*}
which implies
\begin{align*}
& \pi_{j} \ \ >_{\mathcal{K}_{\exp}^*} \ \ 0
		& \forall \ i \in \mathcal{I}, \ j \in \mathcal{J}_{i} 
		\nonumber \\
& (\varpi_{>})_{j} \ \ >_{\mathcal{K}_{\exp}^*} \ \ 0
		& \forall \ i \in \mathcal{I}_{>}, \ j \in \mathcal{J}_{i} 
		\nonumber \\
& (\varpi_{\leq})_{j} \ \ >_{\mathcal{K}_{\exp}^*} \ \ 0
		& \forall \ i \in \mathcal{I}_{\leq}, \ j \in \overline{\mathcal{J}}_{i} 
		\nonumber \\
& (\omega_{\leq})_{i} \ \ >_{\mathcal{K}_{\exp}^*} \ \ 0
		& \forall \ i \in \mathcal{I}_{\leq} 
		\nonumber \\
& (\tau_{\leq})_{i} \ \ >_{\mathcal{K}_{\exp}^*} \ \ 0
		& \forall \ i \in \mathcal{I}_{\leq} 
		\nonumber 
\end{align*}
Thus, \ref{eqn:static NL dual1} has a feasible solution such that the conic variables are all in the interior of $\mathcal{K}_{\text{exp}}^*$. 

Now, since \ref{eqn:fluid NL2} is feasible, the dual of~\ref{eqn:fluid NL2} is bounded.
Also, since \ref{eqn:fluid NL2} is the same as \textsf{$\mathsf{FP_3^{NL}}$} with additional linear inequality constraints, the dual of~\ref{eqn:fluid NL2} has additional signed variables beyond the dual of \textsf{$\mathsf{FP_3^{NL}}$}.
Then a feasible solution of the dual of \textsf{$\mathsf{FP_3^{NL}}$} in the interior of the dual cones, combined with these additional signed variables set to $0$, gives a feasible solution of the dual of \ref{eqn:fluid NL2} in the interior of the dual cones.
Then it follows from the conic duality theorem (Theorem~1.4.2 in \cite{ben2001lectures}) that \ref{eqn:fluid NL2} has an optimal solution.
\end{proof}

 \end{appendix}

%------------------------------------------------------------------------
% Acknowledgments
%------------------------------------------------------------------------

% Acknowledgments here
% \ACKNOWLEDGMENT{}

\end{document}